\documentclass[smallextended,natbib,runningheads]{svjour3}
\journalname{Annals of the Institute of Statistical Mathematics}
\smartqed  

\usepackage{amsfonts}
\usepackage{amsmath}
\usepackage{amssymb}
\usepackage{natbib}
\usepackage{graphicx}
\usepackage[colorlinks=true,linkcolor=RawSienna,citecolor=RawSienna,urlcolor=RawSienna,bookmarksopen=true,pdfstartview=FitB]{hyperref}
\usepackage[left=1.5in,right=1.5in,top=1.5in,bottom=1.5in]{geometry}
\usepackage[dvipsnames]{xcolor}
\usepackage{float}
\usepackage{xr-hyper}%
\begin{document}

\title{Uniformly consistent proportion estimation for composite hypotheses via integral equations: ``the case of Gamma random variables"}

\titlerunning{Proportion estimation via integral equations: Gamma random variables}        

\author{Xiongzhi Chen}


\institute{Xiongzhi Chen \at
              Department of Mathematics and Statistics, Washington
State University, Pullman, WA 99164, USA. \\
              \email{xiongzhi.chen@wsu.edu}             \\
The online version of this article contains supplementary material.}

\date{Received: date / Revised: date}

\maketitle

\begin{abstract}
We consider estimating the proportion of random variables for two types of
composite null hypotheses: (i) the means of the random variables
belonging to a non-empty, bounded interval; (ii) the means of the
random variables belonging to an unbounded interval that is not the whole real
line. For each type of composite null hypotheses, uniformly consistent
estimators of the proportion of false null hypotheses are constructed
for random variables whose distributions are members of the Gamma family. Further, uniformly
consistent estimators of certain functions of a bounded null on the means
are provided for the random variables mentioned earlier.
These functions are continuous and of bounded variation.
The estimators are constructed via solutions to Lebesgue-Stieltjes integral
equations and harmonic analysis, do not rely on a concept of p-value, and have various applications.
\keywords{Composite null hypothesis \and harmonic analysis \and Lebesgue-Stieltjes integral equations \and proportion of false null hypotheses}
\end{abstract}

\section{Introduction}

\label{sec1intro}

The proportion of false null hypotheses (i.e., the ``signal proportion'' or ``alternative proportion'') and its dual, the proportion of true
null hypotheses (i.e., ``null proportion''), play important roles in statistical modelling and in control
and estimation of the false discovery rate (FDR, \cite{Benjamini:1995}) and the false nondiscovery
rate (FNR, \cite{Genovese:2002}). For example, simultaneously testing or
estimating the means of Gaussian random variables can be equivalently phrased
as simultaneously testing or estimating the regression parameters of a linear
model, and if we are interested in regression parameters whose values are
within a specified range, then the proportion of such parameters or some other
average norm of the regression parameters somehow characterizes how well an
estimator of the parameter vector or a testing procedure on the parameters can
perform; see, e.g., \cite{Abramovich:2006}. On the other hand, when
simultaneously testing whether the means or medians of many test statistics
are equal to a specific value or within a specific range, information on the proportion of such means or
medias helps better control the FDR and FNR and may improve the power of a
testing procedure, and estimators of the proportions are often used to
construct ``adaptive FDR procedures'' and ``adaptive FNR procedures'' that are
(asymptotically) conservative and often more powerful than their non-adaptive
counterparts; see, e.g., \cite{Storey:2004, Benjamini:2006, Sarkar:2006,Hu:2010,Nandi:2018,Durand:2019,Ignatiadis:2021} for adaptive FDR procedures for testing a point null hypothesis.
Further, these proportions are key components of the Bayesian two-component mixture model of \cite{Efron:2001}, its extensions
\cite{Ploner:2006,Cai:2009,Liu:2016}, and their corresponding decision rules
such as the ``local FDR'' \cite{Efron:2001}, ``q-value'' \cite{Storey:2003c} and the optimal discovery procedure \cite{Storey:2007}.
Without information on the proportions, these models and decision rules cannot be implemented in practice.
Finally, these proportions can be used to asymptotically ensure that a certain fraction of false null hypotheses are rejected among a set of rejections, as done by \cite{Jin:2008}.
However, neither the proportion of false null hypotheses nor its dual is
known, and it is very important to accurately estimate the proportions.

There are many estimators of the proportion of false (or true) null
hypotheses; see, e.g., those of
\cite{Storey:2002,Storey:2004,Genovese:2004,Kumar:2016,swanepoel1999,Meinshausen:2006,Jin:2008,Jin:2007,Jin:2010a,Cai:2010,Chen:2018a}.
However, they all estimate the proportion of parameters that are
equal to a specific value, i.e., they are almost all for the setting of a
point null hypothesis. Further, these estimators employ the Bayesian two-component
mixture model (which then forces test statistics whose parameters are tested
to be identically distributed marginally), require p-values to be identically
distributed under the alternative hypothesis or uniformly distributed under
the null hypothesis, require the distributions of test statistics to be
members of a location-shift family, or are usually inconsistent. A more
detailed discussion on these proportion estimators and others (including their underlying assumptions, advantages and disadvantages) is given by
\cite{Chen:2018a}.

In contrast to the setting of a point null hypothesis, there are many
situations where composite null hypotheses are tested and where the null
hypothesis assumes that parameters being tested are greater than or smaller than a fixed
value or are within a specific, finite range. For example, in gene expression
studies, a biologist may only be interested in genes whose expression
fold-changes exceed a threshold, or in genes that are down regulated or up
regulated, whereas in studies of drug side effects, only drugs whose side
effects are within a specific range can be recommended for ({urgent}) usage.
These are related to two types of commonly used composite null hypotheses:
``one-sided null'', i.e., a one-sided, unbounded interval to which the null
parameter belongs, and ``bounded null'', i.e., a non-empty, bounded interval to
which the null parameter belongs. Note that a bounded null is a much more
refined assessment on a parameter than a point null hypothesis. For multiple
testing such composite hypotheses, it is quite useful to construct adaptive
FDR and FNR procedures that employ consistent proportion estimators that are specifically
designed for these hypotheses. However, such proportion estimators are scarce and often not consistent.

\subsection{Review of major existing works}

We will abbreviate ``point null hypothesis'' as ``point null'', refer to an
estimator of the proportion of true (or false) null hypotheses as a ``null (or
alternative) proportion estimator'', and take a frequentist perspective. First, there does not seem to be a
proportion estimator that is designed for a one-side null. Storey's estimator
of \cite{Storey:2004} and the ``MR'' estimator of \cite{Meinshausen:2006} were motivated by and designed for
proportion estimation for a point null, are based on p-values, and require
p-value to be uniformly distributed under the null hypothesis. In addition,
the consistency of the ``MR'' estimator relies on the assumption that p-values
under the alternative hypothesis are identically distributed. Doubtlessly,
these two estimators can be applied to proportion estimation whenever p-values
can be defined. Specifically, they can be applied to a one-sided or bounded
null since the p-value for testing either of these nulls can be defined, e.g.,
by Definition 2.1 in Chapter of \cite{Dickhaus} as a generalization of the definition of
p-value in Section 3.3 of \cite{Lehmann:2005}. However, the p-value for
testing a composite null is not uniformly distributed under the null
hypothesis, and these two estimators are not theoretically guaranteed to
function properly in this setting. So, it seems natural to regard
Storey's estimator and the MR estimator as being inapplicable to proportion
estimation for a bounded null since they both require p-value to be uniformly distributed under the null hypothesis.

Second, Jin's estimator, which was constructed in Section 6 of \cite{Jin:2008}%
, can be applied to a special case of a bounded null (i.e., a symmetric
bounded null where the null parameter set is a symmetric interval around $0$)
but is not applicable to a one-sided null. Specifically, this estimator was
only constructed for Gaussian family via Fourier transform, and the
construction fails for a one-sided null since the Fourier transform for the
indicator function of an unbounded interval is undefined. We point out that
Jin's estimator can be used to estimate proportions induced by a suitable
function of the magnitude of the mean parameter of Gaussian random variables,
where the function has a compact support that is a symmetric interval around
$0$ and is even and continuous on its support; see Section 6 of
\cite{Jin:2008} for more details on this. However, the consistency of these
proportion estimators has not been proved in the sense of definition
(\ref{defConsistency}) to be introduced next.

Thirdly, there is another line of research that uses randomized p-values for
composite nulls; see, e.g., \cite{Dickhaus:2013,Hoang:2022a,Hoang:2022b}, who
for a composite null first define p-value under the
``least favorable parameter configurations
(LFCs)'' and then define randomized p-value by randomizing
p-value under LFCs. With such randomized p-values for a composite null,
proportion estimators that are based on p-values such as that of
\cite{Schweder:1982} can be used to estimate the proportion of true null
hypotheses, and this was done by \cite{Hoang:2022b}. Note that the estimator of \cite{Schweder:1982} implicitly requires p-values to be
uniformly distributed under the true null hypothesis. However, the definition
of randomized p-value for a composite null proposed by
\cite{Dickhaus:2013,Hoang:2022a} requires a specific set of assumptions to
yield a valid p-value or a valid p-value that can be practically computed.
Further, only for testing a one-side null on the means of Gaussian random variables
with unit variance or testing the sign of the means of Student t random variables
induced by Gaussian random variables, \cite{Dickhaus:2013,Hoang:2022a,Hoang:2022b} showed
that their randomized p-values are valid and provided methods on
how to compute them, whereas for a bounded null they provided no statistical
models for which their definition of randomized p-value leads to valid p-value
or its computation in practice. Regardless, to use randomized p-value proposed
by these authors, we need to check whether needed assumptions are satisfied
for each statistical model we employ and devise an algorithm to compute such
p-value in practice correspondingly.

\subsection{Main contributions and summary of results}

In this work, we employ the frequentist paradigm as a complement to the Bayesian one, take the parameter as the mean of a random
variable, and consider proportion estimation for multiple testing one-sided
nulls, bounded nulls, and a suitable functional of a bounded null on the
parameters, respectively. For each multiple testing scenario, we construct
uniformly consistent proportion estimators as solutions to Lebesgue-Stieltjes
integral equations and via Fourier analysis and Mellin transform. These new estimators are not based on p-values, and they do not require random variables, test
statistics or p-values to be identically distributed marginally, identically
distributed under the alternative hypothesis, or uniformly distributed under
the null hypothesis. Since both bounded and one-sided nulls have finite
boundary points, estimators of \cite{Chen:2018a} that are for a point null are
used to deal with proportions related solely to these points. However, we
provide new constructions of proportion estimators that utilize Dirichlet
integral to approximate the indicator function of a bounded or one-sided null.
These constructions are quite different from those in \cite{Jin:2008,Jin:2007,Jin:2010a,Cai:2010,Chen:2018a}, and
further demonstrate the power of integral equations and harmonic analysis in statistical estimation.

Our main contributions are summarized as follows:

\begin{itemize}
\item ``Construction I'': Construction of
proportion estimator for testing a bounded null on the means of random variables whose distributions are members of
the Gamma family, and under independence between these random variables the speed of
convergence and uniform consistency class of the estimator.

\item ``Construction II'': Construction of
proportion estimator for testing a one-sided null on the means of random variables whose distributions are members of
the Gamma family, and under independence between these random variables the speed of
convergence and uniform consistency class of the estimator.

\item Extension of Construction I to estimate the ``proportions''
induced by a function of the parameter that is continuous and of bounded variation (see \autoref{SecExtensions} for details), and under independence between these random variables the speeds
of convergence and uniform consistency classes of these estimators.
\end{itemize}

Our new estimators are the first uniformly consistent proportion estimators for a one-sided null or
a bounded null in the literature, and do not rely on a concept of p-value for testing a composite null.
Note here ``consistency'' is defined via the
``ratio'' (see the definition in (\ref{defConsistency}) to be introduced later) rather than the ``difference'' between an estimator and the alternative
proportion, in order to account for a diminishing alternative proportion
that converges to
$0$ as the number $m$ of hypotheses tends to $\infty$, since for such a proportion, an estimate that converges to zero in
probability as $m \to \infty$ is consistent in the classic definition of ``consistency'' that is based on such
difference, and consistency in terms of (\ref{defConsistency}) implies consistency in the classic sense.
Briefly speaking, each new estimator is constructed as follows:
\begin{itemize}
  \item First, we construct a ``nonrandom bridge'' $\varphi_{m} \left(t,\boldsymbol{\mu}\right)$ as a linear combination of Dirichlet integrals, for which these Dirichlet integrals are used to approximate the indicator function of a null parameter set, so that $\lim_{t \to \infty}\varphi_{m} \left(t,\boldsymbol{\mu}\right) = \pi_{1,m}$ for any $m$ and $m$-dimensional parameter vector $\boldsymbol{\mu}$, where $\pi_{1,m}$ is the alternative proportion.
  \item Second, we construct a new estimator $\hat{\varphi}_{m} \left(t,\mathbf{z}\right)$ as an unbiased estimator of $\varphi_{m} \left(t,\boldsymbol{\mu}\right)$, i.e.,
  $\varphi_{m} \left(t,\boldsymbol{\mu}\right) = \mathbb{E}\left(\hat{\varphi}_{m} \left(t,\mathbf{z}\right)\right)$, for any $m$ and $m$-dimensional data vector $\mathbf{z}$ whose parameter vector is $\boldsymbol{\mu}$.
\end{itemize}
Due to this construction, how fast $\hat{\varphi}_{m} \left(t_m,\mathbf{z}\right)$ is able to achieve (uniform) consistency in terms of (\ref{defConsistency}), i.e., how fast
$\left\vert \pi_{1,m}^{-1}\hat{\varphi}_{m} \left(t_m,\mathbf{z}\right) -1 \right\vert \to 0$ in probability for a sequence $\left\{t_m,m\ge1\right\}$ as $m \to \infty$
 is determined by how fast the slower of $\left\vert\pi_{1,m}^{-1}\varphi_{m} \left(t_m,\boldsymbol{\mu}\right) -1\right\vert \to 0$ and $\pi_{1,m}^{-1}\left \vert \hat{\varphi}_{m} \left(t_m,\mathbf{z}\right) - \varphi_{m} \left(t_m,\boldsymbol{\mu}\right) \right \vert \to 0$ in probability as $m \to \infty$.
 Here $t_m,m\ge1$ is referred to as the speed of convergence of the estimator $\hat{\varphi}_{m} \left(t_m,\mathbf{z}\right)$, and the faster $t_m \to \infty$ as $m \to \infty$, the faster
 the speed of convergence is.
For estimating the proportion of false null hypotheses under independence for
both bounded and one-sided nulls, the maximal speeds of convergence of our new
estimators are of the same order as those for the proportion estimators for a
point null in \cite{Chen:2018a} {for Gamma random variables}. However, due to the use of Dirichlet integral to approximate the indicator function of a
bounded or one-sided null, the sparsest alternative proportions the new
estimators can consistently estimate under independence are of smaller order
than their maximal speeds of convergence.

Such a difference between the sparsest alternative proportions that these estimators can consistently estimate under
a point null, bounded null or one-sided null respectively is the main consequence of the speed of convergence of $\varphi_{m} \left(t,\boldsymbol{\mu}\right)$  to $\pi_{1,m}$,
 and this speed is related to that of Dirichlet integral to an indicator function
for a one-sided or bounded null but is independent of the alternative
proportion for a point null. In other words, we have discovered the
universal phenomenon that, for an alternative proportion estimator that is
constructed via a solution to some Lebesgue-Stieltjes equation as an
approximator to the indicator of (a transform of) the parameter set under the
alternative hypothesis, the sparsest alternative proportion such an estimator
is able to consistently estimate can never be of larger order than the maximal speed
of convergence of the solution to its targeted indicator function.

As a by-product, we provide the speed of convergence of Dirichlet-type
integrals (see \autoref{lm:Dirichlet} {in the appendix} and \autoref{ThmDirichletInterval}),
upper bounds on the moments of Gamma distributions (see
Lemma 1 {in the supplementary material}), and a classification result on natural exponential family (see \autoref{MainRes} {in the appendix}), which are of independent interest. Our simulation
studies show that the new estimators often perform much better than the MR
estimator and Storey's estimator.

The computer codes for implementing the proposed proportion estimators {are written in the \textrm{R} language and} can be found on the author's website at \url{http://archive.math.wsu.edu/faculty/xchen/welcome.php} or \url{https://xiongzhichen.github.io/}.

\subsection{Notations and conventions}

The notations and conventions we will use throughout are stated as follows:
$C$ denotes a generic, positive constant whose values may differ at different
occurrences; $O\left(  \cdot\right)  $ and $o\left(  \cdot\right)  $ are
respectively Landau's big O and small o notations; $\mathbb{E}$, $\mathbb{V}$
and $\mathsf{cov}$ are respectively the expectation, variance and covariance
with respect to the probability measure $\Pr$; $\mathbb{R}$ and $\mathbb{C}$
are respectively the set of real and complex numbers; $\Re$ denotes the real
part of a complex number; $\mathbb{N}$ denotes the set of non-negative
integers, and $\mathbb{N}_{+}=\mathbb{N}\setminus\left\{  0\right\}  $;
$\nu$ the Lebesgue
measure, and when no confusion arises, the usual notation $d\cdot$ for
differential will be used in place of $\nu\left(  d\cdot\right)  $; for a
real-valued (measurable) function $f$ defined on some (measurable)
$A\subseteq\mathbb{R}$, $\left\Vert f\right\Vert _{p}=\left\{  \int%
_{A}\left\vert f\left(  x\right)  \right\vert ^{p}\nu\left(  dx\right)
\right\}  ^{1/p}$ and $L^{p}\left(  A\right)  =\left\{  f:\left\Vert
f\right\Vert _{p}<\infty\right\}  $ for $1\leq p < \infty$, $\left\Vert
f\right\Vert _{\infty}$ is its essential supremum, and $\left\Vert
f\right\Vert _{\mathrm{TV}}$ is the total variation of $f$ on $A$ {when $A$ is a non-empty closed interval of $\mathbb{R}$}; for a set
$A\subseteq\mathbb{R}^{d}$, $\left\vert A\right\vert $ is the cardinality of
$A$, and $1_{A}$ the indicator of $A$; $\partial_{\cdot}$ denotes the
derivative with respect to the subscript; {for a set $\mathcal{S}$, $\mathcal{S}^{\mathbb{N}}$ is the
$\aleph$-Cartesian product of $\mathcal{S}$, where $\aleph$ is the cardinality
of $\mathbb{N}$, i.e., an element of $\mathcal{S}^{\mathbb{N}}$ is a countably infinite, ordered sequence of elements of $\mathcal{S}$};
{for a Lebesgue
measurable, bivariate function $f$ defined on $\mathbb{R}^{2}$, the integral
$\int_{\mathbb{R}}dy\int_{\mathbb{R}}f\left(  x,y\right)  dx$ is understood as the
iterated integral $\int_{\mathbb{R}}dy\left(  \int_{\mathbb{R}}f\left(
x,y\right)  dx\right)  $}.

\subsection{Organization of article}

The rest of the article is organized as follows. We
formulate in \autoref{SecPre} the problem of proportion estimation, provide the needed
background, and give an overview of the main constructions to be proposed. We develop in \autoref{secNEF} uniformly consistent proportion
estimators for multiple testing the means of random variables from
the Gamma family,
and extend in \autoref{SecExtensions} one such construction to estimate
``proportions'' induced by a continuous function of a bounded null that is of bounded variation. We end the article with a discussion in \autoref{SecConcAndDisc}. {
Auxiliary results are given in \autoref{AppProofA}, and long technical proofs in the supplementary material.}

\section{Preliminaries}

\label{SecPre}

We formulate in \autoref{secModel} the problem of proportion estimation and in
\autoref{SecStrategy} the strategy to proportion estimation via Lebesgue-Stieltjes integral equations which further
generalizes that in \cite{Chen:2018a}, and {provide in \autoref{SecBackground} a
very brief background on the Gamma family and in \autoref{secResultsPoint} and \autoref{SecIllustration} an overview of the constructions.}

\subsection{The estimation problem}

\label{secModel}

Let $\Theta_{0}$ be a subset of $\mathbb{R}$ that has a non-empty interior and
non-empty complement $\Theta_{1}=\mathbb{R}\setminus\Theta_{0}$. For each
$i\in\left\{  1,\ldots,m\right\}  $, let $z_{i}$ be a random variable with
mean $\mu_{i}$, such that, for some integer $m_{0}$ between $0$ and
$m$, $\mu_{i}\in\Theta_{0}$ for each  $i\in I_{0,m}=\left\{  1,\ldots
,m_{0}\right\}  $ and $\mu_{i}\in\Theta_{1}$ for each $i\in I_{1,m}=\left\{
m_{0}+1,\ldots,m\right\}  $.
Consider simultaneously testing the null hypothesis $H_{i0}:\mu_{i}\in
\Theta_{0}$ versus the alternative hypothesis $H_{i1}:\mu_{i}\in\Theta_{1}$
for all $i\in\left\{  1,\ldots,m\right\}  $. Then the proportion
of true null hypothesis (``null proportion'' for short) is defined as $\pi_{0,m}=m^{-1} \left\vert I_{0,m}\right\vert =m^{-1}\sum_{i=1}^{m}1_{\Theta_{0}}\left(  \mu_{i}\right)$, and the proportion of
false null hypotheses (``alternative proportion'' for short) $\pi_{1,m}=1-\pi_{0,m}$. In other
words,%
\begin{equation}
\pi_{1,m}=m^{-1} \left\vert I_{1,m}\right\vert = m^{-1}\sum_{i=1}^{m}1_{\Theta_{1}}\left(  \mu_{i}\right).    \label{defPi}%
\end{equation}
Our target is to consistently estimate $\pi_{1,m}$ as
$m\rightarrow\infty$, and we will focus on the ``bounded
null'' $\Theta_{0}=\left(  a,b\right)  $ for some fixed,
finite $a,b\in U$ with $a<b$ and the ``one-sided
null'' $\Theta_{0}=\left(  -\infty,b\right)  $, both of which
are composite nulls. However, the strategy to be introduced next in
\autoref{SecStrategy} to achieve this target applies to general $\Theta_{0}$
(and hence general $\Theta_{1}$).

\subsection{The strategy via solutions to Lebesgue-Stieltjes integral
equations}

\label{SecStrategy}

Let $\mathbf{z}=\left(  z_{1},\ldots,z_{m}\right)  ^{\top}$ and
$\boldsymbol{\mu}=\left(  \mu_{1},...,\mu_{m}\right)  ^{\top}$. Denote by
$F_{\mu_{i}}$ the CDF of $z_{i}$ for $i\in\left\{  1,\ldots,m\right\}  $ and
suppose each $F_{\mu_{i}}$ is a member of a set $\mathcal{F}$ of CDFs such
that $\mathcal{F}=\left\{  F_{\mu}:\mu\in U\right\}  $ for some non-empty, known $U$
in $\mathbb{R}$. For the rest of the paper, we assume that each $F_{\mu}$ is
uniquely determined by $\mu$ and that $U$ has a non-empty interior. Recall the
definition of $\pi_{1,m}$ in (\ref{defPi}). The strategy to estimate
$\pi_{1,m}$ via Lebesgue-Stieltjes integral equations approximates each indicator function $1_{\Theta_{1}}\left(
\mu_{i}\right)  $, and is stated below.

Suppose for each fixed $\mu\in U$, we can approximate the indicator function
$1_{\Theta_{0}}\left(  \mu\right)  $ by

\begin{itemize}
\item[C1)] A ``discriminant function'' $\psi\left(  t,\mu\right)  $ satisfying $\lim_{t\rightarrow\infty}%
\psi\left(  t,\mu\right)  =1_{\Theta_{0}}\left(  \mu\right)$, and

\item[C2)] A ``matching function'' $K:\mathbb{R}^{2}\rightarrow\mathbb{R}$ that does not depend on any $\mu
\in\Theta_{1}$ and satisfies the Lebesgue-Stieltjes integral equation%
\begin{equation}
\psi\left(  t,\mu\right)  =\int K\left(  t,x\right)  dF_{\mu}\left(  x\right) { = \mathbb{E}_{Z \sim F_{\mu}}\left[K\left(t,Z\right)\right]}
, \forall\mu\in U. \label{eq3}%
\end{equation}

\end{itemize}

\noindent Then the ``average discriminant function''
\begin{equation}
\varphi_{m}\left(  t,\boldsymbol{\mu}\right)  =m^{-1}\sum_{i=1}^{m}\left\{
1-\psi\left(  t,\mu_{i}\right)  \right\}  \label{eq4a}%
\end{equation}
satisfies $\lim_{t\rightarrow\infty}\varphi_{m}\left(  t,\boldsymbol{\mu
}\right)  =\pi_{1,m}$ for any fixed $m$ and $\boldsymbol{\mu}$.
Further, the ``empirical matching function''
\begin{equation}
\hat{\varphi}_{m}\left(  t,\mathbf{z}\right)  =m^{-1}\sum_{i=1}^{m}\left\{
1-K\left(  t,z_{i}\right)  \right\}  \label{eq4b}%
\end{equation}
satisfies $\mathbb{E}\left\{  \hat{\varphi}_{m}\left(  t,\mathbf{z}\right)
\right\}  =\varphi_{m}\left(  t,\boldsymbol{\mu}\right)  $ for any fixed $m,t$
and $\boldsymbol{\mu}$. Namely, $\hat{\varphi}_{m}\left(  t,\mathbf{z}\right)
$ is an unbiased estimator of $\varphi_{m}\left(  t,\boldsymbol{\mu}\right)
$. We will reserve the notation $\left(  \psi,K\right)  $ for a pair of
discriminant function and matching function and the notations $\varphi
_{m}\left(  t,\boldsymbol{\mu}\right)  $ and $\hat{\varphi}_{m}\left(
t,\mathbf{z}\right)  $ as per (\ref{eq4a}) and (\ref{eq4b}) unless otherwise
noted. The concept of discriminant function and matching function originates
from, is inspired by, and extends the concept of ``phase
functions'' in the work of \cite{Jin:2008}. The pair $\left(
\psi,K\right)  $ presented here has those in \cite{Chen:2018a} and
\cite{Jin:2008} as special cases, and is the most general form for the purpose
of proportion estimation. It converts proportion estimation into solving a
specific Lebesgue-Stieltjes integral equation.
Note that in
general it is very hard and sometimes impossible to find a function
$K:\mathbb{R}\rightarrow\mathbb{R}$ such that $\int K\left(  x\right)
dF_{\mu}\left(  x\right)  =1_{\Theta_{0}}\left(  \mu\right)  ,\forall\mu\in U$
when $\Theta_{0}\neq\mathbb{R}$ due to the discontinuity of the indicator
function $1_{\Theta_{0}}\left(  \cdot\right)  $. Take for example a null set $\Theta_{0}\neq\mathbb{R}$ and the family of
Gaussian distributions with mean $\mu\in\mathbb{R}$ and variance $1$, i.e.,
$f_{0}\left(  \mu-x\right)  $, where $f_{0}\left(  x\right)  =\left(
2\pi\right)  ^{-1/2}\exp\left(  -2^{-1}x^{2}\right)  $. Then%
\[
\int K\left(  x\right)  dF_{\mu}\left(  x\right)  =
\int K\left(  x\right)  f_{0}\left(  \mu-x\right)  dx=
\left(  K\ast f_{0}\right)  \left(  \mu\right)  .
\]
Since $f_{0}$ is uniformly bounded on $\mathbb{R}$, if $K\in L^{1}\left(
\mathbb{R}\right)  $, then $K\ast f_{0}$ is uniformly continuous on
$\mathbb{R}$, and it will never be equal to $1_{\Theta_{0}}\left(  \cdot\right)
$.

When the difference
\begin{equation}
e_{m}\left(  t\right)  =\hat{\varphi}_{m}\left(  t,\mathbf{z}\right)
-\varphi_{m}\left(  t,\boldsymbol{\mu}\right)  \label{eq2d}%
\end{equation}
is small for large $t$, $\hat{\varphi}_{m}\left(  t,\mathbf{z}\right)  $ will
accurately estimate $\pi_{1,m}$. Since $\varphi_{m}\left(  t,\boldsymbol{\mu
}\right)  =\pi_{1,m}$ or $\hat{\varphi}_{m}\left(  t,\mathbf{z}\right)
=\pi_{1,m}$ rarely happens, $\hat{\varphi}_{m}\left(  t,\mathbf{z}\right)  $
usually employs an increasing sequence $\left\{  t_{m}\right\}  _{m\geq1}$
with $\lim_{m\rightarrow\infty}t_{m}=\infty$ in order to achieve consistency
in the sense that%
\begin{equation}
\left\vert \pi_{1,m}^{-1}\hat{\varphi}_{m}\left(  t_{m},\mathbf{z}\right)
-1\right\vert \rightsquigarrow0\text{ as }m\rightarrow\infty, \label{defConsistency}%
\end{equation}
where ``$\rightsquigarrow$'' denotes ``convergence in probability''.
Following the convention set by \cite{Chen:2018a}, we refer to $t_{m}$ as the
``speed of convergence'' of $\hat{\varphi
}_{m}\left(  t_{m},\mathbf{z}\right)  $. Throughout the paper, consistency of
a proportion estimator is defined via (\ref{defConsistency}) to accommodate
the scenario $\lim_{m\rightarrow\infty}\pi_{1,m}=0$. Further, the accuracy of
$\hat{\varphi}_{m}\left(  t_{m},\mathbf{z}\right)  $ in terms of estimating
$\pi_{1,m}$ and its speed of convergence depend on how fast $\pi_{1,m}%
^{-1}e_{m}\left(  t_{m}\right)  $ converges to $0$ and how fast $\pi
_{1,m}^{-1}\varphi_{m}\left(  t_{m},\boldsymbol{\mu}\right)  $ converges to
$1$. This general principle also applies to the works of \cite{Jin:2008,Jin:2007} and
\cite{Chen:2018a}.

By duality, $\textcolor{black}{\psi_{m}}\left(  t,\boldsymbol{\mu}\right)
=1-\varphi_{m}\left(  t,\boldsymbol{\mu}\right)  $ satisfies $\pi_{0,m}%
=\lim_{t\rightarrow\infty}\psi_{m}\left(  t,\boldsymbol{\mu}\right)
$ for any fixed $m$ and $\boldsymbol{\mu}$, and $\textcolor{black}{\hat{\psi}_{m}}
\left(  t,\mathbf{z}\right)  =1-\hat{\varphi}_{m}\left(  t,\mathbf{z}\right)
$ satisfies $\mathbb{E}\left[  \hat{\psi}_{m}\left(  t,\mathbf{z}%
\right)  \right]  =\psi_{m}\left(  t,\boldsymbol{\mu}\right)  $ for
any fixed $m,t$ and $\boldsymbol{\mu}$. Moreover, $\hat{\psi}_{m}
\left(  t,\mathbf{z}\right)  $ will accurately estimate $\pi_{0,m}$ when
$e_{m}\left(  t\right)  $ is suitably small for large $t$, and the stochastic
oscillations of $\hat{\psi}_{m}\left(  t,\mathbf{z}\right)  $ and
$\hat{\varphi}_{m}\left(  t,\mathbf{z}\right)  $ are the same and is
quantified by $e_{m}\left(  t\right)  $.


{In order to better quantify the consistency of the proposed estimators,} the following definition is adapted from \cite{Chen:2018a}:

\begin{definition}
\label{DefUniformConsistency}Given a family $\mathcal{F}$, consider a
non-empty set $\mathcal{Q}\left(  \mathcal{F}\right)
{=\mathcal{Q}_1\left(  \mathcal{F}\right) \times \mathcal{Q}_2\left(  \mathcal{F}\right)}
\subseteq\mathbb{R}^{\mathbb{N}}\times
{\left(\mathcal{P}\left(\mathbb{R}\right)\right)^{\mathbb{N}}}$
and its coordinate projection $\mathcal{Q}_{m}\left(  \mathcal{F}\right)
=\mathcal{Q}_{m,1}\left(  \mathcal{F}\right)  \times\mathcal{Q}_{m,2}\left(
\mathcal{F}\right)  $ onto $\mathbb{R}^{m}\times
{\mathcal{P}\left(\mathbb{R}\right)}$ for each $m\in
\mathbb{N}_{+}$,
{where $\mathcal{P}\left(\mathbb{R}\right)$ is the power set of $\mathbb{R}$, $\mathcal{Q}_{m,1}\left(  \mathcal{F}\right) $ the first $m$ coordinates of $\mathcal{Q}_1\left(  \mathcal{F}\right)$, and $\mathcal{Q}_{m,2}\left(
\mathcal{F}\right) $ the $m$-th coordinate of $\mathcal{Q}_2\left(  \mathcal{F}\right)$}.
Then $\mathcal{Q}\left(  \mathcal{F}\right)  $ is said to be a
\textquotedblleft uniform consistency class\textquotedblright\ for
$\hat{\varphi}_{m}\left(  t,\mathbf{z}\right)  $ and $\hat{\varphi}_{m}\left(
t,\mathbf{z}\right)  $ is said to be \textquotedblleft uniformly
consistent\textquotedblright\ on $\mathcal{Q}\left(  \mathcal{F}\right)  $ if
\begin{equation}
\sup\nolimits_{\boldsymbol{\mu}\mathcal{\in Q}_{m,1}\left(  \mathcal{F}%
\right)  }\left\vert \pi_{1,m}^{-1}\sup\nolimits_{t\in\mathcal{Q}_{m,2}\left(
\mathcal{F}\right)  }\hat{\varphi}_{m}\left(  t,\mathbf{z}\right)
-1\right\vert \rightsquigarrow0\,\ \text{as }m\rightarrow\infty.
\label{defUCC}%
\end{equation}

\end{definition}

In words, a uniform consistency class for an estimator $\hat{\varphi}%
_{m}\left(  t,\mathbf{z}\right)  $ is the asymptotic setting (as $m \to\infty
$) for the family $\mathcal{F}$ under which $\hat{\varphi}_{m}\left(
t,\mathbf{z}\right)  $ at its maximal speed of convergence (that is indicated
by $\sup\nolimits_{t\in\mathcal{Q}_{m,2}\left(  \mathcal{F} \right)  }$ in
\autoref{DefUniformConsistency}) still maintains consistency uniformly over
all settings of all $F_{\mu_{i}}, i=1,\ldots,m$ (that is indicated by
$\sup\nolimits_{\boldsymbol{\mu}\mathcal{\in Q}_{m}\left(  \mathcal{F}\right)
}$ in \autoref{DefUniformConsistency}) and hence uniformly over all settings
of the alternative proportion $\pi_{1,m}$ that are induced by such $F_{\mu
_{i}}, i=1,\ldots,m$.

\subsection{The Gamma family}

\label{SecBackground}

The Gamma family $\mathcal{F}=\mathsf{Gamma}\left(  \theta,\sigma\right)  $ has generating measure (called ``basis'') $\beta$ such that%
\[
\frac{d \beta}{d \nu} \left( x\right)  = {\left\{  \Gamma\left(
\sigma\right)  \right\}  ^{-1} e^{-x} x^{\sigma-1} }1_{\left(  0,\infty\right)  }\left(  x\right),
\]
where $\Gamma$ is the Gamma function and $\sigma >0$ is the scale parameter. So, this family of CDFs is indexed by its natural parameter $\theta \in \Theta = \{\theta:\theta<1\}$, contains CDF $G_{\theta}, \theta \in \Theta$ which has density
\[
f_{\theta}\left(  x\right)  =\left\{  \Gamma\left(  \sigma\right)  \right\}
^{-1}\left(  1-\theta\right)  ^{\sigma}e^{-\left(1-\theta\right) x}x^{\sigma-1}%
1_{\left(  0,\infty\right)  }\left(  x\right)  ,
\]
and has mean function $\mu=\mu\left(  \theta\right) =\int_{0}^{\infty}xf_{\theta}\left(  x\right)  dx =\sigma\left(  1-\theta\right)  ^{-1}$ with $\mu \in U=\mu\left(\Theta\right)$. Note that $f_0 = \frac{d \beta}{d \nu}$.
If we set $\tilde{\theta}=1-\theta$, then we can write $f_{\theta}\left(
x\right)  $ more compactly as%
\[
\tilde{f}_{\tilde{\theta}}\left(  x\right)  =\left\{  \Gamma\left(
\sigma\right)  \right\}  ^{-1}\tilde{\theta}^{\sigma}e^{-\tilde{\theta}%
x}x^{\sigma-1}1_{\left(  0,\infty\right)  }\left(  x\right)
\]
and write its CDF as $\tilde{G}_{\tilde{\theta}}$. Note that $\tilde
{f}_{\tilde{\theta}}\left(  x\right)  =f_{\theta}\left(  x\right)  $ and
$G_{\theta}=$\ $\tilde{G}_{\tilde{\theta}}$ and that $\beta=G_{0}%
=\tilde{G}_{1}$. With these notational
adjustments, we quickly see that the Gamma family has the following
``scaling-invariance property'': if $X$ has
CDF $\tilde{G}_{\tilde{\theta}}$, then $Y=\tilde{\theta}X$ has CDF $\beta$.

The Gamma family has moment functions
\[
\tilde{c}_{n}\left(  \theta\right)  =\int x^{n}G_{\theta}\left(  dx\right)=\Gamma\left(  n+\sigma\right)  \left(
1-\theta\right)  ^{-n}\left\{  \Gamma\left(  \sigma\right)  \right\}  ^{-1} \quad \text{for} \quad n {\in \mathbb{N}},
\]
where we recognize that $\int x^{z}G_{\theta}\left(  dx\right),z \in \mathbb{C}$ is the Mellin transform of $G_{\theta}$.
We see $\tilde{c}_{n}\left(  \theta\right)  =\xi^{n}\left(  \theta\right)
\zeta\left(  \theta\right)  \tilde{a}_{n}$ for each $n\in\mathbb{N}$ and
$\theta\in\Theta$, where $\xi\left(  \theta\right)  =\left(  1-\theta\right)  ^{-1}$, $\tilde
{a}_{n}=\Gamma\left(  n+\sigma\right)  \left\{  \Gamma\left(  \sigma\right)
\right\}  ^{-1}$ and $\zeta \left(\theta\right)=\zeta_{0}$ for the constant $\zeta_{0}=1$.
The special structure for the moment functions $\tilde{c}_{n}$ is a consequence of the scaling-invariance property of the Gamma family.
Further, $\tilde{a}_{1}=\sigma$, $f_{\theta}\left(  x\right)  =O\left(  x^{\sigma
-1}\right)  $ as $x\rightarrow0+$, and $\Psi\left(  t,\theta\right)  =\sum_{n=0}^{\infty}\frac{t^{n}\xi
^{n}\left(  \theta\right)  }{\tilde{a}_{n}n!}$ is absolutely convergent
pointwise in $\left(  t,\theta\right)  \in\mathbb{R}\times\Theta$. Note that $\mu\left(  \theta\right)
=\xi\left(  \theta\right)  \zeta\left(  \theta\right)  \tilde{a}_{1}$, that $\mu=\mu\left(  \theta\right)$ has a unique inverse function $\theta=\theta(\mu) \in \Theta$,
and that each $G_{\theta}$ corresponds uniquely to a non-degenerate CDF $F_{\mu}$ and vice versa, i.e., $\mathcal{F}=\{G_{\theta}: \theta \in \Theta\}=\{F_{\mu}: \mu \in U\}$.

For random variables whose CDFs are members of the Gamma family, testing their means $\mu_{i}$'s simultaneously is equivalent to testing
natural parameters $\theta_{i}=\theta\left(  \mu_{i}\right)  $ simultaneously. So, for the
constructions of proportion estimators via the strategy in
\autoref{SecStrategy}, $\psi$ will also be regarded as functions of $\theta$. Specifically, $\psi$ defined by (\ref{eq3}) becomes%
\[
\psi\left(  t,\theta\right)  =\int K\left(  t,x\right)  dG_{\theta}\left(
x\right)  \text{ for }G_{\theta}\in\mathcal{F}.
\]
Let $\boldsymbol{\theta}=\left(  \theta_{1},\ldots,\theta_{m}\right)^{\top}  $. Then
accordingly $\varphi
_{m}\left(  t,\boldsymbol{\theta}\right)  =m^{-1}\sum_{i=1}^{m}\left\{
1-\psi\left(  t,\theta_{i}\right)  \right\}  $ becomes the counterpart of
(\ref{eq4a}).

The Gamma family is a natural exponential family (NEF) that contains the family of exponential distributions and the family of central Chi-square distributions, and it has a sequence of very special moment functions $\{\tilde{c}_n,n \ge 0\}$ such that $\tilde{c}_{n}\left(  \theta\right)  =\xi^{n}\left(  \theta\right)
\tilde{a}_{n}$ for each $n\in\mathbb{N}$ and
$\theta\in\Theta$ (given earlier), which we refer to as ``separable moments''. \autoref{AppProofA} provides a very short background on NEFs and \autoref{MainRes} there shows that the only NEF with separable moments is the Gamma family.
As revealed by \cite{Chen:2018a} and will be seen later, the key appeal of the Gamma family in constructing proportion estimators for both a one-sided null and a bounded null
respectively is that its scaling-invariant property utilizes the multiplicative structure of the multiplicative group $\left(  \mathbb{R}_{>0},\times\right)  $, i.e., the set $\mathbb{R}_{>0}$ of positive real numbers with multiplication ``$\times$'',
and couples well with convolution on a multiplicative group, Mellin transform, and Dirichlet integrals to provide solutions to the Lebesgue-Stieltjes integral equation (\ref{eq3}).

\subsection{Constructions of proportion estimators for a point null and Dirichlet integrals}
\label{secResultsPoint}

Since the bounded null $\Theta _{0}=\left(  a,b\right)  $ and the one-sided null $\Theta_{0}=\left(  -\infty,b\right)  $ have two
finite boundary points $a$ and $b$, in order to consistently estimate
$\pi_{1,m}=m^{-1}\sum_{i=1}^{m}1_{\Theta_{1}}\left(  \mu_{i}\right)  $
using the strategy in \autoref{SecStrategy} for these two composite
hypotheses, we can use the proportion estimators of \cite{Chen:2018a} for a
point null to specifically account for the proportion of $\theta_{i}$'s that are equal to $\theta\left(
a\right)  $ or $\theta\left(  b\right)  $ when $\mathcal{F}$ is the Gamma family. Consider a point
null $\Theta_{0}=\left\{\theta_{0}\right\}  $ for a fixed $\theta_{0}\in\Theta$.
We can state the constructions of \cite{Chen:2018a} of
uniformly consistent estimators of $\pi_{1,m}$ for the point null when
$\mathcal{F}$ is the Gamma family as
follows in terms of a discriminant function and a matching function that satisfy conditions C1) and C2):

\begin{theorem}
\label{ThmPoinNull}Let $\omega$ be an even, bounded, probability density
function on $\left[  -1,1\right]  $. For
$\theta^{\prime}\in\Theta$, define
\[
K_{3,0}\left(  t,x;\theta^{\prime}\right)  ={\int_{\left[  -1,1\right]  }}%
\sum_{n=0}^{\infty}\left(  \tilde{a}_{n}n!\right)  ^{-1}\left(  -tsx\right)
^{n}\cos\left\{  2^{-1}\pi n+ts\xi\left(  \theta^{\prime}\right)  \right\}
\omega\left(  s\right)  ds
\]
and let%
\[
\psi_{3,0}\left(  t,\theta;\theta^{\prime}\right)  =\int K_{3,0}\left(
t,x;\theta^{\prime}\right)  dG_{\theta}\left(  x\right)  .%
\]
Then
\[
\psi_{3,0}\left(  t,\theta;\theta^{\prime}\right)  ={\int_{\left[
-1,1\right]  }}\cos\left[  ts\left\{  \xi\left(  \theta^{\prime}\right)
-\xi\left(  \theta\right)  \right\}  \right]  \omega\left(  s\right)  ds.
\]
For the point null $\Theta_{0}=\left\{
\theta_{0}\right\}  $ with $\theta_{0}\in\Theta$, $\left(  \psi,K\right)
=\left(  \psi_{3,0}\left(  t,\theta;\theta_{0}\right)  ,K_{3,0}\left(
t,x;\theta_{0}\right)  \right)  $ when $\mathcal{F}$ is the Gamma family. In
particular,
$\psi_{3,0}\left(  t,\theta_{0};\theta_{0}\right)  =1$ for all $t$.
\end{theorem}

In this work, we will further assume that $\omega$ is of bounded variation
unless otherwise noted. For example, the triangular density $\omega\left(s\right)=\left(1-\vert s\vert\right)1_{\left[-1,1\right]}\left(s\right)$ or the uniform density $\omega\left(s\right)=0.5 \times 1_{\left[-1,1\right]}\left(s\right)$ can be used. \autoref{ThmPoinNull} will be used by the constructions to be introduced in
\autoref{secNEF} and \autoref{SecExtensions}.

Also, we need two Dirichlet integrals%

\[
\mathcal{D}_{1}\left(  t,\mu;a,b\right)  =\frac{1}{\pi}\int_{\left(
\mu-b\right)  t}^{\left(  \mu-a\right)  t}\frac{\sin y}{y}dy\text{ \ \ and \
}\mathcal{D}_{2}\left(  t,\mu;b\right)  =\frac{1}{\pi}\int_{0}^{t}\frac
{\sin\left\{  \left(  \mu-b\right)  y\right\}  }{y}dy
\]
and%

\begin{equation}
\mathcal{D}_{1,\infty}\left(  \mu;a,b\right)  =\lim_{t\rightarrow\infty
}\mathcal{D}_{1}\left(  t,\mu;a,b\right)  =\left\{
\begin{array}
[c]{lll}%
1 & \text{if} & a<\mu<b\\
2^{-1} & \text{if} & \mu=a\text{ or }\mu=b\\
0 & \text{if} & \mu<a\text{ or }\mu>b
\end{array}
\right.  ,\label{EqDirichlet1}%
\end{equation}
and
\begin{equation}
\mathcal{D}_{2,\infty}\left(  \mu;b\right)  =\lim_{t\rightarrow\infty
}\mathcal{D}_{2}\left(  t,\mu;b\right)  =\left\{
\begin{array}
[c]{lll}%
2^{-1} & \text{if} & \mu>b\\
0 & \text{if} & \mu=b\\
-2^{-1} & \text{if} & \mu<b
\end{array}
\right.  .\label{EqDB1}%
\end{equation}
{These integrals will be used here to construct discriminant functions.}
The identities (\ref{EqDirichlet1}) and (\ref{EqDB1}) will later be proved by \autoref{lm:Dirichlet} in \autoref{AppProofA}.

\subsection{An overview on the new constructions}
\label{SecIllustration}

We are ready to give an overview on the main (but very subtle) constructions to
be introduced in later sections. Even though conceptually these constructions
can be more elegantly stated and better understood in terms of complex analysis (as
their proofs reveal), we describe them using terms of real analysis wherever
feasible. Let $\circledast$ denote the
multiplicative convolution with respect to $\left(  \mathbb{R}_{>0}%
,\times\right)  $.
For two functions $\tilde{f}$ and $\tilde{g}$, recall
their multiplicative convolution as $(  \tilde{f}\circledast\tilde
{g})  \left(  z\right)  =\int\tilde{g}\left(  y\right)  \tilde{f}\left(
zy^{-1}\right)  \tilde{\nu}\left(  dy\right)  $, where $\tilde{\nu}$ is a Haar
measure on the multiplicative group $\left(  \mathbb{R}_{>0},\times\right)  $;
see, e.g., Section 1.1.6 of \cite{rudin2017fourier} for such a
definition.

Roughly speaking, the pair $\left(  \psi,K\right)  $ for the Gamma
family is obtained as follows. Pick a discriminant function $\psi\left(
t,\theta\right)  $ and form its complex version $\psi^{\dagger}\left(
t,\theta\right)  $ such that $\psi\left(  t,\theta\right)  =\Re\left\{
\psi^{\dagger}\left(  t,\theta\right)  \right\}  $. Define $\psi^{\dagger
}\left(  t,\theta\right)  =\int K^{\dagger}\left(  t,x\right)  dG_{\theta
}\left(  x\right)  $ for some complex-valued function $K^{\dagger}\left(
t,x\right)  $. Recall the facts about Gamma family stated in
\autoref{SecBackground}, including $\mu=\mu\left(  \theta\right)
=\sigma\left(  1-\theta\right)  ^{-1}$, $\beta=G_{0}=\tilde{G}_{1}$,
$G_{\theta}=$\ $\tilde{G}_{\tilde{\theta}}$, $\tilde{\theta}=1-\theta$ and
that $X$ has CDF $\tilde{G}_{\tilde{\theta}}$ implying $Y=\tilde{\theta}X$ has
CDF $\beta$. We see that
\begin{equation}
\psi^{\dagger}\left(  t,\theta\right) =\int K^{\dagger
}\left(  t,x\right)  d\tilde{G}_{\tilde{\theta}}\left(  x\right)  =\int
K^{\dagger}\left(  t,x\tilde{\theta}^{{-1}}\right)  f_{0}\left(  x\right)
dx = \tilde{\psi}^{\dagger}\left(  t,\tilde{\theta}^{-1}\right)  .\label{ExplainGammaStrategy}%
\end{equation}
for some function $\tilde{\psi}^{\dagger} (  t,\tilde{\theta}^{-1} )$.
Note that (\ref{ExplainGammaStrategy}) is not the multiplicative convolution $K^{\dagger}
\circledast f_{0}$ with respect to the Lebesgue measure $\nu$, and that Fourier analysis with respect to $\nu$ and $\times$ cannot be directly applied to invert $\tilde{\psi}^{\dagger}$ obtain $K^{\dagger}$ (and then set $K=\Re\left\{  K^{\dagger}\right\}  $).
However, $f_{0}\left(  x\right)  dx=h_{0}\left(  y\right)  \nu
_{-1}\left(  dy\right)  $ with $y=x^{-1}$ for $x\in\mathbb{R}_{>0}$, where
$h_{0}\left(  y\right)  =y^{-\sigma}e^{-1/y}\left\{  \Gamma\left(
\sigma\right)  \right\}  ^{-1}1_{\left(  0,\infty\right)  }\left(  y\right)  $
and $\nu_{-1}\left(  A\right)  = \int_{A}y^{-1}dy$ for measurable
$A\subseteq\mathbb{R}_{>0}$ is the Haar measure on the group $\left(
\mathbb{R}_{>0},\times\right)  $.
So,
(\ref{ExplainGammaStrategy}) is just%
\[
\tilde{\psi}^{\dagger}\left(  t,\tilde{\theta}^{-1}\right)  =\int K^{\dagger
}\left(  t,y^{-1}\tilde{\theta}^{-1}\right)  h_{0}\left(  y\right)  \nu
_{-1}\left(  dy\right)  =\left(  K^{\dagger}\circledast h_{0}\right)  \left(
t,\tilde{\theta}^{-1}\right)  ,
\]
and taking the real parts of $\tilde{\psi}^{\dagger}$ and $K^{\dagger}$ as
$\tilde{\psi}$ and $K$ respectively, we obtain
\[
\tilde{\psi}\left(  t,\tilde{\theta}^{-1}\right)  =\int K\left(
t,y^{-1}\tilde{\theta}^{-1}\right)  h_{0}\left(  y\right)  \nu_{-1}\left(
dy\right)  =\left(  K\circledast h_{0}\right)  \left(  t,\tilde{\theta}%
^{-1}\right)  = \psi\left(t,\theta\right) .
\]
Consequently, $K$ may be found via the Fourier inversion of $\tilde{\psi}$ or $\psi$ (if this inversion is well-defined) with
respect to $\nu_{-1}$ on $\left(  \mathbb{R}_{>0},\times\right)  $.
For details on this construction and related derivations, please see Theorem 5 of \cite{Chen:2018a} and its proof, and \autoref{ThmConstructionMoments} and \autoref{TmVNEF} and their proofs.

\section{Constructions for Gamma family}

\label{secNEF}

Recall $\mathcal{F}=\left\{  F_{\mu}:\mu\in U\right\}  =\left\{  G_{\theta
}:\theta\in\Theta\right\}  $ for the Gamma family. We will refer to as
``Construction I'' the construction of estimators of
$\pi_{1,m}$ for the bounded null $\Theta_{0}=\left(  a,b\right)  \cap U$ for fixed, finite $a,b\in U$ with
$a<b$, and as ``Construction II'' the
construction of estimators of $\pi_{1,m}$ for the one-sided
null $\Theta_{0}=\left(  -\infty,b\right)  \cap U$ for a
fixed, finite $b\in U$. Both constructions utilize \autoref{ThmPoinNull}, Dirichlet integrals,
and their integral representations provided in \autoref{AppProofA}.
Note that the bounded null $\Theta_{0}=\left(  a,b\right)  \cap
U$ and the one-sided null $\Theta_{0}=\left(  -\infty,b\right)  \cap U$ have to
be convex sets. We write
$\theta\left(  \mu\right)  $ as $\theta_{\mu}$, so that $\theta_{a}=\theta\left(  a\right)  $ and
$\theta_{b}=\theta\left(  b\right)  $, and $\mu_{0}=\mu\left(  \theta
_{0}\right)  $, $a=\mu\left(  \theta_{a}\right)  $ and $b=\mu\left(
\theta_{b}\right)  $. Recall from \autoref{SecBackground} that $\mu\left(  \theta\right)  =\xi\left(
\theta\right)  \zeta\left(  \theta\right)  \tilde{a}_{1}$ for $\theta\in
\Theta$ with $\xi\left(  \theta\right)  =\left(  1-\theta\right)  ^{-1}$, $\zeta\left(  \theta\right)  \equiv\zeta_{0}=1$ and $\tilde{a}_{1}=\sigma$.

\subsection{The case of a bounded null}

Construction I for the bounded null for the Gamma family is provided below:

\begin{theorem}
\label{ThmConstructionMoments}When $\mathcal{F}$ is the Gamma family, set
\begin{equation}
K_{1}\left(  t,x\right)  =\frac{1}{2\pi}\int_{a}^{b}tdy\int_{-1}^{1}\sum
_{n=0}^{\infty}\frac{\left(  tsx\tilde{a}_{1}\right)  ^{n}\cos\left(
2^{-1}n\pi-tsy\right)  }{\tilde{a}_{n}n!}ds. \label{eq13g}%
\end{equation}
Then
\begin{equation}
\psi_{1}\left(  t,\theta\right)  ={\int}K_{1}\left(  t,x\right)  dG_{\theta
}\left(  x\right)
= \frac{1}{\pi}\int_{\left(  \mu\left(\theta\right)-b\right)  t}^{\left(  \mu\left(\theta\right)-a\right)
t}\frac{\sin y}{y}dy,
\label{psiGammaBounded}
\end{equation}
and the desired $\left(  \psi,K\right)  $ is
\begin{equation}
\left\{
\begin{array}
[c]{l}%
K\left(  t,x\right)  =K_{1}\left(  t,x\right)  -2^{-1}\left\{  K_{3,0}\left(
t,x;\theta_{a}\right)  +K_{3,0}\left(  t,x;\theta_{b}\right)  \right\} \\
\psi\left(  t,\theta\right)  =\psi_{1}\left(  t,\theta\right)  -2^{-1}\left\{
\psi_{3,0}\left(  t,\theta;\theta_{a}\right)  +\psi_{3,0}\left(
t,\theta;\theta_{b}\right)  \right\}
\end{array}
\right.  . \label{IV-b}%
\end{equation}

\end{theorem}

Note that the limit as $t \to \infty$ of $\psi_{1}\left(t,\theta\right)$ in (\ref{psiGammaBounded}) is the Dirichlet integral in (\ref{EqDirichlet1}) with $\mu$ as $\mu\left(\theta\right)$, which implies $\lim_{t \to \infty} \psi\left(  t,\theta\right) = 1_{\left(a,b\right)}\left(\mu\right)$ for $\psi\left(  t,\theta\right)$ in (\ref{IV-b}).
For the uniform consistency of the proposed estimator $\hat{\varphi}_{m}\left(  t,\mathbf{z}\right)$, set $\left\Vert 1-\boldsymbol{\theta}\right\Vert _{\infty}
=\max_{1\leq i\leq m}\left\{  1-\theta_{i}\right\}  $,
\begin{equation}
u_{3,m}=\min_{1\leq i\leq m}\left\{  1-\theta_{i}\right\}  \text{ }\ \text{and
}\tilde{u}_{3,m}=\min_{\tau\in\left\{  a,b\right\}  }\min_{\left\{
j:\theta_{j}\neq\theta_{\tau}\right\}  }\left\vert \xi\left(  \theta_{\tau
}\right)  -\xi\left(  \theta_{i}\right)  \right\vert . \label{eq12e}%
\end{equation}
Then we have:

\begin{theorem}
\label{ConcentrationIII}Suppose $\left\{  z_{j}\right\}  _{j=1}^{m}$ are
independent Gamma random variables with parameters $\left\{  \left(
\theta_{i},\sigma\right)  \right\}  _{i=1}^{m}$. Then,
when $t$ is positive and sufficiently large,
\[
\mathbb{V}\left\{  e_{m}\left(  t\right)  \right\}  \leq\frac{C\left(
1+t^{2}\right)  }{m^{2}}\exp\left(  \frac{4t\max\left\{  \sigma,1\right\}
}{u_{3,m}}\right)  \sum_{i=1}^{m}\left(  \frac{t}{1-\theta_{i}}\right)
^{3/4-\sigma}.
\]
Further, when $\sigma\geq 3/4$
\begin{equation}
\mathcal{Q}\left( \mathcal{F}\right)  =\left\{
\begin{array}
[c]{c}%
 t_m=4^{-1}(\max\left\{  \sigma,1\right\})^{-1}\gamma u_{3,m}\ln m,t_m^{-1}\left(  1+\tilde{u}_{3,m}
^{-1}\right)  =o\left(  \pi_{1,m}\right)  ,\\
t_m\rightarrow\infty,\left\Vert 1-\boldsymbol{\theta}\right\Vert _{\infty
}^{\sigma-3/4}t_m^{11/4-\sigma}=o\left(  m^{1-\gamma}\pi_{1,m}^{2}\right)
\end{array}
\right\}  \label{eq21}%
\end{equation}
for each $\gamma\in \left(  0,1\right)  $ is a uniform consistency class, for which $\gamma=1$ can be set when $\sigma>11/4$, and
when $\sigma\leq3/4$, a uniform consistency class is
\begin{equation}
\mathcal{Q}\left(\mathcal{F}\right)  =\left\{
\begin{array}
[c]{c}%
t_m=4^{-1}\gamma u_{3,m}\ln m,t_m^{-1}\left(  1+\tilde{u}_{3,m}^{-1}\right)
=o\left(  \pi_{1,m}\right)  ,\\
t_m\rightarrow\infty,\left(  \ln m\right)  ^{11/4-\sigma}u_{3,m}
^{2}=o\left(  m^{1-\gamma}\pi_{1,m}^{2}\right)
\end{array}
\right\}  \label{eq21a}%
\end{equation}
for each $\gamma\in\left(  0,1\right)  $.
\end{theorem}

\autoref{ConcentrationIII} reveals that the uniform consistency class of the estimator $\hat{\varphi}_{m}\left(  t,\mathbf{z}\right)$ changes with the scale parameter $\sigma$ of the Gamma family, which is similar to that for the estimator in the setting of a
location-shift family that has a scale parameter {in \cite{Chen:2018a}}. However, the expression (\ref{eq21a}) states that when $\sigma \le 3/4$, the maximal speed of convergence, {$t_{m}=4^{-1}\gamma u_{3,m}\ln m$,} of $\hat{\varphi}_{m}\left(  t,\mathbf{z}\right)$ that maintains uniform consistency is independent of $\sigma$,
{whereas $\left(  \gamma\ln m\right)  ^{11/4-\sigma}u_{3,m}%
^{2}=o\left(  m^{1-\gamma}\pi_{1,m}^{2}\right)  $, depending on $\sigma$,
controls the magnitude of $\pi_{1,m}$ that can be uniformly consistently
estimated at this speed $t_{m}$. }
Since $\theta<1$ for the Gamma family, $u_{3,m}$ measures how close the parameter $\theta_i$ of a
$G_{\theta_{i}}$ is to the singularity parameter $1$ for which a Gamma density is undefined, and
it is sensible to assume $\liminf_{m\rightarrow\infty}u_{3,m}>0$. On the other
hand, $\sigma\xi\left(  \theta\right)  =\mu\left(  \theta\right)  $ for all
$\theta\in\Theta$. So, $\tilde{u}_{3,m}$ measures the minimal
difference between the means $\mu\left(  \theta_{i}\right)  $ of
$G_{\theta_{i}}$ for $\theta_{i}\notin\left\{  \theta_{a},\theta_{b}\right\}
$ and the means $\mu\left(  \theta_{a}\right)  $ and $\mu\left(  \theta
_{b}\right)  $, and $\tilde{u}_{3,m}$ cannot be too small relative to $t$ as
$t\rightarrow\infty$ in order for the estimator induced by $K_{3,0}\left(
t,x;\theta_{a}\right)  $ and $K_{3,0}\left(  t,x;\theta_{b}\right)  $ (already given by \autoref{ThmPoinNull}) in
(\ref{IV-b}) to consistently estimate the proportions of means that are equal
to $\mu\left(  \theta_{a}\right)  $ or $\mu\left(  \theta_{b}\right)  $; see
Theorem 9 of \cite{Chen:2018a}. Finally, $\left\Vert 1-\boldsymbol{\theta
}\right\Vert _{\infty}$ measures the maximal range of $\left\{  \theta
_{i}\right\}  _{i=1}^{m}$ from $1$, and when $\sigma\geq11/4$, the larger
$\left\Vert 1-\boldsymbol{\theta}\right\Vert _{\infty}$ is, the more likely the estimator will have a slower
maximal speed of convergence to achieve consistency.

\subsection{The case of a one-sided null}

We present Construction II for the one-sided null for the Gamma family as

\begin{theorem}
\label{TmVNEF}When $\mathcal{F}$ is the Gamma family, set
\[
K_{1}\left(  t,x\right)  =\frac{1}{2\pi}\int_{0}^{1}tdy\int_{-1}^{1}\sum
_{n=0}^{\infty}\cos\left(  2^{-1}\pi n-tysb\right)  \frac{\left(  tys\right)
^{n}\left(  \tilde{a}_{1}x\right)  ^{n}}{n!}\left(  \frac{\tilde{a}_{1}%
x}{\tilde{a}_{n+1}}-\frac{b}{\tilde{a}_{n}}\right)  ds.
\]
Then%
\begin{equation}
\psi_{1}\left(  t,\theta\right)  =\int K_{1}\left(  t,x\right)  dG_{\theta
}\left(  x\right)
= \frac{1}{\pi}\int_{0}^{t}\frac{\sin \left\{  \left(\mu\left(
\theta\right)  -b\right)y\right\}}{y}dy
\label{psiGammaOneSide}
\end{equation}
and the desired $\left(  \psi,K\right)  $ is
\begin{equation}
\left\{
\begin{array}
[c]{l}%
K\left(  t,x\right)  =2^{-1}-K_{1}\left(  t,x\right)  -2^{-1}K_{3,0}\left(
t,x;\theta_{b}\right) \\
\psi\left(  t,\theta\right)  =2^{-1}-\psi_{1}\left(  t,\theta\right)
-2^{-1}\psi_{3,0}\left(  t,\theta;\theta_{b}\right)
\end{array}
\right.  . \label{V-b}%
\end{equation}

\end{theorem}

Note that the limit as $t \to \infty$ of $\psi_{1}\left(t,\theta\right)$ in (\ref{psiGammaOneSide}) is the Dirichlet integral in (\ref{EqDB1}) with $\mu$ as $\mu\left(\theta\right)$, which implies $\lim_{t \to \infty} \psi\left(  t,\theta\right) = 1_{\left(-\infty,b\right)}\left(\mu\right)$ for $\psi\left(  t,\theta\right)$ in (\ref{V-b}).
Recall $u_{3,m}\ $in (\ref{eq12e}) and define $\check{u}_{3,m}=\min_{\left\{
j:\theta_{j}\neq\theta_{b}\right\}  }\left\vert \xi\left(  \theta_{b}\right)
-\xi\left(  \theta_{j}\right)  \right\vert $. A suitable magnitude of
$\check{u}_{3,m}$ is needed for the estimator induced by $K_{3,0}\left(
t,x;\theta_{b}\right)  $ (already given by \autoref{ThmPoinNull}) in (\ref{V-b}) to consistently estimate the
proportion of $\mu\left(  \theta_{i}\right)  $'s that are equal to $\mu\left(
\theta_{b}\right)  $; see Theorem 9 of \cite{Chen:2018a}. We have the uniform
consistency of the proposed estimator $\hat{\varphi}_{m}\left(  t,\mathbf{z}\right)$ as

\begin{theorem}
\label{ThmVNEFConsistency}Assume $\left\{  z_{j}\right\}  _{j=1}^{m}$ are
independent Gamma random variables with parameters $\left\{  \left(
\theta_{i},\sigma\right)  \right\}  _{i=1}^{m}$. Then,
when $t$ is positive and sufficiently large,
\[
\mathbb{V}\left\{  \hat{\varphi}_{m}\left(  t,\mathbf{z}\right)  \right\} {\le} \frac{Ct^{11/4-\sigma}}{m^{2}u_{3,m}^{2}}%
\exp\left(  \frac{4t\max\left\{  1,\sqrt{2}\sigma\right\}  }{u_{3,m}}\right)
\sum_{i=1}^{m}\left(  1-\theta_{i}\right)^{\sigma-3/4}.
\]
Further, when $\sigma\geq 3/4$
\[
\mathcal{Q}\left(  \mathcal{F}\right)  =\left\{
\begin{array}
[c]{c}%
t_{m}=\left(  4\max\left\{  1,\sqrt{2}\sigma\right\}  \right)  ^{-1}%
u_{3,m}\gamma\ln m,t_{m}^{-1}\check{u}_{3,m}^{-1}=o\left(  \pi_{1,m}\right)
,\\
t_{m}\rightarrow\infty,u_{3,m}^{3/4-\sigma}\left(  \ln m\right)
^{11/4-\sigma}\left\Vert 1-\boldsymbol{\theta}\right\Vert _{\infty}%
^{\sigma-3/4}=o\left(  \pi_{1,m}^{2}m^{1-\gamma}\right)
\end{array}
\right\}
\]
is a uniform consistency class for each $\gamma\in\left(  0,1\right)  $, for which $\gamma=1$ can be set when $\sigma>11/4$, and when $\sigma \leq3/4$
\[
\mathcal{Q}\left(  \mathcal{F}\right)  =\left\{
\begin{array}
[c]{c}%
t_{m}=\left(  4\max\left\{  1,\sqrt{2}\sigma\right\}  \right)  ^{-1}%
u_{3,m}\gamma\ln m,t_{m}^{-1}\check{u}_{3,m}^{-1}=o\left(  \pi_{1,m}\right)
,\\
t_{m}\rightarrow\infty,\left(  \ln m\right)  ^{11/4-\sigma}=o\left(
\pi_{1,m}^{2}m^{1-\gamma}\right)
\end{array}
\right\}
\]
is a uniform consistency class for each $\gamma\in\left(  0,1\right)  $.
\end{theorem}

In terms of the uniform consistency class and maximal speed of convergence of the estimator $\hat{\varphi}_{m}\left(  t,\mathbf{z}\right)$,
\autoref{ThmVNEFConsistency} reveals similar things as does \autoref{ConcentrationIII}. In particular, the maximal speed of convergence $\hat{\varphi}_{m}\left(  t,\mathbf{z}\right)$ that maintains uniform consistency in the setting of one-side null does not depends on the
scale parameter $\sigma$ when $\sigma \le 2^{-1}\sqrt{2}$ since $\max\left\{  1,\sqrt{2}\sigma\right\}=1$ for such $\sigma$.
When $\theta=1/2$ and $\sigma$ is a positive, even integer or when $\sigma=1$, the corresponding
Gamma distribution becomes a central Chi-square distribution with degrees of
freedom $2^{-1}\sigma$ or the exponential distribution with mean $(1-\theta)^{-1}$. So, \autoref{ConcentrationIII} and \autoref{ThmVNEFConsistency} can be applied to
proportion estimation for central Chi-square and exponential random variables.

\section{Extension of constructions for a bounded null}

\label{SecExtensions}

We extend the previous constructions of the setting of a bounded null to the setting of estimating the
``induced proportion of true null hypotheses''
, i.e., to estimate
\begin{equation}
\check{\pi}_{0,m}=m^{-1}\sum\nolimits_{\left\{  i\in\left\{  1,\ldots
,m\right\}  :\mu_{i}\in\Theta_{0}\right\}  }\phi\left(  \mu_{i}\right)
\label{eq4c}%
\end{equation}
for a suitable function. In this setting, $\check{\pi}_{0,m}\in\left[
0,1\right]  $ does not necessarily hold. For example, $\phi(x)=\vert x
\vert^{p} $ for some $p>0$ gives $\check{\pi}_{0,m}$ as the ``average $l^{p}%
$-norm'' as a measure of sparsity for the vector $\boldsymbol{\mu}=\left(
\mu_{1},...,\mu_{m}\right)  ^{\top}$, and $\phi(x;c_{\ast})= \min\{\vert x
\vert^{p}, c_{\ast}^{p}\}$ for some $p,c_{\ast}>0$ gives $\check{\pi}_{0,m}$
as the proportion of ``signals'' exceeding a threshold $c_{\ast}$, which were
discussed by \cite{Jin:2008} and covered by his constructions. Such
$\check{\pi}_{0,m}$ has applications in genomics and signal processing as
mentioned in \autoref{sec1intro}.

To formulate the extension, we need a crucial auxiliary result that is different than those provided in \autoref{AuxRes3}. For a $\phi\in L^{1}\left(  \left[  a,b\right]  \right)  $ (with finite $a$ and
$b $ such that $a<b$), define
\begin{equation*}
\mathcal{D}_{\phi}\left(  t,\mu;a,b\right)  =\frac{1}{\pi}\int_{a}^{b}
\frac{\sin\left\{  \left(  \mu-y\right)  t\right\}  }{\mu-y}\phi\left(
y\right)  dy\text{ \ for \ }t,\mu\in\mathbb{R}.
\end{equation*}
Further, we introduce:
\begin{definition}
\label{defDiffratio}
{Let $\phi$ be a function defined on $U.$}  For each $\mu\in U  $ and $z\in\left[  \mu -b,\mu-a\right]  $, define
\begin{equation}
W_{\mu}\left(  z\right)  =z^{-1}\left(  \phi\left(  \mu-z\right)  -\phi\left(
\mu\right)  \right)  \text{ for }\ z\neq0\label{eqWmu}%
\end{equation}
and%
\begin{equation}
W_{\mu}^{-}\left(  z\right)  =\left\{
\begin{array}
[c]{lll}%
W_{\mu}\left(  z\right)   & \text{for} & \mu-a\geq z>0\\
\lim_{z\rightarrow0+}W_{\mu}\left(  z\right)   & \text{for} & z=0\text{ if the
limit exists}%
\end{array}
\right.  \label{eqWmuA}%
\end{equation}
and%
\begin{equation}
W_{\mu}^{+}\left(  z\right)  =\left\{
\begin{array}
[c]{lll}%
-W_{\mu}\left(  z\right)   & \text{for} & \mu-b\leq z<0\\
-\lim_{z\rightarrow0-}W_{\mu}\left(  z\right)   & \text{for} & z=0\text{ if
the limit exists}%
\end{array}
\right.  .\label{eqWmuB}%
\end{equation}
\end{definition}
{
Finally, for the rest of this section, we assume that $U$ is an interval
containing $\left[  a,b\right]  $ as a subset, and take the convention that
the total variation of a function defined on $U$ is the supremum of the total
variation of this function on each compact interval contained in $U$ (to avoid
unneeded complications of defining total variation of a function on a
measurable set).
}
Then we have:
\begin{lemma}
\label{ThmDirichletInterval}If $\phi\in L^{1}\left(  \left[  a,b\right]
\right)  $, then setting $\hat{\phi}\left(  s\right)  =\int_{a}^{b}\phi\left(
y\right)  \exp\left(  -\iota ys\right)  dy$ gives
\begin{equation}
\mathcal{D}_{\phi}\left(  t,\mu;a,b\right)  =\frac{t}{2\pi}\int_{-1}^{1}
\hat{\phi}\left(  ts\right)  \exp\left(  \iota\mu ts\right)  ds. \label{eq4d}%
\end{equation}
On the other hand, if $\phi$ is continuous and of bounded variation on
$\left[  a,b\right]  $, then
\begin{equation}
\mathcal{D}_{\phi,\infty}\left(  \mu;a,b\right) :=
\lim_{t\rightarrow\infty}\mathcal{D}_{\phi}\left(  t,\mu;a,b\right)  =\left\{
\begin{array}
[c]{lll}%
\phi\left(  \mu\right)  & \text{if} & a<\mu<b\\
2^{-1}\phi\left(  \mu\right)  & \text{if} & \mu=a\text{ or }\mu=b\\
0 & \text{if} & \mu<a\text{ or }\mu>b
\end{array}.
\right.  \label{eqA1}%
\end{equation}
If in addition both $W_{\mu}^{+}$ and $W_{\mu}^{-}$ are well-defined and of bounded variation for each fixed $\mu\in U$, then
\begin{equation}
\left\vert \mathcal{D}_{\phi}\left(  t,\mu;a,b\right)  -\mathcal{D}%
_{\phi,\infty}\left(  \mu;a,b\right)  \right\vert \leq\frac{4C_{\mu}\left(
\phi\right)  }{\pi t}+\frac{4\left\Vert \phi\right\Vert _{\infty}}{t}\left(
\frac{{2^{-1}}}{b-a}+\frac{{1}}{\delta_{\mu,a,b}}\right)  ,\label{eqA4}%
\end{equation}
when $\min\left\{  t\delta_{\mu,a,b},t\left(b-a\right)  \right\}  \geq2$, where
$\delta_{\mu,a,b}%
=\min_{\mu\notin\left\{  a,b\right\}  }\left\{ \left\vert \mu
-a\right\vert ,\left\vert \mu-b\right\vert \right\}  $ and
\begin{equation}
C_{\mu}\left(  \phi\right)  =\left\Vert W_{\mu}^{-}\right\Vert _{\infty
}+\left\Vert W_{\mu}^{-}\right\Vert _{\mathrm{TV}}+\left\Vert W_{\mu}%
^{+}\right\Vert _{\infty}+\left\Vert W_{\mu}^{+}\right\Vert _{\mathrm{TV}%
}.\label{defphiBnd}%
\end{equation}

\end{lemma}

\autoref{ThmDirichletInterval} gives the relationship between the Fourier transform of $\phi$ and how to invert its Fourier transform, and will be used to derive \autoref{ThmExtA} below. Note that (\ref{eq4d}) is almost the inverse of the Fourier transform of
$\phi$. We caution that (\ref{eqA1}) does not necessarily hold when $\phi$ is
only continuous, as can be seen from the examples in Chapter VIII of
\cite{Zygmund1959}.
Note also that (\ref{eqA4}) gives the speed of convergence of $\mathcal{D}%
_{\phi}\left(  t,\mu;a,b\right)  $ and helps determine the speed of
convergence of the estimators to be constructed below:

\begin{theorem}
\label{ThmExtA}Let $\phi$ be continuous and of bounded variation on $\left[
a,b\right]  $. Assume $\mathcal{F}$ is the Gamma family, and set
\begin{equation}
K_{1}\left(  t,x\right)  =\frac{t}{2\pi}\int_{a}^{b}\phi\left(  y\right)
dy\int_{-1}^{1}\sum_{n=0}^{\infty}\frac{\left(  tsx\tilde{a}_{1}\right)
^{n}\cos\left(  2^{-1}n\pi-tsy\right)  }{\tilde{a}_{n}n!}ds. \label{eq13c}%
\end{equation}
Then
\begin{equation}
\psi_{1}\left(  t,\theta\right)  ={\int}K_{1}\left(  t,x\right)  dG_{\theta
}\left(  x\right)
= \mathcal{D}_{\phi}\left(  t,\mu\left(\theta\right);a,b\right),
\label{eq13ca}%
\end{equation}
and the desired $\left(  \psi,K\right)  $ for estimating $\check{\pi}_{0,m}$
is%
\begin{equation}
\left\{
\begin{array}
[c]{l}%
K\left(  t,x\right)  =K_{1}\left(  t,x\right)  -2^{-1}\left\{  \phi\left(
a\right)  K_{3,0}\left(  t,x;\theta_{a}\right)  +\phi\left(  b\right)
K_{3,0}\left(  t,x;\theta_{b}\right)  \right\}  \\
\psi\left(  t,\mu\right)  =\psi_{1}\left(  t,\mu\right)  -2^{-1}\left\{
\phi\left(  a\right)  \psi_{3,0}\left(  t,\theta;\theta_{a}\right)
+\phi\left(  b\right)  \psi_{3,0}\left(  t,\theta;\theta_{b}\right)  \right\}
\end{array}
\right.  .\label{eq2c}%
\end{equation}

\end{theorem}

Note that $\lim_{t \to \infty} \psi\left(  t,\mu\right) = 1_{\left(a,b\right)}\left(\mu\right)  \phi \left(\mu\right)$
for $\psi\left(  t,\mu\right)$ in (\ref{eq2c}).
The constructions in \autoref{ThmExtA} can be easily modified to estimate any
linear function of $\check{\pi}_{0,m}$, which will not be discussed here.
Moreover, when $\phi$ is the identity function on $\left[
a,b\right]  $, (\ref{eq13c}) reduces to (\ref{eq13g}).

Define%
\begin{equation}
\hat{\varphi}_{m}\left(  t,\mathbf{z}\right)  =m^{-1}\sum_{i=1}^{m}K\left(
t,z_{i}\right)  \text{ \ and \ }\varphi_{m}\left(  t,\boldsymbol{\mu}\right)
=m^{-1}\sum_{i=1}^{m}\mathbb{E}\left\{  K\left(  t,z_{i}\right)  \right\}
\label{eq21b}%
\end{equation}
with $K$ in (\ref{eq2c}) and set $e_{m}\left(  t\right)
=\hat{\varphi}_{m}\left(  t,\mathbf{z}\right)  -\varphi_{m}\left(
t,\boldsymbol{\mu}\right)  $. Then $\hat{\varphi}_{m}\left(  t,\mathbf{z}
\right)  $ defined by (\ref{eq21b}) estimates $\check{\pi}_{0,m}$ defined by (\ref{eq4c}). Consistency of the estimator
$\hat{\varphi}_{m}\left(  t,\mathbf{z}\right)  $ given by (\ref{eq21b}) can be
obtained for independent $\left\{  z_{i}\right\}  _{i=1}^{m}$ via almost
identical arguments as those for the proof of \autoref{ConcentrationIII}; see
\autoref{ThmFinal} below. For the rest of this
section, we assume that $\phi$ is continuous and of bounded variation on
$\left[  a,b\right]  $. Recall $u_{3,m}$ and
$\tilde{u}_{3,m}$ in (\ref{eq12e}).

Recall $C_{\mu}\left(  \phi\right)$ defined by (\ref{defphiBnd}). We
introduce the following condition:
\begin{description}
\item[C3)] $\phi$ has finite left and right derivatives at each {interior point of $U$}, a finite right derivative at {the left boundary point of $U$}, and a finite left derivative at {the right boundary point of $U$}, and $\left\Vert\phi\right\Vert _{1,\infty}:= \sup_{\mu\in {U} } C_{\mu}\left(  \phi\right) < \infty$.
\end{description}
We remark that the first part of condition C3) implies that $W_{\mu}^{+}$ and $W_{\mu}^{-}$ in \autoref{defDiffratio} are well-defined, that
the second part of condition C3), i.e., $\left \Vert \phi \right \Vert_{1,\infty} < \infty$, implies $W_{\mu}^{+}$ and $W_{\mu}^{-}$ are of bounded variation {on $U$} when they are well-defined,
and that condition C3) implies $\phi$ is Lipschitz on {$U$} and that its derivative $\phi^{\prime}$, if defined everywhere on {$U$}, is of bounded variation on {$U$}.
The two examples of $\phi$ mentioned at the beginning of this section satisfy condition C3) when their exponent $p \ge 1$ {and $U$ is a bounded interval}.

We have the uniform consistency of $\hat{\varphi}_{m}\left(  t,\mathbf{z}\right)$ as:

\begin{theorem}
\label{ThmFinal}Consider the estimator $\hat{\varphi}_{m}\left(
t,\mathbf{z}\right)  $ in (\ref{eq21b}) where $\phi$ is continuous and of bounded variation on $\left[a,b\right]$, and assume $\left\{  z_{i}\right\}
_{i=1}^{m}$ are independent whose CDFs are members of $\mathcal{F}$.
If $\mathcal{F}$ is the Gamma family, then for all positive and large $t$
\[
\mathbb{V}\left\{  e_{m}\left(  t\right)  \right\}  \leq\frac{C\left\Vert
\phi\right\Vert _{\infty}^{2}\left(  1+t^{2}\right)  }{m^{2}}\exp\left(
\frac{4t\max\left\{  \sigma,1\right\}  }{u_{3,m}}\right)  \sum_{i=1}%
^{m}\left(  \frac{t}{1-\theta_{i}}\right)  ^{3/4-\sigma};
\]
if in addition condition C3) holds, then (\ref{eq21}) is a uniform consistency class for each $\gamma
\in\left(  0,1\right)  $ when $\sigma\geq 3/4$, for which $\gamma=1$ can be set when $\sigma > 11/4$ and (\ref{eq21a}) a uniform
consistency class for each $\gamma\in\left(  0,1\right)  $ when $\sigma
\leq3/4$, after replacing $\pi_{1,m}$ in (\ref{eq21}) and (\ref{eq21a}) by
$\check{\pi}_{0,m}$.
\end{theorem}

\section{Discussion}

\label{SecConcAndDisc}

For multiple testing a bounded or one-sided null on the means of
random variables whose CDFs are members of
the Gamma family, we have constructed uniformly consistent estimators of the
corresponding proportion of false null hypotheses via solutions to
Lebesgue-Stieltjes integral equations, for which
consistency, speeds of convergence, and uniform consistency classes have been obtained under independence between these random variables. The strategy proposed in the Discussion
section of \cite{Chen:2018a} or that in Section 2.3 of \cite{Jin:2008}, i.e., choosing a speed of convergence $t_{m}$ that controls the variance of the error term $e_{m}(t_m)$ of a proportion estimator for a finite $m$, can be used to adaptively determine the speed of
convergence (and hence the tuning parameter $\gamma$) for the proposed
estimators. These estimators can be used to develop adaptive versions of the
FDR procedure of \cite{Chen:2018d}, the ``BH'' procedure of
\cite{Benjamini:2001} under the conditions of their Theorem 5.2, the FDR and
FNR procedures of \cite{Sarkar:2006}, or any other conservative FDR or FNR
procedure that is applicable to multiple testing composite null hypotheses.

The constructions and uniform consistency of the proportion estimators
provided here can be easily extended to the setting where the null parameter
set belongs to the algebra generated by bounded, one-sided and point nulls.
Here the term ``algebra'' refers to the
family of sets generated by applying any finite combination of set union,
intersection or complement to these three types of nulls.
This covers null sets $\left(-\infty,b\right]$ and $\left[a,b\right]$ that are more conventional in classic textbooks such as \cite{Lehmann:2005},
even though we have chosen null sets $\left(-\infty,b\right)$ and $\left(a,b\right)$ mainly due to how Dirichlet integrals given in \autoref{secResultsPoint}
converge to their discontinuous, piecewise linear limiting functions.
{Details on how to adapt the constructions and estimators to the null sets $\left(-\infty,b\right]$ and $\left[a,b\right]$
are provided in 
Section 5 of the supplementary material.}

Further, it is possible to
establish the uniform consistency of the proportion estimators in
\autoref{secNEF} and \autoref{SecExtensions} for weakly dependent random
variables that are bivariate bivariate Gamma via their associated Laguerre polynomials.
{However, we were not able to obtain better concentration inequalities for the error of the estimators in this work, and an explanation
on this is provided in \autoref{NotBetterConcentration}. }
Finally, our constructions essentially utilize suitable group structures of the domain of the parameters and arguments of
the CDFs and then apply special transforms that are adapted to such group structures to obtain solutions
to the Lebesgue-Stieltjes integral equation. This general principle is also applicable to probability distributions on Lie groups (which contains, e.g., the unit sphere a special case), and tools of harmonic analysis on Lie groups can be used to construction proportion estimators for these distributions. This has applications
to modeling and analysis of data on Lie groups. We will report on all these in other articles.

\appendix{}

\section{Auxiliary results and their proofs}

\label{AppProofA}

{We will provide in \autoref{NEFGamma} a classification result on NEFs with separable moments,
in \autoref{AuxRes3} a lemma on Dirichlet integrals and two lemmas on the speed of convergence of a Fourier transform,
in \autoref{ProofOfLemma1} a proof of \autoref{ThmDirichletInterval},
in \autoref{OracleSpeeds} speeds of convergence for discriminant functions,
and in \autoref{NotBetterConcentration} an explanation
on not being able to obtain better concentration inequalities.}

\subsection{NEFs with separable moments}
\label{NEFGamma}

{We will show that the only NEF with separable moments is the Gamma family.}
We provide a very brief review on natural exponential family (NEF), whose details can be found in
\cite{Letac:1992}. Let $\beta$ be a positive Radon measure on $\mathbb{R}$
that is not concentrated on one point. Let $L\left(  \theta\right)  =\int
e^{x\theta}\beta\left(  dx\right)  $ for $\theta\in\mathbb{R}$ be its Laplace
transform and $\Theta$ be the maximal open set containing $\theta$ such that
$L\left(  \theta\right)  <\infty$. Suppose $\Theta$ is not empty and let
$\kappa\left(  \theta\right)  =\ln L\left(  \theta\right)  $ be the cumulant
function of $\beta$. Then%
\[
\mathcal{F}=\left\{  G_{\theta}:G_{\theta}\left(  dx\right)  =\exp\left\{
\theta x-\kappa\left(  \theta\right)  \right\}  \beta\left(  dx\right)
,\theta\in\Theta\right\}
\]
forms an NEF with respect to the basis $\beta$. Note that $\Theta$ has a
non-empty interior and is a convex set if it is not empty and that $L$ is analytic on the strip
$A_{\Theta}=\left\{  z\in\mathbb{C}:\Re\left(  z\right)  \in\Theta\right\}  $.

The NEF $\mathcal{F}$ can be equivalently characterized by its mean domain and
variance function. Specifically, the mean function $\mu:\Theta\rightarrow U$
with $U=\mu\left(  \Theta\right)  $ is given by $\mu\left(  \theta\right)
=\frac{d}{d\theta}\kappa\left(  \theta\right)  $, and the variance function is
$V\left(  \theta\right)  =\frac{d^{2}}{d\theta^{2}}\kappa\left(
\theta\right)  $ and can be parametrized by $\mu$ as
\[
V\left(  \mu\right)  =\int\left(  x-\mu\right)  ^{2}F_{\mu}\left(  dx\right)
\text{ for }\mu\in U,
\]
where $\theta=\theta\left(  \mu\right)  $ is the inverse function of $\mu$ and
$F_{\mu}=G_{\theta\left(  \mu\right)  }$. Namely, $\mathcal{F}=\left\{
F_{\mu}:\mu\in U\right\}  $. The pair $\left(  V,U\right)  $ is called the
variance function of $\mathcal{F}$, and it characterizes $\mathcal{F}$ and is analytic on the mean domain $U$.

For the case of a point null when $\mathcal{F}$ is an NEF, \cite{Chen:2018a}
introduced the definition of ``NEF with separable moment
functions at a specific point'', showed that the Gamma family
has such moment functions, and constructed a uniformly consistent proportion
estimator for the Gamma family. In order to construct proportion estimators
for a one-sided or bounded null or a functional of a bounded null when $\mathcal{F}$ is an NEF, we need a
stronger definition as

\begin{definition}
\label{Def3}For each $\theta\in\Theta$ and $G_{\theta}\in\mathcal{F}$ for an
NEF $\mathcal{F}$, denote the moment sequence of $G_{\theta}$ by%
\[
\tilde{c}_{n}\left(  \theta\right)  =\int x^{n}G_{\theta}\left(  dx\right)
\text{ \ for \ }n=0,1,\ldots
\]
If there exist two functions $\zeta,\xi:\Theta\rightarrow\mathbb{R\ }$and a
sequence of real numbers $\left\{  \tilde{a}_{n}\right\}  _{n\geq0}$ that
satisfy the following:

\begin{itemize}
\item $\xi$ is one-to-one, $\zeta\left(  \theta\right)  \neq0$ for all
$\theta\in U$, and $\zeta$ does not depend on any $n\in\mathbb{N}$,

\item $\tilde{c}_{n}\left(  \theta\right)  =\xi^{n}\left(  \theta\right)
\zeta\left(  \theta\right)  \tilde{a}_{n}$ for each $n\in\mathbb{N}$ and
$\theta\in\Theta$,

\item $\Psi\left(  t,\theta\right)  =\sum_{n=0}^{\infty}\frac{t^{n}\xi
^{n}\left(  \theta\right)  }{\tilde{a}_{n}n!}$ is absolutely convergent
pointwise in $\left(  t,\theta\right)  \in\mathbb{R}\times\Theta$,
\end{itemize}

\noindent then the moment sequence $\left\{  \tilde{c}_{n}\left(
\theta\right)  \right\}  _{n\geq0}$ is called
``separable'' and $\mathcal{F}$ is said to have ``separable moments''.
\end{definition}

\autoref{Def3} requires that the moments of $G_{\theta}$ are separable in the
sense of $\tilde{c}_{n}\left(  \theta\right)  =\xi^{n}\left(  \theta\right)
\zeta\left(  \theta\right)  \tilde{a}_{n}$ at all points in $\Theta$ rather
than at just one point in $\Theta$. Note that $\mu\left(  \theta\right)
=\xi\left(  \theta\right)  \zeta\left(  \theta\right)  \tilde{a}_{1}$ for an
NEF with separable moments. In the main text, we have restated that the Gamma family has separable moments.
If there are NEFs with separable moments different than the Gamma family, then
the proportion estimators of \cite{Chen:2018a} and those proposed here
will have wider applications. However, we have a somewhat disappointing result
on how many NEFs with separable moments exist as

\begin{theorem}
\label{MainRes}The only NEF with separable moments is the Gamma family.
\end{theorem}

\begin{proof}
By definition,
\[
\tilde{c}_{n}\left(  \theta\right)  =\frac{1}{L\left(  \theta\right)  }\int
x^{n}e^{\theta x}\beta\left(  dx\right)  =\int x^{n}G_{\theta}\left(
dx\right)  \text{ \ for \ }n=0,1,\ldots
\]
By the assumption $\tilde{c}_{n}\left(  \theta\right)  =\xi^{n}\left(
\theta\right)  \zeta\left(  \theta\right)  \tilde{a}_{n}$ for $n\in\mathbb{N}%
$, we have $\mu\left(  \theta\right)  =\tilde{a}_{1}\xi\left(  \theta\right)
\zeta\left(  \theta\right)  $ and $\tilde{a}_{1}\neq0$. However, $V\left(
\theta\right)  =\frac{d}{d\theta}\mu\left(  \theta\right)  $ and
\begin{equation}
V\left(  \theta\right)  =\tilde{a}_{1}\frac{d}{d\theta}\left\{  \xi\left(
\theta\right)  \zeta\left(  \theta\right)  \right\}  =\xi^{2}\left(
\theta\right)  \zeta\left(  \theta\right)  \tilde{a}_{2}-\xi^{2}\left(
\theta\right)  \zeta^{2}\left(  \theta\right)  \tilde{a}_{1}^{2}.
\label{eq2bx}%
\end{equation}
Since $V\left(  \theta\right)  >0$ and $\tilde{c}_{2}\left(  \theta\right)
>0$, we have $\zeta\left(  \theta\right)  \tilde{a}_{2}>0$ and from
(\ref{eq2bx}) we obtain
\begin{equation}
V\left(  \theta\right)  =\mu^{2}\left(  \theta\right)  \left(  \frac{\tilde
{a}_{2}}{\zeta\left(  \theta\right)  \tilde{a}_{1}^{2}}-1\right)  \text{
\ with \ }\frac{\tilde{a}_{2}}{\zeta\left(  \theta\right)  \tilde{a}_{1}^{2}%
}>1. \label{eq3x}%
\end{equation}
Write $\tilde{h}\left(  \mu\right)  =\zeta\left(  \theta\right)  $ and $\tilde{V}\left(\mu\right)=V\left(\theta\right)$. Then%
\begin{equation}
\tilde{V}\left(  \mu\right)  =\mu^{2}\left(  \frac{\tilde{a}_{2}}{\tilde{h}\left(
\mu\right)  \tilde{a}_{1}^{2}}-1\right)  \text{ \ with \ }\frac{\tilde{a}_{2}%
}{\tilde{h}\left(  \mu\right)  \tilde{a}_{1}^{2}}>1 \label{eq3ax}%
\end{equation}
and $\tilde{V}$, the variance function of the mean $\mu$, is real, positive, and analytic on the mean domain $U$ by Proposition 2.2. of \cite{Letac:1990}.

From (\ref{eq3ax}), we see that $\tilde{V}\left(  \mu\right)  $ is
quadratic in $\mu$ if and only if $\zeta\left(  \theta\right)  $ is constant.
However, setting $n=0$ for $\tilde{c}_{n}\left(  \theta\right)  $, we obtain%
\[
\tilde{c}_{0}\left(  \theta\right)  =\frac{1}{L\left(  \theta\right)  }\int
e^{\theta x}\beta\left(  dx\right)  =1=\zeta\left(  \theta\right)  \tilde
{a}_{0}.
\]
So, $\tilde{a}_{0}\neq0$ and $\zeta\equiv\tilde{a}_{0}^{-1}$. Thus, $\tilde{V}\left(
\mu\right)  $ is quadratic in $\mu$, and $\tilde{V}$ has a double root at zero and no
other roots. By \cite{Letac:1990} and their Table 1, we conclude that the
family $\mathcal{F}$ has to be the Gamma family. \qed
\end{proof}


\subsection{Dirichlet integral and Fourier transform}
\label{AuxRes3}

In order to construct the new proportion estimators and show their
uniform consistency for the settings of a one-sided null, a bounded null and a suitable functional of a bounded null
respectively, we introduce three auxiliary results that are of independent interest. First, we have the speed of convergence of Dirichlet integral as

\begin{lemma}
\label{lm:Dirichlet}$\left\vert \int_{0}^{t}x^{-1}\sin xdx-2^{-1}%
\pi\right\vert \leq2\pi t^{-1}$ for $t\geq2$.
\end{lemma}

\autoref{lm:Dirichlet} implies the {identities (\ref{EqDirichlet1}) and (\ref{EqDB1}) that have been presented in the main text.}


Further, we have the following identities that will be used in the integral
representations of solutions to the Lebesgue-Stieltjes integral equation
(\ref{eq3}) in the constructions of the new proportion estimators:

\begin{lemma}
\label{LmDirichlet}For any $a,b,\mu,t\in\mathbb{R}$ with $a<b$,
\begin{equation}
\int_{\left(  \mu-b\right)  t}^{\left(  \mu-a\right)  t}\frac{\sin v}
{v}dv=\frac{1}{2}\int_{a}^{b}dy\int_{-1}^{1}t\exp\left\{  \iota\left(
\mu-y\right)  ts\right\}  ds. \label{eq1c}%
\end{equation}
On the other hand, for any $b,\mu,t\in\mathbb{R}$,
\begin{equation}
\frac{1}{\pi}\int_{0}^{t}\frac{\sin\left\{  \left(  \mu-b\right)  y\right\}
}{y}dy=\frac{1}{2\pi}\int_{0}^{t}dy\int_{-1}^{1}\left(  \mu-b\right)
\exp\left\{  \iota ys\left(  \mu-b\right)  \right\}  ds. \label{eq1c1}%
\end{equation}

\end{lemma}

Finally, we have the following to be used to derive the speed of the convergence
of $\varphi_m\left(t,\boldsymbol{\mu}\right)$:

\begin{lemma}
\label{lm:OracleSpeed}Let $-\infty<a_{1}<b_{1}<\infty$. If $f:\left[
a_{1},b_{1}\right]  \rightarrow\mathbb{R}$ is of bounded variation, then%
\begin{equation}
\left\vert \int_{\left[  a_{1},b_{1}\right]  }f\left(  s\right)  \cos\left(
ts\right)  ds\right\vert \leq4\left(  \left\Vert f\right\Vert _{\mathrm{TV}%
}+\left\Vert f\right\Vert _{\infty}\right)  \left\vert t\right\vert
^{-1}1_{\left\{  t\neq0\right\}  }\left(  t\right)     \label{eqRateFourier}%
\end{equation}
and
\begin{equation}
\left\vert \int_{\left[  a_{1},b_{1}\right]  }f\left(  s\right)  \sin\left(
ts\right)  ds\right\vert \leq4\left(  \left\Vert f\right\Vert _{\mathrm{TV}%
}+\left\Vert f\right\Vert _{\infty}\right)  \left\vert t\right\vert
^{-1}1_{\left\{  t\neq0\right\}  }\left(  t\right)  . \label{eqRateFourier1}%
\end{equation}
\end{lemma}

Note that better bounds than (\ref{eqRateFourier})\ can be derived when $f$ has higher
order derivatives (by adapting the techniques of \cite{Jackson1920}) but are
not our focus here and will not improve the speeds of convergence of the
proportion estimators to be presented later.

The proofs of these three lemmas are given below in order.

\subsubsection{Proof of \autoref{lm:Dirichlet}}

Pick $\rho \in \left(0,1\right)$ and $R$ such that $R>\rho$. Define the counterclockwise oriented
contour $\tilde{C}=\bigcup\nolimits_{i=1}^{4}C_{i}$, where $C_{1}=\left\{
Re^{\iota x}:0\leq x\leq\pi\right\}  $, $C_{2}=\left\{  \rho e^{\iota x}%
:\pi\geq x\geq0\right\}  $, $C_{3}=\left\{  x:-R\leq x\leq-\rho\right\}  $ and
$C_{4}=\left\{  x:\rho\leq x\leq R\right\}  $. Let $f\left(  z\right)
=e^{\iota z}/z$ for $z\in\mathbb{C\setminus}\left\{  0\right\}  $. Then, by the residue theorem (see, e.g., \cite{Ahlfors:1979}),%
\[
0=\int_{\tilde{C}}f\left(  z\right)  dz=\left(  \int_{C_{1}}+\int_{C_{2}}%
+\int_{C_{3}}+\int_{C_{4}}\right)  f\left(  z\right)  dz.
\]
So,%
\[
-\int_{C_{1}}f\left(  z\right)  dz-\int_{C_{2}}f\left(  z\right)  dz-\iota
\pi=\left(  \int_{C_{3}}+\int_{C_{4}}\right)  f\left(  z\right)  dz-\iota\pi.
\]
However,%
\[
\left\vert \int_{C_{1}}f\left(  z\right)  dz\right\vert \leq2\int_{0}^{\pi
/2}e^{-R\sin x}dx\leq2\int_{0}^{\pi/2}\exp\left(  -2R\pi^{-1}x\right)
dx=\frac{\pi}{R}\left\{  1-\exp\left(  -R\right)  \right\},
\]
where the second inequality in the above display is obtained by using the inequality $\sin x \ge 2 x/\pi $ for $x \in \left[0,\pi/2\right]$ that is due to the concavity of $\sin x$ for $x \in \left[0,\pi/2\right]$,
and%
\[
\left\vert -\int_{C_{2}}f\left(  z\right)  dz-\iota\pi\right\vert \le\pi
\max_{0\leq x\leq\pi}\left\vert \exp\left(  \iota\rho e^{\iota x}\right)
-1\right\vert \leq\pi\frac{\rho}{1-\rho},%
\]
where the last inequality in the above display is obtained by first expanding $\exp\left(  \iota\rho e^{\iota x}\right) -1$ at $0$ into a power series and then bounding the power series by the geometric series $\sum_{i=1}^{\infty} \rho^i$ for $\rho \in \left(0,1\right)$,
and%
\[
\left(  \int_{C_{3}}+\int_{C_{4}}\right)  f\left(  z\right)  dz=2\iota
\int_{\rho}^{R}\frac{\sin x}{x}dx.
\]
Since the ratio $\rho\left(  1-\rho\right)  ^{-1}$ for $\rho\in\left(
0,1\right)  $ upper bounded by $2\rho$ when $\rho\le 2^{-1}$, and $\sin x\leq x$
for all $x\geq0$, setting $\rho=R^{-1}$ gives%
\[
\left\vert \int_{R^{-1}}^{R}x^{-1}\sin xdx-2^{-1}\pi\right\vert \leq \frac{3\pi}{2}
R^{-1}\text{ \ and }\left\vert \int_{0}^{R^{-1}}x^{-1}\sin xdx\right\vert \leq
R^{-1}%
\]
for all $R\geq2$. Thus, $\left\vert \int_{0}^{R}x^{-1}\sin xdx-2^{-1}%
\pi\right\vert \leq2\pi R^{-1}$ for all $R\geq2$. \qed

\subsubsection{Proof of \autoref{LmDirichlet}}

By simple algebra, we have%
\begin{align*}
&  \int_{\left(  \mu-b\right)  t}^{\left(  \mu-a\right)  t}\frac{\sin v}%
{v}dv=\int_{a}^{b}\frac{\sin\left\{  \left(  \mu-y\right)  t\right\}  }{\mu
-y}dy\\
&  =\int_{a}^{b}\frac{\exp\left\{  \iota\left(  \mu-y\right)  t\right\}
-\exp\left\{  -\iota\left(  \mu-y\right)  t\right\}  }{2\iota\left(
\mu-y\right)  }dy\\
&  =\frac{1}{2}\int_{a}^{b}dy\int_{-t}^{t}\exp\left\{  \iota\left(
\mu-y\right)  s\right\}  ds=\frac{1}{2}\int_{a}^{b}dy\int_{-1}^{1}%
t\exp\left\{  \iota\left(  \mu-y\right)  ts\right\}  ds.
\end{align*}
On the other hand,%
\begin{align*}
\frac{1}{\pi}\int_{0}^{t}\frac{\sin\left(  \mu y\right)  }{y}dy  &  =\frac
{1}{2\pi}\int_{0}^{t}\frac{2\iota\sin\left(  \mu y\right)  }{\iota y}%
dy=\frac{1}{2\pi}\int_{0}^{t}\frac{\exp\left(  \iota\mu y\right)  -\exp\left(
-\iota\mu y\right)  }{\iota y}dy\\
&  =\frac{1}{2\pi}\int_{0}^{t}dy\int_{-\mu}^{\mu}\exp\left(  \iota ys\right)
ds=\frac{1}{2\pi}\int_{0}^{t}dy\int_{-1}^{1}\mu\exp\left(  \iota y\mu
s\right)  ds.
\end{align*}
So, replacing $\mu$ by $\mu-b$ in the above identity, we have
the claimed identity. \qed

\subsubsection{Proof of \autoref{lm:OracleSpeed}}

By Jordan's decomposition theorem, $f=g_{1}-g_{2}$, where $g_{1}\left(
x\right)  =2^{-1}g_{0}\left(  x\right)  +2^{-1}f\left(  x\right)  $ and
$g_{2}\left(  x\right)  =2^{-1}g_{0}\left(  x\right)  -2^{-1}f\left(
x\right)  $ are non-decreasing functions on $\left[  a_{1},b_{1}\right]  $ and
$g_{0}\left(  x\right)  $ is the total variation of $f$ on $\left[
a_{1},x\right]  $ for $x\in\left[  a_{1},b_{1}\right]  $. So,%
\begin{equation}
I_{0}=\int_{\left[  a_{1},b_{1}\right]  }f\left(  s\right)  \cos\left(
ts\right)  ds=\int_{\left[  a_{1},b_{1}\right]  }g_{1}\left(  s\right)
\cos\left(  ts\right)  ds-\int_{\left[  a_{1},b_{1}\right]  }g_{2}\left(
s\right)  \cos\left(  ts\right)  ds. \label{eq1a}%
\end{equation}
For the first summand in (\ref{eq1a}), we can apply the second law of the mean
(see, e.g., \cite{Hobson:1909}, \cite{Dixon:1929} or Section 10 of Chapter 1 of \cite{Widder:1946})
to obtain%
\[
I_{g_{1}}=\int_{\left[  a_{1},b_{1}\right]  }g_{1}\left(  s\right)
\cos\left(  ts\right)  ds=g_{1}\left(  a_{1}\right)  \int_{a_{1}}^{s_{0}}%
\cos\left(  ts\right)  ds+g_{1}\left(  b_{1}\right)  \int_{s_{0}}^{b_{1}}%
\cos\left(  ts\right)  ds
\]
for some $s_{0}\in\left[  a_{1},b_{1}\right]  $.
So, when $t\neq0$,
\[
\left\vert g_{1}\left(  a_{1}\right)  \int_{a_{1}}^{s_{0}}\cos\left(
ts\right)  ds\right\vert =\left\vert g_{1}\left(  a_{1}\right)  \right\vert
\frac{\left\vert \sin\left(  s_{0}t\right)  -\sin\left(  ta_{1}\right)
\right\vert }{\left\vert t\right\vert }\leq2\left\vert t\right\vert
^{-1}\left\vert g_{1}\left(  a_{1}\right)  \right\vert
\]
and $\left\vert I_{g_{1}}\right\vert \leq 4 \left\Vert g_{1}\right\Vert
_{\infty}\left\vert t\right\vert ^{-1}$. Applying the same arguments to the second summand $I_{g_{2}}$ in
(\ref{eq1a}) yields
$\left\vert I_{g_{2}}\right\vert \leq 4 \left\Vert
g_{2}\right\Vert _{\infty}\left\vert t\right\vert ^{-1}$ when $t\neq0$.
However,
\[
\max\left\{  \left\Vert g_{1}\right\Vert _{\infty},\left\Vert g_{2}\right\Vert
_{\infty}\right\}  \leq2^{-1}\left\Vert f\right\Vert _{\mathrm{TV}}%
+2^{-1}\left\Vert f\right\Vert _{\infty},
\]
and
\[
\left\vert \int_{\left[  a_{1},b_{1}\right]  }f\left(  s\right)  \cos\left(
ts\right)  ds\right\vert \leq\left\Vert f\right\Vert _{\infty}\left(
b_{1}-a_{1}\right)  \leq\left(  b_{1}-a_{1}\right)  \left(  \left\Vert
f\right\Vert _{\mathrm{TV}}+\left\Vert f\right\Vert _{\infty}\right)  \text{
\ for }t=0\text{.}%
\]
Thus,%
\[
\left\vert I_{0}\right\vert \leq4\left(  \left\Vert f\right\Vert
_{\mathrm{TV}}+\left\Vert f\right\Vert _{\infty}\right)  \left\vert
t\right\vert ^{-1}1_{\left\{  t\neq0\right\}  }\left(  t\right)  +\left(
b_{1}-a_{1}\right)  \left\Vert f\right\Vert _{\infty}1_{\left\{  t=0\right\}
}\left(  t\right)  ,
\]
which is (\ref{eqRateFourier}).
Replacing $\cos\left(ts\right)$ by $\sin\left(ts\right)$ in (\ref{eq1a}) and using identical strategy given above, we can prove (\ref{eqRateFourier1}).\qed

\subsection{Proof of \autoref{ThmDirichletInterval}}
\label{ProofOfLemma1}

Recall%
\begin{equation}
\mathcal{D}_{\phi}\left(  t,\mu;a,b\right)  =\frac{1}{\pi}\int_{a}^{b}%
\frac{\sin\left\{  \left(  \mu-y\right)  t\right\}  }{\mu-y}\phi\left(
y\right)  dy=\frac{1}{\pi}\int_{\mu-b}^{\mu-a}\frac{\sin\left(  tz\right)
}{z}\phi\left(  \mu-z\right)  dz.\label{eqDphi}%
\end{equation}
Firstly,%
\begin{align*}
\mathcal{D}_{\phi}\left(  t,\mu;a,b\right)   &  =\frac{1}{2\pi}\int_{a}%
^{b}\frac{\exp\left\{  \iota\left(  \mu-y\right)  t\right\}  -\exp\left\{
{-}\iota\left(  \mu-y\right)  t\right\}  }{\iota\left(  \mu-y\right)  }%
\phi\left(  y\right)  dy\\
&  =\frac{1}{2\pi}\int_{a}^{b}\phi\left(  y\right)  dy\int_{-t}^{t}%
\exp\left\{  \iota\left(  \mu-y\right)  s\right\}  ds\\
&  = \frac{t}{2\pi}\int_{a}^{b}\phi\left(  y\right)   dy\int_{-1}^{1}\exp\left(  -\iota
yts\right) \exp\left(  \iota\mu {ts}\right)  ds.
\end{align*}
Namely, setting $\hat{\phi}\left(  s\right)  =\int_{a}^{b}\phi\left(
y\right)  \exp\left(  -\iota ys\right)  dy$ yields%
\[
\mathcal{D}_{\phi}\left(  t,\mu;a,b\right)  =\frac{t}{2\pi}\int_{-1}^{1}%
\hat{\phi}\left(  ts\right)  \exp\left(  \iota\mu ts\right)  ds.
\]

The second claim%
\[
\lim_{t\rightarrow\infty}\mathcal{D}_{\phi}\left(  t,\mu;a,b\right)  =\left\{
\begin{array}
[c]{lll}%
\phi\left(  \mu\right)   & \text{if} & a<\mu<b\\
2^{-1}\phi\left(  \mu\right)   & \text{if} & \mu=a\text{ or }\mu=b\\
0 & \text{if} & \mu<a\text{ or }\mu>b
\end{array}
\right.
\]
is a direct consequence of Theorem 7.1 and Theorem 7.2 in Chapter 2 of
\cite{Widder:1946} which contain the key result%
\[
\lim_{t\rightarrow\infty}\frac{1}{\pi}\int_{0}^{\delta}\phi\left(  y\right)
\frac{\sin\left(  ty\right)  }{y}dy=\frac{\phi\left(  0\right)  }{2}%
\]
if $\phi$ is continuous and of bounded variation on $\left[  0,\delta\right]
$ for a fixed $\delta>0$.
Please note the following 3 things though:\ (1) in the
proof of Theorem 7.1 of \cite{Widder:1946}, the inequality $\overline{\lim
}_{T\rightarrow\infty}\left\vert \frac{1}{\pi}\int_{0}^{\delta}\alpha\left(
t\right)  \frac{\sin Tt}{t}dt\right\vert \leq\varepsilon A$ appears there but
the lower limit of this integral should be $\eta$; (2) this $\eta$ can be
chosen to be positive; (3) none of (1) or (2) affects the arguments of the
proof or the conclusion of this Theorem 7.1.

{
For completeness, we provide details to proving the second claim. There will be 3 cases overall: (a)
$\mu=a$ or $\mu=b$; (b) $\mu\in\left(  a,b\right)  $; (c) $\mu<a$ or $\mu>b$,
as shown below:}

\begin{itemize}
\item { The case $\mu=a$ or $\mu=b$. We simply have%
\[
\mathcal{D}_{\phi}\left(  t,a;a,b\right)  =\frac{1}{\pi}\int_{a-b}^{0}%
\frac{\sin\left(  tz\right)  }{z}\phi\left(  a-z\right)  dz=\frac{1}{\pi
}\int_{0}^{b-a}\frac{\sin\left(  tz\right)  }{z}\phi\left(  a+z\right)  dz,
\]
which converges to $2^{-1}\phi\left(  a\right)  $ as $t\rightarrow\infty$, or%
\[
\mathcal{D}_{\phi}\left(  t,b;a,b\right)  =\frac{1}{\pi}\int_{0}^{b-a}%
\frac{\sin\left(  tz\right)  }{z}\phi\left(  b-z\right)  dz,
\]
which converges to $2^{-1}\phi\left(  b\right)  $ as $t\rightarrow\infty$.
}

\item { The case $\mu\in\left(  a,b\right)  $. We have%
\begin{align*}
\mathcal{D}_{\phi}\left(  t,\textcolor{black}{\mu};a,b\right)    & =\frac{1}{\pi}\int_{\mu-b}%
^{0}\frac{\sin\left(  tz\right)  }{z}\phi\left(  \mu-z\right)  dz+\frac{1}%
{\pi}\int_{0}^{\mu-a}\frac{\sin\left(  tz\right)  }{z}\phi\left(
\mu-z\right)  dz\\
& =\frac{1}{\pi}\int_{0}^{b-\mu}\frac{\sin\left(  tz\right)  }{z}\phi\left(
\mu+z\right)  dz+\frac{1}{\pi}\int_{0}^{\mu-a}\frac{\sin\left(  tz\right)
}{z}\phi\left(  \mu-z\right)  dz,
\end{align*}
which converges to $2^{-1}\phi\left(  \mu\right)  +2^{-1}\phi\left(
\mu\right)  =\phi\left(  \mu\right)  $ as $t\rightarrow\infty$.
}

\item { The case $\mu<a$ or $\mu>b$. If $\mu<a$, then $\mu-b<\mu-a<0$; if
$\mu>b$, then $\mu-a>\mu-b>0$. Since $\phi$ is of bounded variation on
$\left[  a,b\right]  $, then $\phi$ is Lebesgue integrable on $\left[
a,b\right]  $, which implies that the function $z^{-1}\phi\left(
\mu-z\right)  $ in $z$ is Lebesgue integrable on $\left[  \mu-a,\mu-b\right]
$ when $\mu<a$ or on $\left[  \mu-b,\mu-a\right]  $ when $\mu>b$. Therefore,
by Riemann-Lebesgue Theorem, see., e.g., \cite{Costin:2016}, we must have
that
\[
\mathcal{D}_{\phi}\left(  t,\mu;a,b\right)  =\frac{1}{\pi}\int_{\mu-b}^{\mu
-a}\frac{\sin\left(  tz\right)  }{z}\phi\left(  \mu-z\right)  dz
\]
converges to $0$ as $t\rightarrow\infty$ when $\mu<a$ or $\mu>b$.
}
\end{itemize}

To show the third claim, set $\mathcal{D}%
_{\phi,\infty}\left(  \mu;a,b\right)  =\lim_{t\rightarrow\infty}%
\mathcal{D}_{\phi}\left(  t,\mu;a,b\right)  $ and define%
\[
\mathcal{\hat{D}}_{\phi}\left(  t,\mu;a,b\right)  =   \frac{\phi\left(  \mu\right)  }{\pi}\int_{\mu-b}^{\mu
-a}\frac{\sin\left(  tz\right)  }{z}dz,
\]
{where now $\phi$ is assumed to be defined on $U$.}
Then $\mathcal{D}_{\phi,\infty}\left(  \mu;a,b\right)  =\lim_{t\rightarrow
\infty}\mathcal{\hat{D}}_{\phi}\left(  t,\mu;a,b\right)  $. Define%
\[
\delta_{\mu,a,b}=\min_{\mu\notin\left\{  a,b\right\}  }\left\{  \left\vert
\mu-a\right\vert ,\left\vert \mu-b\right\vert \right\}  .
\]
Then \autoref{lm:Dirichlet} and the calculations in the proof
{presented in \autoref{SpeedBN}} imply, when $t\delta_{\mu,a,b}\geq2$ and $t\left(  b-a\right)
\geq2$,%
\begin{equation}
\left\vert \mathcal{D}_{\phi,\infty}\left(  \mu;a,b\right)  -\mathcal{\hat{D}%
}_{\phi}\left(  t,\mu;a,b\right)  \right\vert \leq\left\{
\begin{array}
[c]{lcl}%
2t^{-1}\left(  b-a\right)  ^{-1}\left\Vert \phi\right\Vert _{\infty} &
\text{if} & \mu=a\text{ or }\mu=b\\
4t^{-1}\delta_{\mu,a,b}^{-1}\left\Vert \phi\right\Vert _{\infty} & \text{if} &
\mu\in\left(  a,b\right)  \\
4t^{-1}\delta_{\mu,a,b}^{-1}\left\Vert \phi\right\Vert _{\infty} & \text{if} &
\mu<a\text{ or }\mu>b
\end{array}
\right.  ,\label{eqDphiInfityOracle}
\end{equation}
i.e.,%
\begin{equation}
\left\vert \mathcal{D}_{\phi,\infty}\left(  \mu;a,b\right)  -\mathcal{\hat{D}%
}_{\phi}\left(  t,\mu;a,b\right)  \right\vert \leq\frac{4\left\Vert
\phi\right\Vert _{\infty}}{t}\left(  \frac{{2^{-1}}}{b-a}+\frac{{1}}{\delta_{\mu,a,b}%
}\right)  .\label{eqdx7}%
\end{equation}

Now we bound $\left\vert \mathcal{D}_{\phi}\left(  t,\mu;a,b\right)
-\mathcal{\hat{D}}_{\phi}\left(  t,\mu;a,b\right)  \right\vert $. Recall for
each $\mu\in {U}  $ and $z\in\left[  \mu-b,\mu-a\right]  $, we
have%
\[
W_{\mu}\left(  z\right)  =z^{-1}\left(  \phi\left(  \mu-z\right)  -\phi\left(
\mu\right)  \right)  \text{ for }\ z\neq0
\]
and%
\[
W_{\mu}^{-}\left(  z\right)  =\left\{
\begin{array}
[c]{lll}%
W_{\mu}\left(  z\right)   & \text{for} & \mu-a\geq z>0\\
\lim_{z\rightarrow0+}W_{\mu}\left(  z\right)   & \text{for} & z=0\text{ if the
limit exists}%
\end{array}
\right.
\]
and%
\[
W_{\mu}^{+}\left(  z\right)  =\left\{
\begin{array}
[c]{lll}%
-W_{\mu}\left(  z\right)   & \text{for} & \mu-b\leq z<0\\
-\lim_{z\rightarrow0-}W_{\mu}\left(  z\right)   & \text{for} & z=0\text{ if
the limit exists}%
\end{array}
\right.  .
\]
From (\ref{eqDphi}), i.e., the definition of $\mathcal{D}_{\phi}\left(
t,\mu;a,b\right)  $, we have%
\begin{align*}
\mathcal{D}_{\phi}\left(  t,\mu;a,b\right)  -\mathcal{\hat{D}}_{\phi}\left(
t,\mu;a,b\right)    & =\frac{1}{\pi}\int_{\mu-b}^{\mu-a}\sin\left(  tz\right)
\frac{\phi\left(  \mu-z\right)  -\phi\left(  \mu\right)  }{z}dz\\
& =d_{1}\left(  t,\mu,a\right)  +d_{2}\left(  t,\mu,b\right)  ,
\end{align*}
where%
\[
d_{1}\left(  t,\mu,a\right)  =\frac{1}{\pi}\int_{0}^{\mu-a}\sin\left(
tz\right)  \frac{\phi\left(  \mu-z\right)  -\phi\left(  \mu\right)  }%
{z}dz=\frac{1}{\pi}\int_{0}^{\mu-a}\sin\left(  tz\right)  W_{\mu}^{-}\left(
z\right)  dz
\]
and%
\[
d_{2}\left(  t,\mu,b\right)  =\frac{1}{\pi}\int_{0}^{b-\mu}\sin\left(
tz\right)  \frac{\phi\left(  \mu+z\right)  -\phi\left(  \mu\right)  }%
{z}dz=\frac{-1}{\pi}\int_{\mu-b}^{0}\sin\left(  tz\right)  W_{\mu}^{+}\left(
z\right)  dz.
\]
Since both $W_{\mu}^{-}$ and $W_{\mu}^{+}$ are well-defined and of bounded
variation for each fixed $\mu\in U$, \autoref{lm:OracleSpeed} implies%
\begin{equation}
\left\vert \mathcal{D}_{\phi}\left(  t,\mu;a,b\right)  -\mathcal{\hat{D}%
}_{\phi}\left(  t,\mu;a,b\right)  \right\vert \leq\frac{4C_{\mu}\left(
\phi\right)  }{\pi\left\vert t\right\vert }\text{ for }t\neq0,\label{eqdx6}%
\end{equation}
where%
\[
C_{\mu}\left(  \phi\right)  =\left\Vert W_{\mu}^{-}\right\Vert _{\infty
}+\left\Vert W_{\mu}^{-}\right\Vert _{\mathrm{TV}}+\left\Vert W_{\mu}%
^{+}\right\Vert _{\infty}+\left\Vert W_{\mu}^{+}\right\Vert _{\mathrm{TV}}.
\]
Combine (\ref{eqdx7}) and
(\ref{eqdx6}), we have
\begin{align}
&  \left\vert \mathcal{D}_{\phi}\left(  t,\mu;a,b\right)  -\mathcal{D}%
_{\phi,\infty}\left(  \mu;a,b\right)  \right\vert \nonumber\\
&  \leq\left\vert \mathcal{D}_{\phi}\left(  t,\mu;a,b\right)  -\mathcal{\hat
{D}}_{\phi}\left(  t,\mu;a,b\right)  \right\vert +\left\vert \mathcal{D}%
_{\phi,\infty}\left(  \mu;a,b\right)  -\mathcal{\hat{D}}_{\phi}\left(
t,\mu;a,b\right)  \right\vert \nonumber\\
&  \leq\frac{4C_{\mu}\left(  \phi\right)  }{\pi t}+\frac{4\left\Vert
\phi\right\Vert _{\infty}}{t}\left(  \frac{{2^{-1}}}{b-a}+\frac{{1}}{\delta_{\mu,a,b}%
}\right)  \text{ when }\min\left\{  t\delta_{\mu,a,b},t\left(  b-a\right)
\right\}  \geq2.\label{eqdx8}%
\end{align}
\qed

\subsection{Speeds of convergence of discriminant functions}
\label{OracleSpeeds}

{
We will provide speeds of convergence of discriminant functions for both the bounded null (in \autoref{SpeedBN}), one-sided null (in \autoref{SpeedOSN})
and the extension (in \autoref{SpeedExtension}), respectively. Due to the differentiable, monotonic, bijection, $\mu\left(  \theta\right)  =\sigma\left(  1-\theta\right)  ^{-1}$, between the mean
$\mu=\mu\left(  \theta\right)  $ and natural parameter $\theta=\theta\left(
\mu\right)  $ for the Gamma family, these speeds will be provided in terms of the parameters $\mu, \boldsymbol{\mu}, a$ and $b$. Recall $\xi\left(\theta\right)=\left(1-\theta\right)^{-1}$ (which implies $\mu\left(  \theta\right)= \sigma \xi\left(\theta\right)$) and the scale parameter $\sigma >0$.
In order to deal with all three settings, i.e., bounded null, one-sided null and the extension, in accordance with the notations used by the theorems and lemmas in the main text (and to be used in a follow-up work on location-shift families),
we will deliberately introduce definitions in the following three sub-subsections such that these definitions are independent of each other throughout these sub-subsections (unless otherwise noted).
}

\subsubsection{Speed of convergence for the bounded null}
\label{SpeedBN}

Define
\begin{equation}
{\tilde{\psi}_{1,0}\left(  t,\mu;\mu^{\prime}\right)  ={\int_{\left[  -1,1\right]  }%
}\omega\left(  s\right)  \cos\left\{  ts \sigma^{-1} \left(  \mu-\mu^{\prime}\right)
\right\}  ds} \label{NewPsiPoint}
\end{equation}
and
\[
\psi_{1}\left(  t,\mu\right)  =\frac{1}{\pi}%
\int_{\left(  \mu-b\right)  t}^{\left(  \mu-a\right)  t}\frac{\sin y}{y}dy.
\]
Define%
\[
\varphi_{1,m}\left(  t,\boldsymbol{\mu}\right)  =m^{-1}\sum_{i=1}^{m}\psi
_{1}\left(  t,\mu_{i}\right)  \text{ \ and \ }\varphi_{1,0,m}\left(
t,\boldsymbol{\mu};\mu^{\prime}\right)  =m^{-1}\sum_{i=1}^{m}{\tilde{\psi}_{1,0}\left(
t,\mu_{i};\mu^{\prime}\right)}.
\]
{Note that $\tilde{\psi}_{1,0}\left(  t,\mu;\mu^{\prime}\right)$ in (\ref{NewPsiPoint}) also depends on $\sigma$. We have:}
\begin{lemma}
\label{SpeedOracleBoundedNull}
  Set $u_{m}=\min_{\tau\in\left\{  a,b\right\}  }\min_{\left\{
j:\mu_{j}\neq\tau\right\}  }\left\vert \mu_{j}-\tau\right\vert $ and%
\[
\varphi_{m}\left(  t,\boldsymbol{\mu}\right)  =1-\varphi_{1,m}\left(  t,\boldsymbol{\mu
}\right)  +2^{-1}\varphi_{1,0,m}\left(  t,\boldsymbol{\mu};a\right)
+2^{-1}\varphi_{1,0,m}\left(  t,\boldsymbol{\mu};b\right)  .
\]
Then, for positive $t$ such that $t\left(  b-a\right)  \geq2$
and $t u_{m}\geq2$,
\begin{equation}
\left\vert \pi_{1,m}^{-1} \varphi_{1,m}\left(  t,\boldsymbol{\mu}\right) -1\right\vert
\leq\frac{4{\sigma}\left(  \left\Vert \omega\right\Vert _{\mathrm{TV}}+\left\Vert
\omega\right\Vert _{\infty}\right)  +12}{\pi_{1,m}u_{m}t}+\frac{4}{\left(
b-a\right)  t\pi_{1,m}}. \label{eq10cx}%
\end{equation}
\end{lemma}

\begin{proof}
First, \autoref{lm:OracleSpeed} implies%
\[
\left\vert {\tilde{\psi}_{1,0}\left(  t,\mu;\mu^{\prime}\right)}  \right\vert
\leq4\left(  \left\Vert \omega\right\Vert _{\mathrm{TV}}+\left\Vert
\omega\right\Vert _{\infty}\right)  \frac{1}{\left\vert t\left(  \mu
-\mu^{\prime}\right)  \right\vert {\sigma^{-1}}}\text{ for }\mu\neq\mu^{\prime},
\]
and%
\begin{equation}
\max_{\tau\in\left\{  a,b\right\}  }\max_{\left\{  j:\mu_{j}\neq\tau\right\}
}\left\vert {\tilde{\psi}_{1,0}\left(  t,\mu_{j};\tau\right)}   \right\vert \leq
\frac{1}{t {\sigma^{-1}}}\max_{\tau\in\left\{  a,b\right\}  }\frac{4\left(  \left\Vert
\omega\right\Vert _{\mathrm{TV}}+\left\Vert \omega\right\Vert _{\infty
}\right)  }{\min_{\left\{  j:\mu_{j}\neq\tau\right\}  }\left\vert \mu_{j}%
-\tau\right\vert }=\frac{4\left(  \left\Vert \omega\right\Vert _{\mathrm{TV}%
}+\left\Vert \omega\right\Vert _{\infty}\right)  }{t u_{m} {\sigma^{-1}}}.\label{eq10ca}%
\end{equation}
Define%
\[
\left\{
\begin{array}
[c]{l}%
\widetilde{d}_{1,m}=1-m^{-1}\sum\nolimits_{\left\{  j:\mu_{j}\in\left(
a,b\right)  \right\}  }\psi_{1}\left(  t ,\mu_{j}\right)  \\
\widetilde{d}_{2,m}=-m^{-1}\sum\nolimits_{\left\{  j:\mu_{j}=a\right\}  }%
\psi_{1}\left(  t ,\mu_{j}\right)  +2^{-1}m^{-1}\sum\nolimits_{\left\{
j:\mu_{j}=a\right\}  }{\tilde{\psi}}_{1,0}\left(  t ,\mu_{j};a\right)  \\
\widetilde{d}_{3,m}=-m^{-1}\sum\nolimits_{\left\{  j:\mu_{j}=b\right\}  }%
\psi_{1}\left(  t ,\mu_{j}\right)  +2^{-1}m^{-1}\sum\nolimits_{\left\{
j:\mu_{j}=b\right\}  }{\tilde{\psi}}_{1,0}\left(  t ,\mu_{j};b\right)  \\
\widetilde{d}_{4,m}=2^{-1}m^{-1}\sum\nolimits_{\left\{  j:\mu_{j}\neq
a\right\}  }{\tilde{\psi}}_{1,0}\left(  t ,\mu_{j};a\right)  +2^{-1}m^{-1}%
\sum\nolimits_{\left\{  j:\mu_{j}\neq b\right\}  }{\tilde{\psi}}_{1,0}\left(  t %
,\mu_{j};b\right)  \\
\widetilde{d}_{5,m}=-m^{-1}\sum\nolimits_{\left\{  j:\mu_{j}<a\right\}  }%
\psi_{1}\left(  t,\mu_{j}\right)  -m^{-1}\sum\nolimits_{\left\{  j:\mu
_{j}>b\right\}  }\psi_{1}\left(  t,\mu_{j}\right)
\end{array}
\right.  .
\]
Then ${\varphi_{m}\left(  t,\boldsymbol{\mu}\right) }=\sum_{i=1}^{5}\widetilde{d}_{1,m}$ and%
\[
\pi_{1,m}^{-1} {\varphi_{m}\left(  t,\boldsymbol{\mu}\right) } =\pi_{1,m}^{-1}\widetilde{d}_{1,m}+\pi_{1,m}%
^{-1}\widetilde{d}_{2,m}+\pi_{1,m}^{-1}\widetilde{d}_{3,m}+\pi_{1,m}%
^{-1}\widetilde{d}_{4,m}+\pi_{1,m}^{-1}\widetilde{d}_{5,m}.
\]
Let us bound each $\widetilde{d}_{1,m}$ for $1\leq i\leq5$. By (\ref{eq10ca}),
we have $\left\vert \widetilde{d}_{4,m}\right\vert \leq4\left(  \left\Vert
\omega\right\Vert _{\mathrm{TV}}+\left\Vert \omega\right\Vert _{\infty
}\right)  \left(  t u_{m} {\sigma^{-1}}\right)  ^{-1}$. Recall%
\[
\psi_{1}\left(  t,\mu\right)  =\frac{1}{\pi}\int_{\left(  \mu-b\right)
t}^{\left(  \mu-a\right)  t}\frac{\sin y}{y}dy%
\overset{t \to\infty}{\longrightarrow}%
\left\{
\begin{array}
[c]{lll}%
1 & \text{if} & a<\mu<b\\
2^{-1} & \text{if} & \mu=a\text{ or }\mu=b\\
0 & \text{if} & \mu<a\text{ or }\mu>b
\end{array}
\right.  .
\]
Then \autoref{lm:Dirichlet} implies%
\begin{align*}
\left\vert \widetilde{d}_{2,m}\right\vert  & \leq m^{-1}\sum
\nolimits_{\left\{  j:\mu_{j}=a\right\}  }\left\vert \psi_{1}\left(  t %
,\mu_{j}\right)  -2^{-1}\right\vert \\
& =m^{-1}\sum\nolimits_{\left\{  j:\mu_{j}=a\right\}  }\left\vert \frac{1}%
{\pi}\int_{\left(  \mu_{j}-b\right)  t }^{0}\frac{\sin y}{y}dy-2^{-1}%
\right\vert \\
& \leq m^{-1}\sum\nolimits_{\left\{  j:\mu_{j}=a\right\}  }\frac{2}{\left(
b-a\right)  t }\leq\frac{2}{\left(  b-a\right)  t }%
\end{align*}
and%
\[
\left\vert \widetilde{d}_{3,m}\right\vert \leq m^{-1}\sum\nolimits_{\left\{
j:\mu_{j}=b\right\}  }\left\vert \frac{1}{\pi}\int_{0}^{\left(  \mu
_{j}-a\right)  t }\frac{\sin y}{y}dy-2^{-1}\right\vert \leq\frac{2}{\left(
b-a\right)  t }.
\]
To deal with $\widetilde{d}_{1,m}$ with $\mu_{j}\in\left(  a,b\right)  $ and
$\widetilde{d}_{5,m}$ with $\mu_{j}<a$ or $\mu_{j}>b$, we have%
\begin{align*}
\psi_{1}\left(  t ,\mu_{j}\right)    & =\frac{1}{\pi}\int_{0}^{\left(
\mu_{j}-a\right)  t }\frac{\sin y}{y}dy+\frac{1}{\pi}\int_{\left(  \mu
_{j}-b\right)  t }^{0}\frac{\sin y}{y}dy\\
& =\frac{1}{\pi}\int_{0}^{\left(  \mu_{j}-a\right)  t }\frac{\sin y}%
{y}dy-\frac{1}{\pi}\int_{0}^{\left(  \mu_{j}-b\right)  t }\frac{\sin y}%
{y}dy.
\end{align*}
So, when $\mu_{j}\in\left(  a,b\right)  $, \autoref{lm:Dirichlet} implies%
\begin{align*}
\left\vert \psi_{1}\left(  t ,\mu_{j}\right)  -1\right\vert  &
\leq\left\vert \frac{1}{\pi}\int_{0}^{\left(  \mu_{j}-a\right)  t }%
\frac{\sin y}{y}dy-2^{-1}\right\vert +\left\vert \frac{1}{\pi}\int_{\left(
\mu_{j}-b\right)  t }^{0}\frac{\sin y}{y}dy-2^{-1}\right\vert \\
& \leq\frac{2}{\left(  \mu_{j}-a\right)  t }+\frac{2}{\left(  b-\mu
_{j}\right)  t }\leq\frac{4}{u_{m}t },
\end{align*}
whereas when $\mu_{j}<a$ or $\mu_{j}>b$, \autoref{lm:Dirichlet} implies%
\begin{align*}
\left\vert \psi_{1}\left(  t ,\mu_{j}\right)  \right\vert  & =\left\vert
\left(  \frac{1}{\pi}\int_{0}^{\left(  \mu_{j}-a\right)  t }\frac{\sin
y}{y}dy-2^{-1}\right)  -\left(  \frac{1}{\pi}\int_{0}^{\left(  \mu
_{j}-b\right)  t }\frac{\sin y}{y}dy-2^{-1}\right)  \right\vert \\
& \leq\frac{2}{\left\vert \mu_{j}-a\right\vert t }+\frac{2}{\left\vert
b-\mu_{j}\right\vert t }\leq\frac{4}{u_{m}t }\text{.}%
\end{align*}
Thus,
\[
\left\vert \widetilde{d}_{5,m}\right\vert \leq m^{-1}\left(  \sum
\nolimits_{\left\{  j:\mu_{j}<a\right\}  }+\sum\nolimits_{\left\{  j:\mu
_{j}>b\right\}  }\right)  \left\vert \psi_{1}\left(  t ,\mu_{j}\right)
\right\vert \leq\frac{8}{u_{m}t },
\]
and by noticing the null set $\Theta_{0}=\left(  a,b\right)  $ and $1=\pi_{0,m}+\pi_{1,m}$,
we see%
\begin{align}
\left\vert \pi_{1,m}^{-1}\widetilde{d}_{1,m}-1\right\vert  &  =\pi_{1,m}%
^{-1}\left\vert \pi_{0,m}-m^{-1}\sum\nolimits_{\left\{  j:\mu_{j}\in\left(
a,b\right)  \right\}  }\psi_{1}\left(  t ,\mu_{j}\right)  \right\vert
\nonumber\\
&  =\pi_{1,m}^{-1}m^{-1}\sum\nolimits_{\left\{  j:\mu_{j}\in\left(
a,b\right)  \right\}  }\left\vert \psi_{1}\left(  t ,\mu_{j}\right)
-1\right\vert \leq\frac{4}{\pi_{1,m}u_{m}t }.\label{eq10cz}%
\end{align}
Combining the bounds for $\tilde{d}_{i,m}$ for all $2\leq i\leq5$ and (\ref{eq10cz}) gives%
\begin{align}
\left\vert \pi_{1,m}^{-1}{\varphi_{m}\left(  t,\boldsymbol{\mu}\right) }-1\right\vert  &  \leq\left\vert
\pi_{1,m}^{-1}\widetilde{d}_{2,m}+\pi_{1,m}^{-1}\widetilde{d}_{3,m}+\pi
_{1,m}^{-1}\widetilde{d}_{4,m}+\pi_{1,m}^{-1}\widetilde{d}_{5,m}+\pi
_{1,m}^{-1}\widetilde{d}_{1,m}-1\right\vert \nonumber\\
&  \leq\frac{4{\sigma}\left(  \left\Vert \omega\right\Vert _{\mathrm{TV}}+\left\Vert
\omega\right\Vert _{\infty}\right) +12}{\pi_{1,m}u_{m}t }+\frac{4}{\left(
b-a\right)  t \pi_{1,m}}\nonumber\\
&  \le C\left(  \frac{1}{\pi_{1,m}u_{m}t }+\frac{1}{t \pi_{1,m}}\right)
.\label{eq10c}%
\end{align}
\qed
\end{proof}

\subsubsection{Speed of convergence for the one-side null}
\label{SpeedOSN}

Recall%
\[
{\tilde{\psi}_{1,0}\left(  t,\mu;\mu^{\prime}\right)  ={\int_{\left[  -1,1\right]  }%
}\omega\left(  s\right)  \cos\left\{  ts \sigma^{-1} \left(  \mu-\mu^{\prime}\right)
\right\}  ds}.
\]
Define
\[
\psi_{1}\left(  t,\mu\right)  =\frac{1}{\pi}\int%
_{0}^{t}\frac{\sin\left(  \left(  \mu-b\right)  y\right)  }{y}dy
\]
and%
\[
\varphi_{1,m}\left(  t,\boldsymbol{\mu}\right)  =m^{-1}\sum_{i=1}^{m}\psi
_{1}\left(  t,\mu_{i}\right)  \text{ \ and \ }\varphi_{1,0,m}\left(
t,\boldsymbol{\mu};b\right)  =m^{-1}\sum_{i=1}^{m}{\tilde{\psi}}_{1,0}\left(  t,\mu
_{i};b\right)  .
\]
Then we have:

\begin{lemma}
\label{SpeedOneSidedNull}

 Define $\tilde{u}_{m}=\min_{\left\{  j:\mu_{j}\neq b\right\}
}\left\vert \mu_{j}-b\right\vert $ and
\[
\varphi_{m}\left(  t,\boldsymbol{\mu}\right)  =2^{-1}+\varphi_{1,m}\left(
t,\boldsymbol{\mu}\right)  +2^{-1}\varphi_{1,0,m}\left(  t,\boldsymbol{\mu
};{b}\right)  .
\]
Then, when $t\tilde{u}_{m}\geq2$,%
\begin{equation}
\left\vert \pi_{1,m}^{-1}\varphi_{m}\left(  t,\boldsymbol{\mu}\right)
-1\right\vert \leq\frac{4+2{\sigma}\left(  \left\Vert \omega\right\Vert _{\infty
}+\left\Vert \omega\right\Vert _{\mathrm{TV}}\right)  }{t\tilde{u}_{m}%
\pi_{1,m}}.\label{eqdx4}%
\end{equation}
\end{lemma}

\begin{proof}
Define%
\[
\left\{
\begin{array}
[c]{l}%
\bar{d}_{1,m}=m^{-1}\sum\nolimits_{\left\{  i:\mu_{i}>{b}\right\}  }\left(
2^{-1}+\psi_{1}\left(  t,\mu_{i}\right)  \right)  \\
\bar{d}_{2,m}=m^{-1}\sum\nolimits_{\left\{  i:\mu_{i}={b}\right\}  }\left(
2^{-1}+\psi_{1}\left(  t,\mu_{i}\right)  +2^{-1}{\tilde{\psi}}_{1,0}\left(  t,\mu
_{i};{b}\right)  \right)  \\
\bar{d}_{3,m}=m^{-1}\sum\nolimits_{\left\{  i:\mu_{i}<{b}\right\}  }\left(
2^{-1}+\psi_{1}\left(  t,\mu_{i}\right)  \right)  \\
\bar{d}_{4,m}=2^{-1}m^{-1}\sum\nolimits_{\left\{  i:\mu_{i}\neq {b}\right\}
}{\tilde{\psi}}_{1,0}\left(  t,\mu_{i};{b}\right)
\end{array}
\right.  .
\]
Then%
\[
\varphi_{m}\left(  t,\boldsymbol{\mu}\right)  =2^{-1}+\varphi_{1,m}\left(
t,\boldsymbol{\mu}\right)  +2^{-1}\varphi_{1,0,m}\left(  t,\boldsymbol{\mu
};{b}\right)  =\bar{d}_{1,m}+\bar{d}_{2,m}+\bar{d}_{3,m}+\bar{d}_{4,m}.
\]
Recall the alternative set $\Theta_{1}=[{b},\infty)$ and $\pi_{1,m}=m^{-1}%
\sum_{i=1}^{m}1_{\Theta_{1}}\left(  \mu_{i}\right)  $. Since ${\tilde{\psi}_{1,0}}\left(
t,{b};{b}\right)  =1$ for $\forall t$ and%
\[
\psi_{1}\left(  t,\mu\right)  =\frac{1}{\pi}\int_{0}^{t}\frac{\sin\left(  {\left(\mu-b\right)}
y\right)  }{y}dy%
\overset{t \to\infty}{\longrightarrow}%
\left\{
\begin{array}
[c]{lll}%
2^{-1} & \text{if} & \mu>{b}\\
0 & \text{if} & \mu={b}\\
-2^{-1} & \text{if} & \mu<{b}
\end{array}
\right.  ,
\]
then $\bar{d}_{2,m}=m^{-1}\sum\nolimits_{\left\{  i:\mu_{i}={b}\right\}  }1$ and%
\begin{align*}
\pi_{1,m}^{-1}\left(  \bar{d}_{1,m}+\bar{d}_{2,m}\right)  -1  & =\pi
_{1,m}^{-1}\left(  m^{-1}\sum\nolimits_{\left\{  i:\mu_{i}\geq {b}\right\}
}1+m^{-1}\sum\nolimits_{\left\{  i:\mu_{i}>{b}\right\}  }\left(  \psi_{1}\left(
t,\mu_{i}\right)  -2^{-1}\right)  \right)  -1\\
& =\pi_{1,m}^{-1}m^{-1}\sum\nolimits_{\left\{  i:\mu_{i}>{b}\right\}  }\left(
\psi_{1}\left(  t,\mu_{i}\right)  -2^{-1}\right)  .
\end{align*}
Further, \autoref{lm:OracleSpeed} and \autoref{lm:Dirichlet} imply%
\[
\left\{
\begin{array}
[c]{l}%
\left\vert \pi_{1,m}^{-1}\left(  \bar{d}_{1,m}+\bar{d}_{2,m}\right)
-1\right\vert \leq2\pi_{1,m}^{-1}\left(  t\tilde{u}_{m}\right)  ^{-1}\text{
when }t\tilde{u}_{m}\geq2\\
\left\vert \bar{d}_{3,m}\right\vert \leq2\left(  \left\vert t\right\vert
\tilde{u}_{m}\right)  ^{-1}\text{ \ when }t\tilde{u}_{m}\geq2\\
\left\vert \bar{d}_{4,m}\right\vert \leq2\left(  \left\vert t\right\vert
\tilde{u}_{m} {\sigma^{-1}}\right)  ^{-1}\left(  \left\Vert \omega\right\Vert _{\infty
}+\left\Vert \omega\right\Vert _{\mathrm{TV}}\right)  \text{ when }t\neq0
\end{array}
\right.  ,
\]
where $\tilde{u}_{m}={\min_{\left\{  j:\mu_{j}\neq b\right\}
}\left\vert \mu_{j}-b\right\vert }$. Thus, when $t\tilde{u}_{m}\geq2$,
\begin{equation}
\left\vert \pi_{1,m}^{-1}\varphi_{m}\left(  t,\boldsymbol{\mu}\right)
-1\right\vert \leq\frac{4+2{\sigma}\left(  \left\Vert \omega\right\Vert _{\infty
}+\left\Vert \omega\right\Vert _{\mathrm{TV}}\right)  }{t\tilde{u}_{m}%
\pi_{1,m}}\leq\frac{C}{t\tilde{u}_{m}\pi_{1,m}}.\nonumber 
\end{equation}
\qed
\end{proof}

\subsubsection{Speed of convergence for the extension}
\label{SpeedExtension}

Define%
\[
\psi_{1}\left(  t,\mu\right)  =\mathcal{D}_{\phi}\left(  t,\mu;a,b\right)
=\frac{1}{\pi}\int_{a}^{b}\frac{\sin\left\{  \left(  \mu-y\right)  t\right\}
}{\mu-y}\phi\left(  y\right)  dy
\]
and recall%
\[
{\tilde{\psi}_{1,0}\left(  t,\mu;\mu^{\prime}\right)  ={\int_{\left[  -1,1\right]  }%
}\omega\left(  s\right)  \cos\left\{  ts \sigma^{-1} \left(  \mu-\mu^{\prime}\right)
\right\}  ds}
\]
and $u_{m}=\min_{\tau\in\left\{  a,b\right\}  }\min_{\left\{  j:\mu_{j}%
\neq\tau\right\}  }\left\vert \mu_{j}-\tau\right\vert $ and $\left\Vert
\phi\right\Vert _{1,\infty}=\sup_{\mu\in \textcolor{black}{U}  }C_{\mu}\left(
\phi\right)  $.

Then we have:
\begin{lemma}
\label{SpeedOracleExt}
Define%
\[
\varphi_{m}\left(  t,\boldsymbol{\mu}\right)  =m^{-1}\sum_{i=1}^{m}\left[
\psi_{1}\left(  t,\mu_{i}\right)  -2^{-1}\left\{  \phi\left(  a\right)
{\tilde{\psi}}_{1,0}\left(  t,\mu;a\right)  +\phi\left(  b\right)  {\tilde{\psi}}_{1,0}\left(
t,\mu;b\right)  \right\}  \right]  .
\]
Then, for $t$ such that $tu_{m}\geq2$ and $t\left(  b-a\right)  \geq2$,%
\begin{equation}
\left\vert \check{\pi}_{0,m}^{-1}\varphi_{m}\left(  t,\boldsymbol{\mu}\right)
-1\right\vert \leq\frac{\textcolor{black}{20}\left\Vert \phi\right\Vert _{1,\infty}}{\pi
t\check{\pi}_{0,m}}+\frac{\textcolor{black}{20}\left\Vert \phi\right\Vert _{\infty}}{tu_{m}%
\check{\pi}_{0,m}}+\frac{\textcolor{black}{10}\left\Vert \phi\right\Vert _{\infty}}{t\left(
b-a\right)  \check{\pi}_{0,m}}+\frac{4\left(  \left\Vert \omega\right\Vert
_{\mathrm{TV}}+\left\Vert \omega\right\Vert _{\infty}\right)  \left\Vert
\phi\right\Vert _{\infty}}{tu_{m}\check{\pi}_{0,m} {\sigma^{-1}}}.\label{eqdx13xx}%
\end{equation}
\end{lemma}

\begin{proof}

First, \autoref{lm:OracleSpeed} implies%
\begin{equation}
\max_{\tau\in\left\{  a,b\right\}  }\left\vert {\tilde{\psi}}_{1,0}\left(  t,\mu
_{j};\tau\right)  \right\vert \leq\frac{4\left(  \left\Vert \omega\right\Vert
_{\mathrm{TV}}+\left\Vert \omega\right\Vert _{\infty}\right)  }{tu_{m} {\sigma^{-1}}}\text{
for }t\neq0.\label{eqdx12}%
\end{equation}
Define%
\[
\left\{
\begin{array}
[c]{l}%
d_{\phi,1}\left(  t,\boldsymbol{\mu}\right)  =m^{-1}\sum\nolimits_{\left\{
i:\mu_{i}\in\left(  a,b\right)  \right\}  }\psi_{1}\left(  t,\mu_{i}\right)
\\
d_{\phi,2}\left(  t,\boldsymbol{\mu}\right)  =m^{-1}\sum\nolimits_{\left\{
i:\mu_{i}=a\right\}  }\left(  \psi_{1}\left(  t,\mu_{i}\right)  -2^{-1}%
\phi\left(  a\right)  {\tilde{\psi}}_{1,0}\left(  t,\mu_{i};a\right)  \right)  \\
d_{\phi,3}\left(  t,\boldsymbol{\mu}\right)  =m^{-1}\sum\nolimits_{\left\{
i:\mu_{i}=b\right\}  }\left(  \psi_{1}\left(  t,\mu_{i}\right)  -2^{-1}%
\phi\left(  b\right)  {\tilde{\psi}}_{1,0}\left(  t,\mu_{i};b\right)  \right)  \\
d_{\phi,4}\left(  t,\boldsymbol{\mu}\right)  =-m^{-1}\left(  \sum_{\left\{
i:\mu_{i}\neq a\right\}  }+\sum_{\left\{  i:\mu_{i}\neq b\right\}  }\right)
\left\{  2^{-1}\left[  \phi\left(  a\right)  {\tilde{\psi}}_{1,0}\left(  t,\mu
_{i};a\right)  +\phi\left(  b\right)  {\tilde{\psi}}_{1,0}\left(  t,\mu_{i};b\right)
\right]  \right\}  \\
d_{\phi,5}\left(  t,\boldsymbol{\mu}\right)  =m^{-1}\sum_{\left\{  i:\mu
_{i}<a\right\}  }\psi_{1}\left(  t,\mu_{i}\right)  +m^{-1}\sum_{\left\{
i:\mu_{i}>b\right\}  }\psi_{1}\left(  t,\mu_{i}\right)
\end{array}
\right..
\]
Then $\varphi_{m}\left(  t,\boldsymbol{\mu}\right)  =\sum_{j=1}^{5}d_{\phi
,j}\left(  t,\boldsymbol{\mu}\right)  $. \textcolor{black}{Now applying \autoref{ThmDirichletInterval}
once to deal with each of $d_{\phi,1}$, $d_{\phi,2}$
and $d_{\phi,3}$ but twice for $d_{\phi,5}$, and applying inequality (\ref{eqdx12}) twice to deal
with $d_{\phi,4}$ which has a $2^{-1}$, we get}
\begin{align}
\left\vert \check{\pi}_{0,m}^{-1}\varphi_{m}\left(  t,\boldsymbol{\mu}\right)
-1\right\vert  &  \leq\left\vert \check{\pi}_{0,m}^{-1}d_{\phi,1}\left(
t,\boldsymbol{\mu}\right)  -1\right\vert +\sum_{j=2}^{5}\check{\pi}_{0,m}%
^{-1}\left\vert d_{\phi,j}\left(  t,\boldsymbol{\mu}\right)  \right\vert
\nonumber\\
&  \leq\frac{\textcolor{black}{20}\left\Vert \phi\right\Vert _{1,\infty}}{\pi t\check{\pi}_{0,m}%
}+\frac{\textcolor{black}{20}\left\Vert \phi\right\Vert _{\infty}}{tu_{m}\check{\pi}_{0,m}}%
+\frac{\textcolor{black}{10}\left\Vert \phi\right\Vert _{\infty}}{t\left(  b-a\right)  \check
{\pi}_{0,m}}+\frac{4\left(  \left\Vert \omega\right\Vert _{\mathrm{TV}%
}+\left\Vert \omega\right\Vert _{\infty}\right) \left\Vert \phi\right\Vert _{\infty} }{tu_{m}\check{\pi}_{0,m}{\sigma^{-1}}%
}\nonumber\\
&  \leq\frac{C}{t\check{\pi}_{0,m}}\left(  1+\left\Vert \phi\right\Vert
_{1,\infty}+\frac{1}{u_{m}}\right)  \label{eqdx13}%
\end{align}
for $t$ such that $tu_{m}\geq2$ and $t\left(  b-a\right)  \geq2$, where
$\left\Vert \phi\right\Vert _{1,\infty}=\sup_{\mu\in\left[  a,b\right]
}C_{\mu}\left(  \phi\right)  $. Note that
(\ref{eqdx13}) is the analogue of (\ref{eq10c}).\qed
\end{proof}

\subsection{Unavailability of better concentration inequalities}
\label{NotBetterConcentration}

In order to construct the estimators, we have employed the following empirical
processes \textcolor{black}{(see the proofs of Theorem 3 and Theorem 7 in the supplementary material)}%
\[
\left\{
\begin{array}
[c]{l}%
S_{1,m}\left(  t,y\right)  =m^{-1}\sum_{i=1}^{m}\left[  w_{1}\left(
t,z_{i},y\right)  -\mathbb{E}\left\{  w_{1}\left(  t,z_{i},y\right)  \right\}
\right]  \\
\Delta_{3,m,1}\left(  t,y,\mathbf{z}\right)  =m^{-1}\sum_{i=1}^{m}\left[
w_{3,1}\left(  t,z_{i},y\right)  -\mathbb{E}\left\{  w_{3,1}\left(
t,z_{i},y\right)  \right\}  \right]  \\
\Delta_{3,m,2}\left(  t,y,\mathbf{z}\right)  =m^{-1}\sum_{i=1}^{m}\left[
w_{3,2}\left(  t,z_{i},y\right)  -\mathbb{E}\left\{  w_{3,2}\left(
t,z_{i},y\right)  \right\}  \right]
\end{array}
\right.  ,
\]
where, for $t\geq0,x>0$ and $y\in\left[  0,1\right]  $ or $\left[  a,b\right]
$,%
\[
\left\{
\begin{array}
[c]{l}%
w_{1}\left(  t,x,y\right)  =\Gamma\left(  \sigma\right)  \sum_{n=0}^{\infty
}\frac{\left(  tx\sigma\right)  ^{n}\cos\left(  2^{-1}n\pi-ty\right)
}{n!\Gamma\left(  n+\sigma\right)  }\\
w_{3,1}\left(  t,x,y\right)  =\Gamma\left(  \sigma\right)  \sum_{n=0}^{\infty
}\cos\left(  2^{-1}\pi n-tyb\right)  \frac{\left(  ty\right)  ^{n}}{n!}%
\frac{\left(  \sigma x\right)  ^{n+1}}{\Gamma\left(  \sigma+n+1\right)  }\\
w_{3,2}\left(  t,x,y\right)  =\Gamma\left(  \sigma\right)  \sum_{n=0}^{\infty
}\cos\left(  2^{-1}\pi n-tyb\right)  \frac{\left(  ty\right)  ^{n}\left(
\sigma x\right)  ^{n}}{n!\Gamma\left(  \sigma+n\right)  }%
\end{array}
\right.  .
\]

Even thought for each fixed $y$, each of the functions $w_{1}\left(
t,x,y\right)  $, $w_{3,1}\left(  t,x,y\right)  $ and $w_{3,2}\left(
t,x,y\right)  $ is analytic in $tx$, none of them has an explicit, analytic
expression, none of them is a polynomial function or Lipschitz function of
$tx$, and we were unable to obtain relative tight upper bounds for any of them
based on results in existing literature. So, results on concentrations of
empirical processes induced by polynomial functions or Lipschitz functions
cannot be applied to any of the above three empirical processes. In fact, it
may well be an open problem to derive good concentration inequalities for such
processes. Further, as shown in the supplementary material, bounding the variance of any of the
three empirical processes already involves complicated asymptotics of Bessel
functions and another special function, none of which has an explicit,
analytic expression. So, it is very complicated, if not impossible, to
relatively tightly bound higher order moments of any of these empirical
processes, in order to derive better concentration inequalities for them than
already provided by this work. This explains why we were not able to obtain
better concentration inequalities for the error of the estimators.

\section*{Acknowledgements}

This research was funded by the New Faculty Seed Grant provided by Washington
State University. I am grateful to Prof. G\'{e}rard Letac for his constant
guidance and feedback on an earlier version of the manuscript and to Prof.
Jiashun Jin for his warm encouragements. In particular, the second part of the
proof of \autoref{MainRes} was inspired by Prof. G\'{e}rard Letac's suggestion
on examining the ratio $\tilde{c}_{n+1}\left(  \theta\right)  /\tilde{c}%
_{n}\left(  \theta\right)  $ via the relationship
$\frac{d}{d\theta}\int x^{n}G_{\theta}\left(  dx\right)  =\int x^{n+1}%
G_{\theta}\left(  dx\right) $ for the Gamma family introduced in \autoref{SecBackground}.
I am also very grateful to the two anonymous reviewers for their very detailed and critical comments that helped greatly improve
the clarity and presentation of the manuscript.

\bibliographystyle{dcu}


\newpage
\renewcommand*{\thefootnote}{\fnsymbol{footnote}}
\numberwithin{figure}{section}

\begin{center}
{\Large{Supplementary Material to ``Uniformly consistent proportion estimation for composite hypotheses via integral equations: `the case of Gamma random variables'"\\
{\small Xiongzhi Chen (xiongzhi.chen@wsu.edu)}\\
{\small Department of Mathematics and Statistics, Washington State University}}}
\bigskip
\end{center}

\autoref{AppProofsNEFA} contains proofs related to Construction I, \autoref{AppProofsNEFB} proofs related to Construction II,
\autoref{AppProofsExt} proofs related to the extension of Construction I, \autoref{SecNumericalStudies} a simulation study on the
proposed estimators, and \autoref{SecClosedNulls} how to adapt the constructions and estimators to the null sets $\left(-\infty,b\right]$ and $\left[a,b\right]$.


\section{Proofs Related to Construction I}

\label{AppProofsNEFA}

\subsection{Proof of \autoref{ThmConstructionMoments}}

Recall $\tilde{c}_{n}\left(  \theta\right)  =\int x^{n}dG_{\theta}\left(
x\right)  $ and define%
\begin{equation}
K_{1}^{\dag}\left(  t,x\right)  =\frac{1}{2\pi\zeta_{0}}\int_{a}^{b}%
tdy\int_{-1}^{1}\exp\left(  -\iota tsy\right)  \sum_{n=0}^{\infty}%
\frac{\left(  \iota tsx\zeta_{0}\tilde{a}_{1}\right)  ^{n}}{\tilde{a}_{n}%
n!}ds\text{.} \label{eq1e}%
\end{equation}
By assumption, $\tilde{c}_{n}\left(  \theta\right)  =\xi^{n}\left(
\theta\right)  \zeta\left(  \theta\right)  \tilde{a}_{n}=\zeta_{0}\xi
^{n}\left(  \theta\right)  \tilde{a}_{n}$, where $\zeta_{0}\equiv\zeta\equiv
1$. So, $\mu\left(  \theta\right)  =\xi\left(  \theta\right)  \zeta\left(
\theta\right)  \tilde{a}_{1}=\xi\left(  \theta\right)  \zeta_{0}\tilde{a}_{1}$
and
\begin{align*}
\psi_{1}\left(  t,\theta\right)   &  = \int K_{1}^{\dag
}\left(  t,x\right)  dG_{\theta}\left(  x\right) \\
&  =\frac{1}{2\pi\zeta_{0}}\int_{a}^{b}
tdy\int_{-1}^{1} \exp\left(  -\iota tsy\right)\sum_{n=0}^{\infty}\frac{\left(  \iota ts\zeta_{0}\tilde
{a}_{1}\right)  ^{n}}{\tilde{a}_{n}n!}\tilde{c}_{n}\left(  \theta\right)  ds\\
&  =\frac{\zeta\left(  \theta\right)  }{2\pi\zeta_{0}}\int_{a}^{b}
 tdy\int_{-1}^{1} \exp\left(  -\iota tsy\right) \sum_{n=0}^{\infty}\frac{\left(  \iota
ts\zeta_{0}\tilde{a}_{1}\right)  ^{n}}{n!}\xi^{n}\left(  \theta\right)  ds\\
&  =\frac{1}{2\pi}\int_{a}^{b}tdy\int_{-1}^{1}\exp\left[  \iota ts\left\{
\mu\left(  \theta\right)  -y\right\}  \right]  ds.
\end{align*}
Since $\psi_{1}$ is real, $\psi_{1}=\mathbb{E}\left\{  \Re\left(  K_{1}^{\dag
}\right)  \right\}  $. However,%
\[
K_{1}\left(  t,x\right)  =\Re\left\{  K_{1}^{\dag}\left(  t,x\right)
\right\}  =\frac{1}{2\pi\zeta_{0}}\int_{a}^{b}tdy\int_{-1}^{1}\sum
_{n=0}^{\infty}\frac{\left(  tsx\zeta_{0}\tilde{a}_{1}\right)  ^{n}\cos\left(
2^{-1}n\pi-tsy\right)  }{\tilde{a}_{n}n!}ds.
\]
Since $\mu\left(  \theta\right)  $ is smooth and strictly increasing in
$\theta\in\Theta$, $a\leq\mu\leq b$ if and only if $\theta_{a}\leq\theta
\leq\theta_{b}$. By \autoref{ThmPoinNull}, the pair $\left(  K,\psi\right)  $
in (\ref{IV-b}) is as desired.\qed

\subsection{Proof of \autoref{ConcentrationIII}}

In order the present the proof, we quote Lemma 4 of \cite{Chen:2018a} as
follows: for a fixed $\sigma>0$, let%
\begin{equation}
\tilde{w}\left(  z,x\right)  =\sum_{n=0}^{\infty}\frac{\left(  zx\right)
^{n}}{n!\Gamma\left(  \sigma+n\right)  }\text{ \ for }z,x>0. \label{eq15a}%
\end{equation}
If $Z$ has CDF $G_{\theta}$ from the Gamma family with scale parameter
$\sigma$, then%
\begin{equation}
\mathbb{E}\left[  \tilde{w}^{2}\left(  z,Z\right)  \right]  \leq C\left(
\frac{z}{1-\theta}\right)  ^{3/4-\sigma}\exp\left(  \frac{4z}{1-\theta
}\right)  \label{eq15}%
\end{equation}
for positive and sufficiently large $z$.

Now we present the arguments. Firstly, we
will obtain an upper bound for $\mathbb{V}\left\{  \hat{\varphi}_{m}\left(
t,\mathbf{z}\right)  \right\}  $. For Gamma family, $\zeta\left(
\theta\right)  \equiv\zeta_{0}=1$, $\tilde{a}_{1}=\sigma$ and $\mu\left(
\theta\right)  =\sigma\xi\left(  \theta\right)  $. Define%
\[
w_{1}\left(  t,x,y\right)  =\Gamma\left(  \sigma\right)  \sum_{n=0}^{\infty
}\frac{\left(  tx\sigma\right)  ^{n}\cos\left(  2^{-1}n\pi-ty\right)
}{n!\Gamma\left(  n+\sigma\right)  }\text{ \ for }t\geq0\text{ and }x>0,
\]
and set $S_{1,m}\left(  t, y\right)  =m^{-1}\sum_{i=1}^{m}\left[  w_{1}\left(
t,z_{i},y\right)  -\mathbb{E}\left\{  w_{1}\left(  t,z_{i},y\right)  \right\}
\right]  $.
Recall $\tilde{a}_{n}=\frac{\Gamma\left(  n+\sigma\right)
}{\Gamma\left(  \sigma\right)  }$.
Then%
\[
K_{1}\left(  t,x\right)  =\frac{1}{2\pi}\int_{a}^{b}tdy\int_{-1}^{1}%
w_{1}\left(  ts,x,y\right)  ds.
\]
Define $\tilde{V}_{1,m}=\mathbb{V}\left\{  \hat{\varphi}_{1,m}\left(
t,\mathbf{z}\right)  \right\}  $, where $\hat{\varphi}_{1,m}\left(
t,\mathbf{z}\right)  =m^{-1}\sum_{i=1}^{m}K_{1}\left(  t,z_{i}\right)  $ and
$\varphi_{1,m}\left(  t,\boldsymbol{\theta}\right)  =\mathbb{E}\left\{
\hat{\varphi}_{1,m}\left(  t,\mathbf{z}\right)  \right\}  $. Then,
applying H\"{o}lder's inequality to $\int_{-1}^{1}\left\vert S_{1,m}\left(
ts,y\right)  \right\vert ds$ and then to $\int_{a}^{b}dy\left[  \left(  \int%
_{-1}^{1}\left\vert S_{1,m}\left(  ts,y\right)  \right\vert ^{2}ds\right)
^{1/2}\right]  $,%
\begin{align*}
\tilde{V}_{1,m}  & =\mathbb{E}\left[  \left\{  \frac{1}{2\pi}\int_{a}%
^{b}tdy\int_{-1}^{1}S_{1,m}\left(  ts,y\right)  ds\right\}  ^{2}\right]  \\
& \leq\frac{t^{2}}{4\pi^{2}}\mathbb{E}\left[  \left\{  \int_{a}^{b}dy\left[
\sqrt{2}\left(  \int_{-1}^{1}\left\vert S_{1,m}\left(  ts,y\right)  \right\vert
^{2}ds\right)  ^{1/2}\right]  \right\}  ^{2}\right]  \\
& =\frac{t^{2}}{2\pi^{2}}\mathbb{E}\left[  \left\{  \int_{a}^{b}\left(
\int_{-1}^{1}\left\vert S_{1,m}\left(  ts,y\right)  \right\vert ^{2}ds\right)
^{1/2}dy\right\}  ^{2}\right]  \\
& \leq\frac{t^{2}}{2\pi^{2}}\mathbb{E}\left[  \left(  b-a\right)  \int_{a}%
^{b}\int_{-1}^{1}\left\vert S_{1,m}\left(  ts,y\right)  \right\vert
^{2}dy\right]  \\
& =\frac{\left(  b-a\right)  t^{2}}{2\pi^{2}}\mathbb{E}\left\{  \int_{a}%
^{b}dy\int_{-1}^{1}\left\vert S_{1,m}\left(  ts,y\right)  \right\vert
^{2}ds\right\}  ,
\end{align*}
i.e.,
\begin{equation}
\tilde{V}_{1,m}=\mathbb{E}\left[  \left\{  \frac{1}{2\pi}\int_{a}^{b}%
tdy\int_{-1}^{1}S_{1,m}\left(  ts,y\right)  ds\right\}  ^{2}\right]  \leq
\frac{\left(  b-a\right)  t^{2}}{2\pi^{2}}\mathbb{E}\left\{  \int_{a}%
^{b}dy\int_{-1}^{1}\left\vert S_{1,m}\left(  ts,y\right)  \right\vert
^{2}ds\right\}  .\label{eqdx2}%
\end{equation}
Since $\left\vert w_{1}\left(  t,x,y\right)  \right\vert \leq\Gamma\left(
\sigma\right)  \tilde{w}\left(  t\sigma,x\right)  $ uniformly in $\left(
t,x,y\right)  $, the inequality (\ref{eq15}) implies, for $t>0$ sufficiently large,%
\begin{align*}
\tilde{V}_{1,m}  &  \leq Ct^{2}\mathbb{E}\left\{  \int_{a}^{b}dy\int_{-1}%
^{1}\left\vert S_{1,m}\left(  ts,y\right)  \right\vert ^{2}ds\right\}  \leq
\frac{Ct^{2}}{m^{2}}\sum_{i=1}^{m}\mathbb{E}\left[  \tilde{w}^{2}\left(
t\sigma,z_{i}\right)  \right] \\
&  \leq\frac{Ct^{2}}{m^{2}}\sum_{i=1}^{m}\left(  \frac{t}{1-\theta_{i}%
}\right)  ^{3/4-\sigma}\exp\left(  \frac{4t\sigma}{1-\theta_{i}}\right)
\leq\frac{Ct^{2}}{m}V_{1,m},
\end{align*}
where we recall $u_{3,m}=\min_{1\leq i\leq m}\left\{  1-\theta_{i}\right\}  $
and have set%
\[
V_{1,m}=\frac{1}{m}\exp\left(  \frac{4t\sigma}{u_{3,m}}\right)  \sum_{i=1}%
^{m}\left(  \frac{t}{1-\theta_{i}}\right)  ^{3/4-\sigma}.
\]

Recall for $\tau\in\left\{  a,b\right\}  $%
\[
K_{3,0}\left(  t,x;\theta_{\tau}\right)  =\frac{\Gamma\left(  \sigma\right)
}{\zeta_{0}}\int_{\left[  -1,1\right]  }\sum_{n=0}^{\infty}\dfrac{\left(
-tsx\right)  ^{n}\cos\left\{  2^{-1}\pi n+ts\xi\left(  \theta_{\tau}\right)
\right\}  }{n!\Gamma\left(  n+\sigma\right)  }\omega\left(  s\right)  ds.
\]
Define $\hat{\varphi}_{3,0,m}\left(  t,\mathbf{z};\tau\right)  =m^{-1}%
\sum_{i=1}^{m}K_{3,0}\left(  t,z_{i};\theta_{\tau}\right)  $ and
$\varphi_{3,0,m}\left(  t,\boldsymbol{\theta};\tau\right)  =\mathbb{E}\left\{
\hat{\varphi}_{3,0,m}\left(  t,\mathbf{z};\tau\right)  \right\}  $. Then
Theorem 8 of \cite{Chen:2018a} implies, for $t>0$ sufficiently large,%
\[
\mathbb{V}\left\{  \hat{\varphi}_{3,0,m}\left(  t,\mathbf{z};\tau\right)
\right\}  \leq Cm^{-1}V_{0,m}\text{ \ with \ }V_{0,m}=\frac{1}{m}\exp\left(
\frac{4t}{u_{3,m}}\right)  \sum_{i=1}^{m}\frac{t^{3/4-\sigma}}{\left(
1-\theta_{i}\right)  ^{3/4-\sigma}}.
\]
So, for $t>0$ sufficiently large,
\begin{equation}
\mathbb{V}\left\{  \hat{\varphi}_{m}\left(  t,\mathbf{z}\right)  \right\}
\leq Cm^{-1}V_{0,m}+Ct^{2}m^{-1}V_{1,m}\leq Cm^{-1}\left(  1+t^{2}\right)
\tilde{V}_{1,m}^{\ast}, \label{eq20a}%
\end{equation}
where%
\begin{equation}
\tilde{V}_{1,m}^{\ast}=\frac{1}{m}\exp\left(  \frac{4t\max\left\{
\sigma,1\right\}  }{u_{3,m}}\right)  \sum_{i=1}^{m}\left(  \frac{t}%
{1-\theta_{i}}\right)  ^{3/4-\sigma}. \label{eq20a1}%
\end{equation}

Secondly, we provide a uniform consistency class. If $\sigma\geq 3/4$, then
(\ref{eq20a1}) induces%
\begin{equation}
\tilde{V}_{1,m}^{\ast}\leq t^{3/4-\sigma}\exp\left(
\frac{4 \max\left\{  \sigma,1\right\} t}{u_{3,m}}\right)  \left\Vert 1-\boldsymbol{\theta}\right\Vert
_{\infty}^{\sigma-3/4} \label{eq20}%
\end{equation}
for $t>0$ sufficiently large, where $\left\Vert 1-\boldsymbol{\theta}\right\Vert _{\infty}=\max_{1\leq i\leq
m}\left(  1-\theta_{i}\right)  $. Let $\varepsilon>0$ be a constant and set
$t_m=\left(  4 \max\left\{  \sigma,1\right\}\right)  ^{-1}u_{3,m}\gamma\ln m$ for any fixed $\gamma
\in \left(  0,1\right)  $. Then, (\ref{eq20a}) and (\ref{eq20}) imply, for all $m$ such that $t_m$ is sufficiently large,%
\begin{equation}
\Pr\left\{  \frac{\left\vert \hat{\varphi}_{m}\left(  t_m,\mathbf{z}\right)
-\varphi_{m}\left(  t_m,\boldsymbol{\theta}\right)  \right\vert }{\pi_{1,m}}%
\geq\varepsilon\right\}  \leq\frac{C\left\Vert 1-\boldsymbol{\theta
}\right\Vert _{\infty}^{\sigma-3/4}}{\varepsilon^{2}m^{1-\gamma}\pi_{1,m}^{2}%
}\left(  u_{3,m}\ln m\right)  ^{11/4-\sigma}. \label{eq18c}%
\end{equation}
Note that $\gamma=1$ can be set when $\sigma>11/4$ since $\lim_{m\rightarrow
\infty}\left(  \ln m\right)  ^{11/4-\sigma}=0$ for all such $\sigma$.
In contrast, if $\sigma\leq3/4$, then (\ref{eq20a1}) implies%
\[
\tilde{V}_{1,m}^{\ast}\leq \left(  \frac{t}{u_{3,m}%
}\right)  ^{3/4-\sigma}\exp\left(  \frac{4t}{u_{3,m}}\right)
\]
for all $t>0$ sufficiently large. Set $t_m=4^{-1}u_{3,m}\gamma\ln m$ for any fixed $\gamma\in\left(  0,1\right)
$. Then, for all $m$ such that $t_m$ is sufficiently large%
\begin{equation}
\Pr\left\{  \frac{\left\vert \hat{\varphi}_{m}\left(  t_m,\mathbf{z}\right)
-\varphi_{m}\left(  t_m,\boldsymbol{\theta}\right)  \right\vert }{\pi_{1,m}}%
\geq\varepsilon\right\}  \leq\frac{C\left(  \ln m\right)  ^{11/4-\sigma
}u_{3,m}^{2}}{\varepsilon^{2}m^{1-\gamma}\pi_{1,m}^{2}}. \label{eq18b}%
\end{equation}

To determine a uniform consistency class, we only need to incorporate the
speed of convergence of $\varphi_m\left(t,\boldsymbol{\mu}\right)$ to $\pi_{1,m}$. Recall for $\tau\in\left\{  a,b\right\}
$
\[
\psi_{3,0}\left(  t,\theta;\theta_{\tau}\right)  =\int_{\left[  -1,1\right]
}\cos\left[  ts\left\{  \xi\left(  \theta_{\tau}\right)  -\xi\left(
\theta\right)  \right\}  \right]  \omega\left(  s\right)  ds,
\]
{
which is exactly
\[
\tilde{\psi}_{1,0}\left(  t,\mu;\mu^{\prime}\right)  ={\int_{\left[  -1,1\right]  }%
}\omega\left(  s\right)  \cos\left\{  ts \sigma^{-1} \left(  \mu-\mu^{\prime}\right)
\right\}  ds
\]
that is defined by (\ref{NewPsiPoint}) (in \autoref{SpeedOracleBoundedNull}) but evaluated at $\mu^{\prime}=\tau$
since  $\xi\left(  \theta\right)  =\left(  1-\theta\right)  ^{-1}$ and
$\mu\left(  \theta\right)  =\sigma\left(  1-\theta\right)  ^{-1} = \sigma \xi\left(  \theta\right)$.
}%
{Recall}
$\tilde{u}_{3,m}=\min_{\tau\in\left\{  a,b\right\}  }\min_{\left\{
j:\theta_{j}\neq\theta_{\tau}\right\}  }\left\vert \xi\left(  \theta_{\tau
}\right)  -\xi\left(  \theta_{i}\right)  \right\vert $,
\[
\psi_{1}\left(  t,\theta\right)  ={\int}K_{1}\left(  t,x\right)  dG_{\theta
}\left(  x\right)  =\frac{1}{\pi}\int_{\left(  \mu\left(  \theta\right)
-b\right)  t}^{\left(  \mu\left(  \theta\right)  -a\right)  t}\frac{\sin y}%
{y}dy,
\]
and $u_{m}=\min_{\tau\in\left\{  a,b\right\}  }\min_{\left\{  j:\mu_{j}%
\neq\tau\right\}  }\left\vert \mu_{j}-\tau\right\vert $, {the last of which is defined in \autoref{SpeedOracleBoundedNull}}.
Noticing $\mu_i =\mu\left(  \theta_{i}\right),i=1,\ldots,m$, we have $u_{m}=\sigma\tilde{u}_{3,m}$.
{
Due to the differentiable, monotonic, bijection, $\mu\left(  \theta\right)  =\sigma\left(  1-\theta\right)  ^{-1}$, between the mean
$\mu=\mu\left(  \theta\right)  $ and natural parameter $\theta=\theta\left(
\mu\right)  $ for the Gamma family,
\autoref{SpeedOracleBoundedNull}} implies, when $tu_{m} = {t \sigma \tilde{u}_{3,m}} \geq2$ and $t\left(  b-a\right)  \geq2$,%
\begin{equation}
\left\vert \frac{\varphi_{m}\left(  t,\boldsymbol{\theta}\right)  }{\pi_{1,m}%
}-1\right\vert \leq C\left(  \frac{1}{t\pi_{1,m}}+\frac{1}{t\tilde{u}_{3,m}%
\pi_{1,m}}\right)  \text{ \ for }t{\sigma}\tilde{u}_{3,m}\geq
2.\label{eqGammaBoundedNullB}%
\end{equation}
So, $\pi_{1,m}%
^{-1}\varphi_{m}\left(  t_m,\boldsymbol{\theta}\right)  \rightarrow1$ if
$t_m^{-1}\left(  1+\tilde{u}_{3,m}^{-1}\right)  =o\left(  \pi_{1,m}\right)  $.
Therefore, by (\ref{eq18c}) a uniform consistency class when $\sigma \ge 3/4$ is%
\[
\mathcal{Q}\left(  \mathcal{F}\right)  =\left\{
\begin{array}
[c]{c}%
t_m=4^{-1}(\max\left\{  \sigma,1\right\})^{-1}\gamma u_{3,m}\ln m,t_m^{-1}\left(  1+\tilde{u}_{3,m}%
^{-1}\right)  =o\left(  \pi_{1,m}\right)  ,\\
t_m\rightarrow\infty,\left\Vert 1-\boldsymbol{\theta}\right\Vert _{\infty
}^{\sigma-3/4}t_m^{11/4-\sigma}=o\left(  m^{1-\gamma}\pi_{1,m}^{2}\right)
\end{array}
\right\}
\]
for each $\gamma\in \left(  0,1\right)$, for which $\gamma=1$ can be set when $\sigma>11/4$, and by (\ref{eq18b}) a uniform consistency class
when $\sigma\leq3/4$ is%
\[
\mathcal{Q}\left(  \mathcal{F}\right)  =\left\{
\begin{array}
[c]{c}%
 t_m=4^{-1}\gamma u_{3,m}\ln m,t_m^{-1}\left(  1+\tilde{u}_{3,m}^{-1}\right)
=o\left(  \pi_{1,m}\right)  ,\\
t_m\rightarrow\infty,\left(  \ln m\right)  ^{11/4-\sigma}u_{3,m}%
^{2}=o\left(  m^{1-\gamma}\pi_{1,m}^{2}\right)
\end{array}
\right\}
\]
for each $\gamma\in\left(  0,1\right)  $.\qed

\section{Proofs Related to Construction II}

\label{AppProofsNEFB}

\subsection{Proof of \autoref{TmVNEF}}

Recall $\tilde{c}_{n}\left(  \theta\right)  =\int x^{n}dG_{\theta}\left(
x\right)  =\zeta_{0}\xi^{n}\left(  \theta\right)  \tilde{a}_{n}$ for the
constant $\zeta_{0}=1$ and $\mu\left(  \theta\right)  =\tilde{c}_{1}$. Set%
\[
K_{4,0}^{\dagger}\left(  t,x\right)  =\frac{1}{2\pi\zeta_{0}}\int_{0}%
^{t}dy\int_{-1}^{1}\exp\left(  -\iota ysb\right)  \sum_{n=0}^{\infty}%
\frac{\left(  \iota ys\right)  ^{n}\left(  \zeta_{0}\tilde{a}_{1}x\right)
^{n+1}}{\tilde{a}_{n+1}n!}ds.
\]
Then%
\begin{equation}
K_{4,0}^{\dagger}\left(  t,x\right)  =\frac{1}{2\pi\zeta_{0}}\int_{0}%
^{1}tdy\int_{-1}^{1}\exp\left(  -\iota tysb\right)  \sum_{n=0}^{\infty}%
\frac{\left(  \iota tys\right)  ^{n}\left(  \zeta_{0}\tilde{a}_{1}x\right)
^{n+1}}{\tilde{a}_{n+1}n!}ds. \label{eqk40}
\end{equation}
Further,%
\begin{align*}
\int K_{4,0}^{\dagger}\left(  t,x\right)  dG_{\theta}\left(  x\right)   &
=\frac{\zeta\left(  \theta\right)  }{2\pi\zeta_{0}}\int_{0}^{t}dy\int_{-1}%
^{1}\exp\left(  -\iota ysb\right)  \sum_{n=0}^{\infty}\frac{\left(  \iota
ys\right)  ^{n}}{n!}\left(  \zeta_{0}\tilde{a}_{1}\right)  ^{n+1}\xi
^{n+1}\left(  \theta\right)  ds\\
&  =\frac{1}{2\pi}\int_{0}^{t}\mu\left(  \theta\right)  dy\int_{-1}^{1}%
\exp\left(  -\iota ysb\right)  \exp\left(  \iota ys\mu\left(  \theta\right)
\right)  ds\\
&  =\frac{1}{2\pi}\int_{0}^{t}\mu\left(  \theta\right)  dy\int_{-1}^{1}%
\exp\left[  \iota ys\left\{  \mu\left(  \theta\right)  -b\right\}  \right]
ds.
\end{align*}
On the other hand, set%
\[
K_{4,1}^{\dagger}\left(  t,x\right)  =-\frac{1}{2\pi\zeta_{0}}\int_{0}%
^{t}dy\int_{-1}^{1}b\exp\left(  -\iota ysb\right)  \sum_{n=0}^{\infty}%
\frac{\left(  \iota ys\right)  ^{n}\left(  \zeta_{0}\tilde{a}_{1}x\right)
^{n}}{\tilde{a}_{n}n!}ds.
\]
Then%
\begin{equation}
K_{4,1}^{\dagger}\left(  t,x\right)  =-\frac{1}{2\pi\zeta_{0}}\int_{0}%
^{1}tdy\int_{-1}^{1}b\exp\left(  -\iota tysb\right)  \sum_{n=0}^{\infty}%
\frac{\left(  \iota tys\right)  ^{n}\left(  \zeta_{0}\tilde{a}_{1}x\right)
^{n}}{\tilde{a}_{n}n!}ds. \label{eqk41}
\end{equation}
Further,%
\begin{align*}
\int K_{4,1}^{\dagger}\left(  t,x\right)  dG_{\theta}\left(  x\right)   &
=-\frac{b\zeta\left(  \theta\right)  }{2\pi\zeta_{0}}\int_{0}^{t}dy\int%
_{-1}^{1}\exp\left(  -\iota ysb\right)  \sum_{n=0}^{\infty}\frac{\left(  \iota
ys\right)  ^{n}}{n!}\left(  \zeta_{0}\tilde{a}_{1}\right)  ^{n}\xi^{n}\left(
\theta\right)  ds\\
&  =-\frac{b}{2\pi}\int_{0}^{t}dy\int_{-1}^{1}\exp\left(  -\iota ysb\right)
\exp\left(  \iota ys\mu\left(  \theta\right)  \right)  ds\\
&  =-\frac{b}{2\pi}\int_{0}^{t}dy\int_{-1}^{1}\exp\left[  \iota ys\left\{
\mu\left(  \theta\right)  -b\right\}  \right]  ds.
\end{align*}

Set $K_{1}^{\dagger}\left(  t,x\right)  =K_{4,0}^{\dagger}\left(  t,x\right)
+K_{4,1}^{\dagger}\left(  t,x\right)  $. Then%
\[
K_{1}^{\dagger}\left(  t,x\right)  =\frac{1}{2\pi\zeta_{0}}\int_{0}^{1}%
tdy\int_{-1}^{1}\exp\left(  -\iota tysb\right)  \sum_{n=0}^{\infty}%
\frac{\left(  \iota tys\right)  ^{n}\left(  \zeta_{0}\tilde{a}_{1}x\right)
^{n}}{n!}\left(  \frac{\zeta_{0}\tilde{a}_{1}x}{\tilde{a}_{n+1}}-\frac
{b}{\tilde{a}_{n}}\right)  ds.
\]
and%
\[
\psi_{1}\left(  t,\theta\right)  =\int K_{1}^{\dagger}\left(  t,x\right)
dG_{\theta}\left(  x\right)  =\frac{1}{2\pi}\int_{0}^{t}\left\{  \mu\left(
\theta\right)  -b\right\}  dy\int_{-1}^{1}\exp\left[  \iota ys\left\{
\mu\left(  \theta\right)  -b\right\}  \right]  ds.
\]
Since $\psi_{1}\left(  t,\theta\right)  $ is real-valued, we also have
$\psi_{1}\left(  t,\theta\right)  =\int K_{1}\left(  t,x\right)  dG_{\theta
}\left(  x\right)  $, where%
\begin{align*}
K_{1}\left(  t,x\right)   &  =\Re\left\{  K_{1}^{\dagger}\left(  t,x\right)
\right\} \\
&  =\frac{1}{2\pi\zeta_{0}}\int_{0}^{1}tdy\int_{-1}^{1}\sum_{n=0}^{\infty}%
\cos\left(  2^{-1}\pi n-tysb\right)  \frac{\left(  tys\right)  ^{n}\left(
\zeta_{0}\tilde{a}_{1}x\right)  ^{n}}{n!}\left(  \frac{\zeta_{0}\tilde{a}%
_{1}x}{\tilde{a}_{n+1}}-\frac{b}{\tilde{a}_{n}}\right)  ds.
\end{align*}
Now set $K\left(  t,x\right)  =2^{-1}-K_{1}\left(  t,x\right)  -2^{-1}%
K_{3,0}\left(  t,x;\theta_{b}\right)  $ with%
\[
K_{3,0}\left(  t,x;\theta_{b}\right)  =\frac{1}{\zeta_{0}}\int_{\left[
-1,1\right]  }\sum_{n=0}^{\infty}\frac{\left(  -tsx\right)  ^{n}\cos\left\{
\frac{\pi}{2}n+ts\xi\left(  \theta_{b}\right)  \right\}  }{\tilde{a}_{n}%
n!}\omega\left(  s\right)  ds
\]
given by \autoref{ThmPoinNull}. Then%
\begin{align*}
\psi\left(  t,\theta\right)  =\int K\left(  t,x\right)  dG_{\theta}\left(
x\right)   &  =2^{-1}-\int_{0}^{t}\left\{  \mu\left(  \theta\right)
-b\right\}  dy\int_{-1}^{1}\exp\left[  \iota ys\left\{  \mu\left(
\theta\right)  -b\right\}  \right]  ds\\
&  -2^{-1}\int_{\left[  -1,1\right]  }\cos\left[  ts\left\{  \xi\left(
\theta_{b}\right)  -\xi\left(  \theta\right)  \right\}  \right]  \omega\left(
s\right)  ds.
\end{align*}
By \autoref{ThmPoinNull} the pair $\left(  K,\psi\right)  $ in (\ref{V-b}) is
as desired.\qed

\subsection{Proof of \autoref{ThmVNEFConsistency}}

We need the following:

\begin{lemma}
\label{Lm:GammaBessel}For a fixed $\sigma>0$, let%
\[
\tilde{w}_{2}\left(  t,x\right)  =\Gamma\left(  \sigma\right)  \sum
_{n=0}^{\infty}\frac{t^{n}}{n!}\frac{x^{n+1}}{\Gamma\left(  \sigma+n+1\right)
}\text{ for }t,x>0\text{.}%
\]
If $Z$ has CDF $G_{\theta}$ from the Gamma family with scale parameter
$\sigma$, then%
\begin{equation}
\mathbb{E}\left[  \tilde{w}_{2}^{2}\left(  z,Z\right)  \right]  \leq
\frac{Cz^{3/4-\sigma}}{\left(  1-\theta\right)  ^{11/4-\sigma}}\exp\left(
\frac{8z/\sqrt{2}}{1-\theta}\right)  \label{eq15b}%
\end{equation}
for positive and sufficiently large $z$.
\end{lemma}

The proof of \autoref{Lm:GammaBessel} is provided in
\autoref{ProofOfGammaMoments}. Now we provide the arguments. First, we obtain an upper bound on $\mathbb{V}\left\{
\hat{\varphi}_{m}\left(  t,\mathbf{z}\right)  \right\}  $. Note that
$\zeta_{0}=1$ and $\tilde{a}_{1}=\sigma$. For $y\in\left[  0,1\right]  $ and
$t,x>0$, define%
\[
w_{3,1}\left(  t,x,y\right)  =\Gamma\left(  \sigma\right)  \sum_{n=0}^{\infty
}\cos\left(  2^{-1}\pi n-tyb\right)  \frac{\left(  ty\right)  ^{n}}{n!}%
\frac{\left(  \sigma x\right)  ^{n+1}}{\Gamma\left(  \sigma+n+1\right)  }%
\]
and%
\[
w_{3,2}\left(  t,x,y\right)  =\Gamma\left(  \sigma\right)  \sum_{n=0}^{\infty
}\cos\left(  2^{-1}\pi n-tyb\right)  \frac{\left(  ty\right)  ^{n}\left(
\sigma x\right)  ^{n}}{n!\Gamma\left(  \sigma+n\right)  }.
\]
Then, uniformly for $s\in\left[  -1,1\right]  $ and $y\in\left[  0,1\right]
$,
\begin{equation}
\left\vert w_{3,1}\left(  ts,x,y\right)  \right\vert \leq\tilde{w}%
_{3,1}\left(  t\sigma,x\right)  =\sigma\Gamma\left(  \sigma\right)  \sum
_{n=0}^{\infty}\frac{\left\vert t\sigma\right\vert ^{n}}{n!}\frac{\left\vert
x\right\vert ^{n+1}}{\Gamma\left(  \sigma+n+1\right)  } \label{eq18d}%
\end{equation}
and%
\begin{equation}
\left\vert w_{3,2}\left(  ts,x,y\right)  \right\vert \leq\tilde{w}%
_{3,2}\left(  t\sigma,x\right)  =\Gamma\left(  \sigma\right)  \sum
_{n=0}^{\infty}\frac{\left\vert t\sigma\right\vert ^{n}\left\vert x\right\vert
^{n}}{n!\Gamma\left(  \sigma+n\right)  }. \label{eq18e}%
\end{equation}
Notice $\tilde{a}_{n}=\Gamma\left(  n+\sigma\right)  /\Gamma\left(
\sigma\right)  $. Recall the functions $K_{4,0}^{\dagger}\left(  t,x\right)  $
and $K_{4,1}^{\dagger}\left(  t,x\right)  $ defined by (\ref{eqk40}) and
(\ref{eqk41}) in the proof of \autoref{TmVNEF} such that $K_{1}\left(
t,x\right)  =\Re\left\{  K_{4,0}^{\dagger}\left(  t,x\right)  \right\}
+\Re\left\{  K_{4,1}^{\dagger}\left(  t,x\right)  \right\}  $. Let
$K_{4,0}\left(  t,x\right)  =\Re\left\{  K_{4,0}^{\dagger}\left(  t,x\right)
\right\}  $ and $K_{4,1}\left(  t,x\right)  =\Re\left\{  K_{4,1}^{\dagger
}\left(  t,x\right)  \right\}  $. Then%
\[
K_{4,0}\left(  t,x\right)  =\frac{1}{2\pi}\int_{0}^{1}tdy\int_{-1}^{1}%
w_{3,1}\left(  ts,x,y\right)  dy
\]
and%
\[
K_{4,1}\left(  t,x\right)  =\frac{-b}{2\pi}\int_{0}^{1}tdy\int_{-1}^{1}%
w_{3,2}\left(  ts,x,y\right)  dy.
\]

Set $\hat{S}_{3,m,1}\left(  ts,y,\mathbf{z}\right)  =m^{-1}\sum_{i=1}^{m}w_{3,1}\left(
ts,z_i,y\right)  $, $\hat{S}_{3,m,2}\left(  ts,y,\mathbf{z}\right)  =-bm^{-1}\sum_{i=1}%
^{m}w_{3,2}\left(  ts,z_i,y\right)  $ and%
\[
\hat{S}_{3,m}\left(  ts,y,\mathbf{z}\right)  =\hat{S}_{3,m,1}\left(  ts,y,\mathbf{z}\right)
{+}\hat{S}_{3,m,2}\left(  ts,y,\mathbf{z}\right)  .
\]
Recall
$\bigskip\hat{\varphi}_{1,m}\left(  t,\mathbf{z}\right)  =m^{-1}\sum_{i=1}%
^{m}K_{1}\left(  t,z_{i}\right)  $ and\ $\varphi_{1,m}\left(
t,\boldsymbol{\theta}\right)  =m^{-1}\sum_{i=1}^{m}\mathbb{E}\left\{
K_{1}\left(  t,z_{i}\right)  \right\}  $.
Then
\begin{align*}
\hat{\varphi}_{1,m}\left(  t,\mathbf{z}\right)   &  =m^{-1}\sum_{i=1}^{m}%
K_{1}\left(  t,z_{i}\right)  =m^{-1}\sum_{i=1}^{m}K_{4,0}\left(
t,z_{i}\right)  +m^{-1}\sum_{i=1}^{m}K_{4,1}\left(  t,z_{i}\right)  \\
&  =m^{-1}\sum_{i=1}^{m}\frac{1}{2\pi}\int_{0}^{1}tdy\int_{-1}^{1}%
w_{3,1}\left(  ts,z_{i},y\right)  dy+m^{-1}\sum_{i=1}^{m}\frac{-b}{2\pi}%
\int_{0}^{1}tdy\int_{-1}^{1}w_{3,2}\left(  ts,x,y\right)  dy\\
&  =\frac{t}{2\pi}\int_{0}^{1}dy\int_{-1}^{1}\left[  \hat{S}_{3,m,1}\left(
ts,y,\mathbf{z}\right)  dy+\hat{S}_{3,m,2}\left(  ts,y,\mathbf{z}\right)
\right]  ds\\
&  =\frac{t}{2\pi}\int_{0}^{1}dy\int_{-1}^{1}\hat{S}_{3,m}\left(
ts,y,\mathbf{z}\right)  ds,
\end{align*}
and setting%
\[
\left\{
\begin{array}
[c]{l}%
\Delta_{3,m,j}\left(  ts,y,\mathbf{z}\right)  =\hat{S}_{3,m,j}\left(
ts,y,\mathbf{z}\right)  -\mathbb{E}\left(  \hat{S}_{3,m,j}\left(
ts,y,\mathbf{z}\right)  \right)  ,j=1,2\\
\Delta_{3,m}\left(  ts,y,\mathbf{z}\right)  =\hat{S}_{3,m}\left(
ts,y,\mathbf{z}\right)  -\mathbb{E}\left(  \hat{S}_{3,m}\left(
ts,y,\mathbf{z}\right)  \right)
\end{array}
\right.
\]
gives
\begin{align*}
\hat{\varphi}_{1,m}\left(  t,\mathbf{z}\right)  -\varphi_{1,m}\left(
t,\boldsymbol{\theta}\right)    & =\frac{t}{2\pi}\int_{0}^{1}dy\int_{-1}%
^{1}\Delta_{3,m,1}\left(  ts,y,\mathbf{z}\right)  ds+\frac{t}{2\pi}\int%
_{0}^{1}dy\int_{-1}^{1}\Delta_{3,m,2}\left(  ts,y,\mathbf{z}\right)  ds\\
& =\frac{t}{2\pi}\int_{0}^{1}dy\int_{-1}^{1}\Delta_{3,m}\left(
ts,y,\mathbf{z}\right)  ds.
\end{align*}

Therefore, using the same technique that obtained (\ref{eqdx2}), we get%
\begin{align}
\mathbb{V}\left\{  \hat{\varphi}_{1,m}\left(  t,\mathbf{z}\right)  \right\}
&  \leq\frac{t^{2}}{2\pi^{2}}\mathbb{E}\left(  \int_{0}^{1}dy\int_{-1}%
^{1}\Delta_{3,m}^{2}\left(  ts,y,\mathbf{z}\right)  ds\right)  \nonumber\\
&  =\frac{t^{2}}{2\pi^{2}}\int_{0}^{1}dy\int_{-1}^{1}\mathbb{E}\left(
\Delta_{3,m}^{2}\left(  ts,y,\mathbf{z}\right)  \right)  ds\nonumber\\
&  \leq\frac{t^{2}}{\pi^{2}}\int_{0}^{1}dy\int_{-1}^{1}\left[  \mathbb{E}%
\left(  \Delta_{3,m,1}^{2}\left(  ts,y,\mathbf{z}\right)  \right)
+\mathbb{E}\left(  \Delta_{3,m,2}^{2}\left(  ts,y,\mathbf{z}\right)  \right)
\right]  ds,\label{eq18a}%
\end{align}
where to obtain the first inequality in (\ref{eq18a}) we have used the fact
(due to H\~{o}lder's inequality)%
\begin{align}
&  \mathbb{E}\left[  \left(  \int_{a_{1}}^{b_{1}}dy\int_{a_{2}}^{b_{2}%
}\left\vert X\left(  s,y\right)  \right\vert ds\right)  ^{2}\right]
\nonumber\\
&  \leq\mathbb{E}\left[  \left(  \int_{a_{1}}^{b_{1}}dy\sqrt{b_{2}-a_{2}%
}\left[  \int_{a_{2}}^{b_{2}}\left\vert X\left(  s,y\right)  \right\vert
^{2}ds\right]  ^{1/2}\right)  ^{2}\right]  \nonumber\\
&  \leq\prod\nolimits_{j=1}^{2}\left(  b_{j}-a_{j}\right)  \mathbb{E}\left(
\int_{a_{1}}^{b_{1}}dy\int_{a_{2}}^{b_{2}}\left\vert X\left(  s,y\right)
\right\vert ^{2}ds\right) \nonumber\\
& =\prod\nolimits_{j=1}^{2}\left(  b_{j}%
-a_{j}\right)  \left(  \int_{a_{1}}^{b_{1}}dy\int_{a_{2}}^{b_{2}}%
\mathbb{E}\left(  \left\vert X\left(  s,y\right)  \right\vert ^{2}\right)
ds\right)  \label{eqFact}%
\end{align}
for a random variable $X\left(  s,y\right)  $ with parameters $\left(
s,y\right)  $ and finite constants $a_{j}<b_{j}$ with $j=1,2$,
 and to obtain the second inequality in (\ref{eq18a}) we have
used the fact $\left(  a_{\ast}+b_{\ast}\right)  ^{2}\leq2a_{\ast}%
^{2}+2b_{\ast}^{2}$ for $a_{\ast},b_{\ast}\in\mathbb{R}$. Note that
$\mathbb{V}\left(  \hat{S}_{3,m,j}\left(  ts,y,\mathbf{z}\right)  \right)
=\mathbb{E}\left(  \Delta_{3,m,j}^{2}\left(  ts,y,\mathbf{z}\right)  \right)
$ for $j=1,2$.
By the inequalities (\ref{eq18d}), (\ref{eq18e}), (\ref{eq15}) and
\autoref{Lm:GammaBessel}, we have, for $t>0$ sufficiently large,%
\begin{align}
\mathbb{V}\left(  \hat{S}_{3,m,1}\left(  ts,y,\mathbf{z}\right)  \right) &
\leq\frac{1}{m^{2}}\sum_{i=1}^{m}\mathbb{E}\left\{  \tilde{w}_{3,1}^{2}\left(
t\sigma,z_{i}\right)  \right\}  \leq\frac{C}{m^{2}}\sum_{i=1}^{m}%
\frac{t^{3/4-\sigma}}{\left(  1-\theta_{i}\right)  ^{11/4-\sigma}}\exp\left(
\frac{8\sigma t/\sqrt{2}}{1-\theta_{i}}\right) \nonumber\\
&  \leq V_{3,1,m}=\frac{C}{m^{2}}\exp\left(  \frac{8\sigma t/\sqrt{2}}%
{u_{3,m}}\right)  \sum_{i=1}^{m}\frac{t^{3/4-\sigma}}{\left(  1-\theta
_{i}\right)  ^{11/4-\sigma}} \label{eq17}%
\end{align}
and%
\begin{align}
\mathbb{V}\left(  \hat{S}_{3,m,2}\left(  ts,y,\mathbf{z}\right)  \right)
  &  \leq\frac{b^{2}}{m^{2}}\sum_{i=1}^{m}\mathbb{E}\left\{
\tilde{w}_{3,2}^{2}\left(  t\sigma,z_{i}\right)  \right\}  \leq\frac{b^{2}%
}{m^{2}}\sum_{i=1}^{m}\frac{t^{3/4-\sigma}}{\left(  1-\theta_{i}\right)
^{3/4-\sigma}}\exp\left(  \frac{4\sigma t}{1-\theta_{i}}\right) \nonumber\\
&  \leq V_{3,2,m}=\frac{C}{m^{2}}\exp\left(  \frac{4\sigma t}{u_{3,m}}\right)
\sum_{i=1}^{m}\frac{t^{3/4-\sigma}}{\left(  1-\theta_{i}\right)  ^{3/4-\sigma
}}, \label{eq17b}%
\end{align}
where $u_{3,m}=\min_{1\leq i\leq m}\left\{  1-\theta_{i}\right\}  $. Combining
(\ref{eq18a}), (\ref{eq17}) and (\ref{eq17b}) gives, for $t>0$ sufficiently large,%
\[
\mathbb{V}\left\{  \hat{\varphi}_{1,m}\left(  t,\mathbf{z}\right)  \right\}
\leq\frac{Ct^{11/4-\sigma}}{m^{2}}\exp\left(  \frac{4\sqrt{2}\sigma t}%
{u_{3,m}}\right)  \sum_{i=1}^{m}l\left(  \theta_{i},\sigma\right)
\]
where%
\begin{equation}
l\left(  \theta_{i},\sigma\right)  =\max\left\{  \left(  1-\theta_{i}\right)
^{\sigma-11/4},\left(  1-\theta_{i}\right)  ^{\sigma-3/4}\right\}  .
\label{eq17e}%
\end{equation}

Recall%
\[
K_{3,0}\left(  t,x;\theta_{b}\right)  =\frac{\Gamma\left(  \sigma\right)
}{\zeta_{0}}\int_{\left[  -1,1\right]  }\sum_{n=0}^{\infty}\dfrac{\left(
-tsx\right)  ^{n}\cos\left\{  2^{-1}\pi n+ts\xi\left(  \theta_{b}\right)
\right\}  }{n!\Gamma\left(  n+\sigma\right)  }\omega\left(  s\right)  ds.
\]
and $\hat{\varphi}_{3,0,m}\left(  t,\mathbf{z};\theta_{b}\right)  =m^{-1}%
\sum_{i=1}^{m}K_{3,0}\left(  t,z_{i};\theta_{b}\right)  $ and $\varphi
_{3,0,m}\left(  t,\boldsymbol{\theta};\tau\right)  =\mathbb{E}\left\{
\hat{\varphi}_{3,0,m}\left(  t,\mathbf{z};\theta_{b}\right)  \right\}  $. Then
Theorem 8 of \cite{Chen:2018a} asserts, for $t>0$ sufficiently large,%
\[
\mathbb{V}\left\{  \hat{\varphi}_{3,0,m}\left(  t,\mathbf{z};\theta
_{b}\right)  \right\}  \leq\frac{C}{m^{2}}\exp\left(  \frac{4t}{u_{3,m}%
}\right)  \sum_{i=1}^{m}\frac{t^{3/4-\sigma}}{\left(  1-\theta_{i}\right)
^{3/4-\sigma}}.
\]
Recall $K\left(  t,x\right)  =2^{-1}-K_{1}\left(  t,x\right)  -2^{-1}%
K_{3,0}\left(  t,x;\theta_{b}\right)  $. Then, for $t>0$ sufficiently large,
\begin{align}
\mathbb{V}\left\{  \hat{\varphi}_{m}\left(  t,\mathbf{z}\right)  \right\}   &
\leq2\mathbb{V}\left\{  \hat{\varphi}_{1,m}\left(  t,\mathbf{z}\right)
\right\}  + 2^{-1}\mathbb{V}\left\{  \hat{\varphi}_{3,0,m}\left(  t,\mathbf{z}%
;\theta_{b}\right)  \right\} \nonumber\\
&  \leq V_{3,m}=\frac{Ct^{11/4-\sigma}}{m^{2}}\exp\left(  \frac{4t\max\left\{
1,\sqrt{2}\sigma\right\}  }{u_{3,m}}\right)  \sum_{i=1}^{m}l\left(  \theta
_{i},\sigma\right)  , \label{eq17d}%
\end{align}
and%
\begin{equation}
V_{3,m}\leq V_{3,m}^{\ast}=\frac{Ct^{11/4-\sigma}}{m^{2}u_{3,m}^{2}}%
\exp\left(  \frac{4t\max\left\{  1,\sqrt{2}\sigma\right\}  }{u_{3,m}}\right)
\sum_{i=1}^{m}\left(  1-\theta_{i}\right)  ^{\sigma-3/4}\label{eq17ddx}%
\end{equation}
since $l\left(  \theta_{i},\sigma\right)  $ in (\ref{eq17e}) is upper bounded
by $Cu_{3,m}^{-2}\left(  1-\theta_{i}\right)  ^{\sigma-3/4}$ regardless of
whether $\liminf_{m\rightarrow\infty}u_{3,m}=0$ or not.

Secondly, we provide a uniform consistency class.
When $\sigma\geq3/4$, then
(\ref{eq17d}) and (\ref{eq17ddx}) imply
\[
V_{3,m}^{\ast}\leq V_{3,m}^{\ast\dag}=\frac{Ct^{11/4-\sigma}}{mu_{3,m}^{2}%
}\exp\left(  \frac{4t\max\left\{  1,\sqrt{2}\sigma\right\}  }{u_{3,m}}\right)
\left\Vert 1-\boldsymbol{\theta}\right\Vert _{\infty}^{\sigma-3/4},
\]
and that setting $t_{m}=\left(  4\max\left\{  1,\sqrt{2}\sigma\right\}  \right)
^{-1}u_{3,m}\gamma\ln m$ for any fixed $\gamma\in\left(  0,1\right)  $ gives%
\[
\Pr\left\{  \left\vert \frac{\hat{\varphi}_{m}\left(  t_{m},\mathbf{z}\right)
-\varphi_{m}\left(  t_{m},\boldsymbol{\theta}\right)  }{\pi_{1,m}}\right\vert
\geq\varepsilon\right\}  \leq\frac{Cu_{3,m}^{3/4-\sigma}\left(  \ln
m\right)  ^{11/4-\sigma}}{\pi_{1,m}^{2}m^{1-\gamma}\varepsilon^{2}}\left\Vert
1-\boldsymbol{\theta}\right\Vert _{\infty}^{\sigma-3/4},
\]
both for $t$ and $t_{m}$ sufficiently large.
Note that $\gamma=1$ can be set when $\sigma>11/4$ since $\lim_{m\rightarrow
\infty}\left(  \ln m\right)  ^{11/4-\sigma}=0$ for all such $\sigma$.
In contrast, when $\sigma\leq
3/4$, then (\ref{eq17d}) and (\ref{eq17ddx}) imply%
\[
V_{3,m}^{\ast}\leq\tilde{V}_{3,m}^{\ast\ast}=\frac{Ct^{11/4-\sigma}}%
{mu_{3,m}^{2}}\exp\left(  \frac{4t\max\left\{  1,\sqrt{2}\sigma\right\}
}{u_{3,m}}\right)  u_{3,m}^{\sigma-3/4},
\]
and that choosing the same sequence $t_{m}\ $for any fixed $\gamma\in\left(
0,1\right)  $ gives%
\[
\Pr\left\{  \left\vert \frac{\hat{\varphi}_{m}\left(  t_{m},\mathbf{z}\right)
-\varphi_{m}\left(  t_{m},\boldsymbol{\theta}\right)  }{\pi_{1,m}}\right\vert
\geq\varepsilon\right\}  \leq\frac{C\left(  \ln m\right)  ^{11/4-\sigma
}}{\pi_{1,m}^{2}m^{1-\gamma}\varepsilon^{2}},
\]
both for $t$ and $t_{m}$ sufficiently large.

Recall
\[
\psi_{3,0}\left(  t,\theta;\theta_{b}\right)  =\int_{\left[  -1,1\right]
}\cos\left[  ts\left\{  \xi\left(  \theta_{b}\right)  -\xi\left(
\theta\right)  \right\}  \right]  \omega\left(  s\right)  ds,
\]
{
which is exactly
\[
\tilde{\psi}_{1,0}\left(  t,\mu;\mu^{\prime}\right)  ={\int_{\left[  -1,1\right]  }%
}\omega\left(  s\right)  \cos\left\{  ts \sigma^{-1} \left(  \mu-\mu^{\prime}\right)
\right\}  ds
\]
that is defined by (\ref{NewPsiPoint}) (in \autoref{SpeedOracleBoundedNull}) but evaluated at $\mu^{\prime}=b$
since  $\xi\left(  \theta\right)  =\left(  1-\theta\right)  ^{-1}$ and
$\mu\left(  \theta\right)  =\sigma\left(  1-\theta\right)  ^{-1} = \sigma \xi\left(  \theta\right)$.
}%
{Also recall}
\[
\left\{
\begin{array}
[c]{l}%
\psi\left(  t,\theta\right)  =2^{-1}-\psi_{1}\left(  t,\theta\right)
-2^{-1}\psi_{3,0}\left(  t,\theta;\theta_{b}\right)  \\
\psi_{1}\left(  t,\theta\right)  =\frac{1}{\pi}\int_{0}^{t}\frac{\sin\left\{
\left(  \mu\left(  \theta\right)  -b\right)  y\right\}  }{y}dy
\end{array}
\right.  ,
\]
$\check{u}_{3,m}=\min_{\left\{  j:\theta_{j}\neq\theta_{b}\right\}
}\left\vert \xi\left(  \theta_{b}\right)  -\xi\left(  \theta_{j}\right)
\right\vert $, and $\tilde{u}_{m}\left(  b\right)  =\min_{\left\{  j:\mu
_{j}\neq b\right\}  }\left\vert \mu_{j}-b\right\vert $ ({the last of which is defined in \autoref{SpeedOneSidedNull}}), where $\mu_i = \mu \left(\theta_i\right),i=1,\ldots,m$. We have $\tilde{u}_{m}\left(  b\right)=\sigma\check{u}_{3,m}$.
{
Due to the differentiable, monotonic, bijection, $\mu\left(  \theta\right)  =\sigma\left(  1-\theta\right)  ^{-1} $, between the mean
$\mu=\mu\left(  \theta\right)  $ and natural parameter $\theta=\theta\left(
\mu\right)  $ for the Gamma family, \autoref{SpeedOneSidedNull}} implies,
when $t\tilde{u}_{m}\left(  b\right) = {t\sigma\check{u}_{3,m}} \geq2$,
\[
\left\vert \pi_{1,m}^{-1}\varphi_{m}\left(  t,\boldsymbol{\theta}\right)
-1\right\vert \leq\frac{C}{t\check{u}_{3,m}\pi_{1,m}}.
\]
So, $\pi_{1,m}^{-1}\varphi_{m}\left(  t_{m},\boldsymbol{\theta
}\right)  \rightarrow1$ when $t_{m}^{-1}\check{u}_{3,m}^{-1}
=o\left(  \pi_{1,m}\right)  $. Therefore, a uniform consistency class is%
\[
\mathcal{Q}\left(  \mathcal{F}\right)  =\left\{
\begin{array}
[c]{c}%
t_{m}=\left(  4\max\left\{  1,\sqrt{2}\sigma\right\}  \right)  ^{-1}%
u_{3,m}\gamma\ln m,t_{m}^{-1}\check{u}_{3,m}^{-1}=o\left(  \pi_{1,m}\right)
,\\
t_{m}\rightarrow\infty,u_{3,m}^{3/4-\sigma}\left(  \ln m\right)
^{11/4-\sigma}\left\Vert 1-\boldsymbol{\theta}\right\Vert _{\infty}%
^{\sigma-3/4}=o\left(  \pi_{1,m}^{2}m^{1-\gamma}\right)
\end{array}
\right\}
\]
when $\sigma\geq3/4$ for each $\gamma\in\left(  0,1\right)  $, for which $\gamma=1$ can be set when $\sigma>11/4$,, and it is%
\[
\mathcal{Q}\left(  \mathcal{F}\right)  =\left\{
\begin{array}
[c]{c}%
t_{m}=\left(  4\max\left\{  1,\sqrt{2}\sigma\right\}  \right)  ^{-1}%
u_{3,m}\gamma\ln m,t_{m}^{-1}\check{u}_{3,m}^{-1}=o\left(  \pi_{1,m}\right)
,\\
t_{m}\rightarrow\infty,\left(  \ln m\right)  ^{11/4-\sigma}=o\left(
\pi_{1,m}^{2}m^{1-\gamma}\right)
\end{array}
\right\}
\]
when $\sigma\leq3/4$ for each $\gamma\in\left(  0,1\right)  $.\qed

\subsection{Proof of \autoref{Lm:GammaBessel}}

\label{ProofOfGammaMoments}

Recall (\ref{eq15a}), i.e.,%
\[
\tilde{w}\left(  z,x\right)  =\sum_{n=0}^{\infty}\frac{\left(  zx\right)
^{n}}{n!\Gamma\left(  \sigma+n\right)  }\text{ \ for }z,x>0.
\]
From the proof of Lemma 4 of \cite{Chen:2018a}, we have%
\[
\tilde{w}\left(  z,x\right)  =\left(  zx\right)  ^{\frac{1}{4}-\frac{\sigma
}{2}}\exp\left(  2\sqrt{zx}\right)  \left[  1+O\left\{  \left(  zx\right)
^{-1}\right\}  \right]
\]
when $zx\rightarrow\infty$. So, when $zx\rightarrow\infty$,%
\begin{align*}
\tilde{w}_{2}\left(  z,x\right)   &  =\Gamma\left(  \sigma\right)  \sum
_{n=0}^{\infty}\frac{z^{n}}{n!}\frac{x^{n+1}}{\Gamma\left(  \sigma+n+1\right)
}\\
&  \leq\Gamma\left(  \sigma\right)  x\left(  zx\right)  ^{\frac{1}{4}%
-\frac{\sigma}{2}}\exp\left(  2\sqrt{zx}\right)  \left[  1+O\left\{  \left(
zx\right)  ^{-1}\right\}  \right]  .
\end{align*}
Note $f_{\theta}\left(  x\right)  \leq C\left(  1-\theta\right)
^{\sigma}x^{\sigma-1}$ for all $\theta<1$ and $x>0$.
Pick a constant $A>0$ such that
\[
\tilde{w}_{2}\left(  z,x\right)  \leq 2 {\Gamma\left(  \sigma\right)} x\left(  zx\right)  ^{\frac{1}{4}%
-\frac{\sigma}{2}}\exp\left(  2\sqrt{zx}\right)  \text{ for all }zx>A,
\]
and define $A_{1,z}=\left[  0,Az^{-1}\right]  $
and $A_{2,z}=\left(  Az^{-1},\infty\right)  $ for each fixed $z>0$.
Then, $\tilde{w}\left(
z,x\right)  \leq Ce^{zx}=O\left(  1\right)  $ and $\tilde{w}_{2}\left(
z,x\right)  \leq x\Gamma\left(  \sigma\right)  e^{zx} \le C  x  $ on
the set $A_{1,z}$.
Therefore,%
\begin{equation}
\int_{A_{1,z}}\tilde{w}_{2}^{2}\left(  z,x\right)  dG_{\theta}\left(
x\right)  \leq C\left(  1-\theta\right)  ^{\sigma}\int_{A_{1,z}}x^{2}%
x^{\sigma-1}dx\leq C\left(  1-\theta\right)  ^{\sigma}z^{-\left(
\sigma+2\right)  }. \label{eq16f}%
\end{equation}
On the other hand,
\begin{align}
\int_{A_{2,z}}\tilde{w}_{2}^{2}\left(  z,x\right)  dG_{\theta}\left(
x\right)   &  \leq C\int_{A_{2,z}}x^{2}\left(  zx\right)  ^{\frac{1}{2}%
-\sigma}\exp\left(  4\sqrt{zx}\right)  dG_{\theta}\left(  x\right) \nonumber\\
&  =C\int_{A_{2,z}}x^{2}\left(  zx\right)  ^{\frac{1}{2}-\sigma}\sum
_{n=0}^{\infty}\frac{\left(  4\sqrt{zx}\right)  ^{n}}{n!}dG_{\theta}\left(
x\right)  = {C}z^{\frac{1}{2}-\sigma}B_{3}\left(  z\right)  , \label{eq16a}%
\end{align}
where%
\[
B_{3}\left(  z\right)  =\sum_{n=0}^{\infty}\frac{4^{n}z^{n/2}}{n!}\tilde
{c}_{2^{-1}\left(  n+5\right)  }^{\ast}\text{ \ and \ }\tilde{c}%
_{2^{-1}\left(  n+5\right)  }^{\ast}=\int x^{2^{-1}\left(  n+5\right)
-\sigma}dG_{\theta}\left(  x\right)  .
\]

By the formula,%
\[
\frac{\left(  1-\theta\right)  ^{\sigma}}{\Gamma\left(  \sigma\right)  }%
\int_{0}^{\infty}x^{\beta}e^{\theta x}x^{\sigma-1}e^{-x}dx=\frac{\Gamma\left(
\beta+\sigma\right)  }{\Gamma\left(  \sigma\right)  }\frac{\left(
1-\theta\right)  ^{\sigma}}{\left(  1-\theta\right)  ^{\beta+\sigma}}\text{
\ for }\alpha,\beta>0,
\]
we have%
\[
\tilde{c}_{2^{-1}\left(  n+5\right)  }^{\ast}=\frac{\Gamma\left(
2^{-1}n+2^{-1}\times5\right)  }{\Gamma\left(  \sigma\right)  }\frac{\left(
1-\theta\right)  ^{\sigma-\frac{5}{2}}}{\left(  1-\theta\right)  ^{2^{-1}n}}.
\]
By Theorem 1 of \cite{Karatsuba:2001} that implies ``Ramanujan's double inequality'' as%
\[
\left(  8x^{3}+4x^{2}+x+\frac{1}{100}\right)  ^{1/6}<\frac{\Gamma\left(
x+1\right)  }{\sqrt{\pi}\left(  \frac{x}{e}\right)  ^{x}}<\left(
8x^{3}+4x^{2}+x+\frac{1}{30}\right)  ^{1/6}\text{ for }x\geq1
\]
and which implies Stirling's formula,%
\begin{align*}
\frac{\Gamma\left(  \frac{n+5}{2}\right)  }{n!}  &  \leq C\frac{\sqrt
{\pi\left(  n+3\right)  }\left(  \frac{n+3}{2}\right)  ^{\frac{n+3}{2}}%
}{e^{\frac{n+3}{2}}\sqrt{2\pi n}\left(  \frac{n}{e}\right)  ^{n}}\leq
Ce^{\frac{n}{2}}2^{-\frac{n}{2}}\frac{\left(  n+3\right)  ^{n/2}}{n^{n/2}%
}\frac{\left(  n+3\right)  ^{3/2}}{n^{n/2}}\\
&  \leq Ce^{\frac{n}{2}}2^{-\frac{n}{2}}\frac{\left(  n+3\right)  ^{7/4}%
}{n^{n/2}}\leq C2^{-\frac{n}{4}}\frac{1}{\sqrt{n!}},\forall n\geq1.
\end{align*}
Therefore,%
\begin{equation}
B_{3}\left(  z\right)  \leq C\left(  1-\theta\right)  ^{\sigma-\frac{5}{2}%
}\sum_{n=0}^{\infty}\frac{4^{n}z^{n/2}2^{-n/4}}{\left(  1-\theta\right)
^{n/2}}\frac{1}{\sqrt{n!}}=C\left(  1-\theta\right)  ^{\sigma-\frac{5}{2}%
}Q^{\ast}\left(  \frac{16z/\sqrt{2}}{1-\theta}\right)  , \label{eq16}%
\end{equation}
where $Q^{\ast}\left(  z\right)  =\sum_{n=0}^{\infty}\frac{z^{n/2}}{\sqrt{n!}%
}$. By definition (8.01) and identity (8.07) in Chapter 8 of \cite{Olver:1974}%
,%
\begin{equation}
Q^{\ast}\left(  z\right)  =\sqrt{2}\left(  2\pi z\right)  ^{1/4}\exp\left(
2^{-1}z\right)  \left\{  1+O\left(  z^{-1}\right)  \right\}  . \label{eq16c}%
\end{equation}
Combining (\ref{eq16a}) through (\ref{eq16c}) gives%
\[
\int_{A_{2,z}}\tilde{w}_{2}^{2}\left(  z,x\right)  dG_{\theta}\left(
x\right)  \leq C\left(  1-\theta\right)  ^{\sigma-\frac{5}{2}}z^{\frac{1}%
{2}-\sigma}\left(  \frac{z}{1-\theta}\right)  ^{1/4}\exp\left(  \frac
{8z/\sqrt{2}}{1-\theta}\right)
\]
for all positive and sufficiently large $z$. Recall (\ref{eq16f}). Thus, when
$1-\theta>0$, $\sigma>0$ and $z$ is positive and sufficiently large,
\begin{align*}
\mathbb{E}\left[  \tilde{w}_{2}^{2}\left(  z,Z\right)  \right]   &  \leq
\int_{A_{1,z}}\tilde{w}_{2}^{2}\left(  z,x\right)  dG_{\theta}\left(
x\right)  +\int_{A_{2,z}}\tilde{w}_{2}^{2}\left(  z,x\right)  dG_{\theta
}\left(  x\right) \\
&  \leq C\left\{  \left(  1-\theta\right)  ^{\sigma}z^{-\left(  \sigma
+2\right)  }+\frac{z^{3/4-\sigma}}{\left(  1-\theta\right)  ^{11/4-\sigma}%
}\exp\left(  \frac{8z/\sqrt{2}}{1-\theta}\right)  \right\} \\
&  \leq\frac{Cz^{3/4-\sigma}}{\left(  1-\theta\right)  ^{11/4-\sigma}}%
\exp\left(  \frac{8z/\sqrt{2}}{1-\theta}\right)  .
\end{align*}
\qed

\section{Proofs Related to the Extension}

\label{AppProofsExt}

\subsection{Proof of \autoref{ThmExtA}}

Recall $\tilde{c}_{n}\left(  \theta\right)
=\int x^{n}dG_{\theta}\left(  x\right)  =\zeta_{0}\xi^{n}\left(
\theta\right)  \tilde{a}_{n}$ and $\mu\left(  \theta\right)  =\zeta_{0}%
\xi\left(  \theta\right)  \tilde{a}_{1}$ and $\zeta_{0}=1$. Define%
\[
K_{1}^{\dag}\left(  t,x\right)  =\frac{t}{2\pi\zeta_{0}}\int_{a}^{b}%
\phi\left(  y\right)    dy\int_{-1}^{1} \exp\left(  -\iota tsy\right)\sum
_{n=0}^{\infty}\frac{\left(  \iota tsx\zeta_{0}\tilde{a}_{1}\right)  ^{n}%
}{\tilde{a}_{n}n!}ds.
\]
Then,%
\begin{align*}
\psi_{1}\left(  t,\theta\right)   &  =\int K_{1}^{\dag
}\left(  t,x\right)  dG_{\theta}\left(  x\right) \\
&  =\frac{t}{2\pi\zeta_{0}}\int_{-1}^{1}\hat{\phi}\left(  ts\right)
\sum_{n=0}^{\infty}\frac{\left(  \iota ts\right)  ^{n}}{\tilde{a}_{n}%
n!}\left(  \zeta_{0}\tilde{a}_{1}\right)  ^{n}\tilde{c}_{n}\left(
\theta\right)  ds\\
&  =\frac{t}{2\pi}\int_{a}^{b}  \phi\left(
y\right)  dy\int_{-1}^{1} \exp\left(  -\iota tsy\right)\sum_{n=0}^{\infty}\frac{\left(  \iota ts\right)
^{n}}{n!}\mu^{n}\left(  \theta\right)  ds\\
&  =\frac{t}{2\pi}\int_{a}^{b}\phi\left(  y\right)  dy\int_{-1}^{1}\exp\left[
\iota ts\left\{  \mu\left(  \theta\right)  -y\right\}  \right]  ds.
\end{align*}
Since $\psi_{1}$ is real, $\psi_{1}=\mathbb{E}\left\{  \Re\left(  K_{1}^{\dag
}\right)  \right\}  $. However,%
\[
K_{1}\left(  t,x\right)  =\Re\left\{  K_{1}^{\dag}\left(  t,x\right)
\right\}  =\frac{t}{2\pi\zeta_{0}}\int_{a}^{b}\phi\left(  y\right)
dy\int_{-1}^{1}\sum_{n=0}^{\infty}\frac{\left(  tsx\zeta_{0}\tilde{a}%
_{1}\right)  ^{n}\cos\left(  2^{-1}n\pi-tsy\right)  }{\tilde{a}_{n}n!}ds.
\]
By \autoref{ThmPoinNull}, the pair $\left(  K,\psi\right)  $ in (\ref{eq2c})
is as desired.\qed

\subsection{Proof of \autoref{ThmFinal}}

The proof uses almost identical arguments as those for the proof of
\autoref{ConcentrationIII}. Take $t>0$.
Recall%
\[
K_{1}\left(  t,x\right)  = \frac{t}{2\pi}\int_{a}^{b}\phi\left(  y\right)
dy\int_{-1}^{1}\sum_{n=0}^{\infty}\frac{\left(  tsx\tilde{a}_{1}\right)
^{n}\cos\left(  2^{-1}n\pi-tsy\right)  }{\tilde{a}_{n}n!}ds
\]
with $\tilde{a}_{1}=\sigma$ and $\tilde{a}_{n}=\Gamma\left(  n+\sigma\right)
/\Gamma\left(  \sigma\right)  $ and%
\[
\psi_{1}\left(  t,\theta\right)  =\int K_{1}\left(  t,x\right)  dG_{\theta
}\left(  x\right)  = \mathcal{D}_{\phi}\left(  t,\mu\left(  \theta\right)
;a,b\right)  =\frac{1}{\pi}\int_{a}^{b}\frac{\sin\left\{  \left(  \mu\left(
\theta\right)  -y\right)  t\right\}  }{\mu\left(  \theta\right)  -y}%
\phi\left(  y\right)  dy.
\]
Take $t>0$ to be sufficiently large.
Recall the following from the proof
of \autoref{ConcentrationIII}:%
\[
w_{1}\left(  t,x,y\right)  =\Gamma\left(  \sigma\right)  \sum_{n=0}^{\infty
}\frac{\left(  tx\sigma\right)  ^{n}\cos\left(  2^{-1}n\pi-ty\right)
}{n!\Gamma\left(  n+\sigma\right)  }\text{ \ for }t\geq0\text{ and }x>0,
\]
and $S_{1,m}\left(  t,y\right)  =m^{-1}\sum_{i=1}^{m}\left[  w_{1}\left(
t,z_{i},y\right)  -\mathbb{E}\left\{  w_{1}\left(  t,z_{i},y\right)  \right\}
\right]  $. Then%
\[
K_{1}\left(  t,x\right)  =\frac{t}{2\pi}\int_{a}^{b}\phi\left(  y\right)
dy\int_{-1}^{1}w_{1}\left(  ts,x,y\right)  ds
\]
and%
\[
\hat{\varphi}_{1,m}\left(  t,\mathbf{z}\right)  -\mathbb{E}\left(
\hat{\varphi}_{1,m}\left(  t,\mathbf{z}\right)  \right)  =\frac{t}{2\pi}%
\int_{a}^{b}\phi\left(  y\right)  dy\int_{-1}^{1}S_{1,m}\left(  ts,y\right)  ds.
\]
So,%
\begin{equation}
\mathbb{V}\left\{  \hat{\varphi}_{1,m}\left(  t,\mathbf{z}\right)  \right\}
\leq\left\Vert \phi\right\Vert _{\infty}^{2}\tilde{V}_{1,m},\label{eqdx2a}%
\end{equation}
where as in the proof of \autoref{ConcentrationIII}
\[
\tilde{V}_{1,m}=\mathbb{E}\left[  \left\{  \frac{1}{2\pi}\int_{a}^{b}%
tdy\int_{-1}^{1}S_{1,m}\left(  ts,y\right)  ds\right\}  ^{2}\right]  \leq
\frac{Ct^{2}}{m}\frac{1}{m}\exp\left(  \frac{4t\sigma}{u_{3,m}}\right)
\sum_{i=1}^{m}\left(  \frac{t}{1-\theta_{i}}\right)  ^{3/4-\sigma}%
\]
and $u_{3,m}=\min_{1\leq i\leq m}\left\{  1-\theta_{i}\right\}  $.

From the proof of \autoref{ConcentrationIII}, recall, for $\tau\in\left\{
a,b\right\}  $,%
\[
K_{3,0}\left(  t,x;\theta_{\tau}\right)  =\Gamma\left(  \sigma\right)
\int_{\left[  -1,1\right]  }\sum_{n=0}^{\infty}\dfrac{\left(  -tsx\right)
^{n}\cos\left\{  2^{-1}\pi n+ts\xi\left(  \theta_{\tau}\right)  \right\}
}{n!\Gamma\left(  n+\sigma\right)  }\omega\left(  s\right)  ds
\]
and $\hat{\varphi}_{3,0,m}\left(  t,\mathbf{z};\tau\right)  =m^{-1}\sum
_{i=1}^{m}K_{3,0}\left(  t,z_{i};\theta_{\tau}\right)  $. Since%
\[
\left\{
\begin{array}
[c]{l}%
K\left(  t,x\right)  =K_{1}\left(  t,x\right)  -2^{-1}\left\{  \phi\left(
a\right)  K_{{3},0}\left(  t,x;\theta_{a}\right)  +\phi\left(
b\right)  K_{{3},0}\left(  t,x;\theta_{b}\right)  \right\}  \\
\psi\left(  t,\mu\right)  =\psi_{1}\left(  t,\mu\right)  -2^{-1}\left\{
\phi\left(  a\right)  \psi_{{3},0}\left(  t,\mu;\theta_{a}\right)
+\phi\left(  b\right)  \psi_{{3},0}\left(  t,\mu;\theta_{b}\right)
\right\}
\end{array}
\right.  ,
\]
then the bound derived in the proof of \autoref{ConcentrationIII}, i.e.,%
\[
\mathbb{V}\left\{  \hat{\varphi}_{3,0,m}\left(  t,\mathbf{z};\tau\right)
\right\}  \leq Cm^{-1}V_{0,m}\text{ \ with \ }V_{0,m}=\frac{1}{m}\exp\left(
\frac{4t}{u_{3,m}}\right)  \sum_{i=1}^{m}\frac{t^{3/4-\sigma}}{\left(
1-\theta_{i}\right)  ^{3/4-\sigma}}%
\]
together with (\ref{eqdx2a}), implies
\begin{equation}
\mathbb{V}\left\{  \hat{\varphi}_{m}\left(  t,\mathbf{z}\right)  \right\}
\leq V_{2,m}^{\dag}=\frac{C\left\Vert \phi\right\Vert _{\infty}^{2}\left(
1+t^{2}\right)  }{m^{2}}\exp\left(  \frac{4t\max\left\{  \sigma,1\right\}
}{u_{3,m}}\right)  \sum_{i=1}^{m}\left(  \frac{t}{1-\theta_{i}}\right)
^{3/4-\sigma}. \label{eq20b}%
\end{equation}
So, when $\sigma\geq 3/4$, (\ref{eq20b}) implies%
\begin{equation}
\Pr\left\{  \frac{\left\vert \hat{\varphi}_{m}\left(  t_m,\mathbf{z}\right)
-\varphi_{m}\left(  t_m,\boldsymbol{\theta}\right)  \right\vert }{\check{\pi
}_{0,m}}\geq\varepsilon\right\}  \leq\frac{C\left\Vert 1-\boldsymbol{\theta
}\right\Vert _{\infty}^{\sigma-3/4}}{\varepsilon^{2}m^{1-\gamma}\check{\pi
}_{0,m}^{2}}\left(  u_{3,m}\ln m\right)  ^{11/4-\sigma} \label{eq19a}%
\end{equation}
by setting $t_m=4^{-1}\left(  \max\left\{  \sigma,1\right\}  \right)^{-1}u_{3,m}\gamma\ln m$ for any fixed $\gamma
\in \left(  0,1\right) $, whereas, when $\sigma\leq3/4$, (\ref{eq20b}) implies%
\begin{equation}
\Pr\left\{  \frac{\left\vert \hat{\varphi}_{m}\left(  t_m,\mathbf{z}\right)
-\varphi_{m}\left(  t_m,\boldsymbol{\theta}\right)  \right\vert }{\check{\pi
}_{0,m}}\geq\varepsilon\right\}  \leq\frac{C\left(  \ln m\right)
^{11/4-\sigma}u_{3,m}^{2}}{\varepsilon^{2}m^{1-\gamma}\check{\pi}_{0,m}^{2}}
\label{eq19b}%
\end{equation}
by setting $t_m=4^{-1}u_{3,m}\gamma\ln m$ for any fixed $\gamma\in\left(
0,1\right)  $,
both for $t_m$ sufficiently large. Note that $\gamma=1$ can be set in (\ref{eq19a}) when $\sigma>11/4$ since $\lim_{m\rightarrow
\infty}\left(  \ln m\right)  ^{11/4-\sigma}=0$ for all such $\sigma$.

Finally, recall $\tilde{u}_{3,m}=\min_{\tau\in\left\{  a,b\right\}
}\min_{\left\{  j:\theta_{j}\neq\theta_{\tau}\right\}  }\left\vert \xi\left(
\theta_{\tau}\right)  -\xi\left(  \theta_{i}\right)  \right\vert $ and for
$\tau\in\left\{  a,b\right\}  $
\[
\psi_{3,0}\left(  t,\theta;\theta_{\tau}\right)  =\int_{\left[  -1,1\right]
}\cos\left[  ts\left\{  \xi\left(  \theta_{\tau}\right)  -\xi\left(
\theta\right)  \right\}  \right]  \omega\left(  s\right)  ds,
\]
{
which is exactly
\[
\tilde{\psi}_{1,0}\left(  t,\mu;\mu^{\prime}\right)  ={\int_{\left[  -1,1\right]  }%
}\omega\left(  s\right)  \cos\left\{  ts \sigma^{-1} \left(  \mu-\mu^{\prime}\right)
\right\}  ds
\]
that is defined by (\ref{NewPsiPoint}) (in \autoref{SpeedOracleBoundedNull}) but evaluated at $\mu^{\prime}=\tau$
since  $\xi\left(  \theta\right)  =\left(  1-\theta\right)  ^{-1}$ and
$\mu\left(  \theta\right)  =\sigma\left(  1-\theta\right)  ^{-1} = \sigma \xi\left(  \theta\right)$.
}%
Further, $u_{m}=\sigma\tilde{u}_{3,m}$. So, {\autoref{SpeedOracleExt}}, i.e.,%
\[
\left\vert \check{\pi}_{0,m}^{-1}\varphi_{m}\left(  t,\boldsymbol{\mu}\right)
-1\right\vert \leq\frac{C}{t\check{\pi}_{0,m}}\left(  1+\left\Vert
\phi\right\Vert _{1,\infty}+\frac{1}{u_{m}%
}\right)
\]
becomes%
\begin{equation}
\left\vert \check{\pi}_{0,m}^{-1}\varphi_{m}\left(  t,\boldsymbol{\theta
}\right)  -1\right\vert \leq\frac{C}{t\check{\pi}_{0,m}}\left(  1+\left\Vert
\phi\right\Vert _{1,\infty}+\frac
{1}{\tilde{u}_{3,m}}\right)  .\label{eq19bx}%
\end{equation}
Thus, $\check{\pi}_{0,m}^{-1}\varphi_{m}\left(  t_m,\boldsymbol{\mu}\right)
\rightarrow1$ when $t_m^{-1}\left(  1+\tilde{u}_{3,m}^{-1}\right)  =o\left(
\check{\pi}_{0,m}\right)  $.
Since (\ref{eq19a}), (\ref{eq19b}) and
(\ref{eq19bx}) asymptotically are identical to (\ref{eq18c}), (\ref{eq18b})
and (\ref{eqGammaBoundedNullB}) respectively, we obtain from (\ref{eq19a}),
(\ref{eq19b}) and (\ref{eq19bx}) the claimed uniform consistency class for
$\sigma\geq 3/4$ and $\sigma\leq3/4$ respectively, for which $\gamma=1$ can be set when $\sigma>11/4$.\qed

\section{Simulation study}

\label{SecNumericalStudies}

We will present a simulation study on the proposed estimators, with a
comparison to the ``MR'' estimator of
\cite{Meinshausen:2006} or Storey's estimator of \cite{Storey:2004} for the
case of a one-sided {null}. For one-sided null $\Theta_{0}=\left(
-\infty,b\right)  \cap U$, when $X_{0}$ is an observation from a random
variable $X$ with CDF $F_{\mu}$, $\mu\in U$, its one-sided p-value is computed
as $1-F_{b}\left(  X_{0}\right)  $.
We will not include a comparison with the {two estimators of \cite{Dickhaus:2013,Hoang:2022a,Hoang:2022b}}, since it is not an aim here to investigate
for {Gamma random variables} whether the definition of randomized p-value of \cite{Dickhaus:2013,Hoang:2022a,Hoang:2022b} leads to valid randomized p-values that can be practically computed.

We numerically implement the solution $\left(\psi,K\right)$ in two cases as follows: (a) if $\psi$ or $K$ is defined by a univariate integral, then the univariate integral is approximated by a Riemann sum based on an equally spaced partition with norm $0.01$ of the corresponding domain of integration; (b) if $\psi$ or $K$ is defined by a double integral, then the double integral is computed as an iterated integral, for which each univariate integral is computed as if it were case (a).
We choose norm $0.01$ for a partition so as to reduce a bit the computational complexity
of the proposed estimators when the number of hypotheses to test is very large.
However, we will not explore here how much more accurate these estimators can
be when finer partitions are used to obtain the Riemman sums, or explore here which density function $\omega(s)$ on $[-1,1]$ should be used to
give the best performances to the proposed estimators among all continuous densities on $[-1,1]$ that are of bounded variation. By default, we will choose
the triangular density $\omega\left(s\right)=\left(1-\vert s\vert\right)1_{\left[-1,1\right]}\left(s\right)$, since numerical evidence in \cite{Jin:2008,Chen:2018a}
shows that this density leads to good performances of the proposed estimators for the setting of a point null.

The MR estimator (designed for continuous p-values) is implemented as follows: let the
ascendingly ordered p-values be $p_{\left(  1\right)  }<p_{\left(  2\right)
}<\cdots<p_{\left(  m\right)  }$ for $m>4$, set $b_{m}^{\ast}=m^{-1/2}%
\sqrt{2\ln\ln m}$, and define
\[
q_{i}^{\ast}=\left(  1-p_{\left(  i\right)  }\right)  ^{-1}\left\{
im^{-1}-p_{\left(  i\right)  }-b_{m}^{\ast}\sqrt{p_{\left(  i\right)  }\left(
1-p_{\left(  i\right)  }\right)  }\right\}  ;
\]
then $\hat{\pi}_{1,m}^{\mathsf{MR}}=\min\left\{  1,\max\left\{  0,\max_{2\leq
i\leq m-2}q_{i}^{\ast}\right\}  \right\}  $ is the MR estimator. Storey's
estimator will be implemented by the \textsf{qvalue} package (version 2.14.1)
via the \textrm{`pi0.method=smoother'} option. All simulations will be done
with \textsf{R} version 3.5.0.

For an estimator $\hat{\pi}_{1,m}$ of $\pi_{1,m}$ or an estimator $\hat{\pi
}_{0,m}$ of $\tilde{\pi}_{0,m}$, its accuracy is measured by the excess
$\tilde{\delta}_{m}=\hat{\pi}_{1,m}\pi_{1,m}^{-1}-1$ or $\tilde{\delta}%
_{m}=\hat{\pi}_{0,m}\tilde{\pi}_{0,m}^{-1}-1$. For each experiment, the mean
$\mu_{m}^{\ast}$ and standard deviation $\sigma_{m}^{\ast}$ of $\tilde{\delta
}_{m}$ is estimated from independent realizations of the experiment. Among two
estimators, the one that has smaller $\sigma_{m}^{\ast}$ is taken to be more
stable, and the one that has both smaller $\sigma_{m}^{\ast}$ and smaller
$\left\vert \mu_{m}^{\ast}\right\vert $ is better. {
In each boxplot in each figure of simulation results to be presented later, the horizontal bar has been programmed to represent the mean of the quantity being plotted and the black dots represent the outliers from the quantity being plotted.}

\subsection{Simulation design and results}

\label{simDesignG}
For $a<b$, let $\mathsf{U}\left(  a,b\right)  $ be the uniform
distribution on the closed interval $\left[  a,b\right]  $.
When implementing the estimator in \autoref{ThmConstructionMoments} or
\autoref{TmVNEF}, the power series in the definition of $K$ in (\ref{IV-b}) or
(\ref{V-b}) is replaced by the partial sum of its first $26$ terms, i.e., the
power series is truncated at $n=25$. However, the double integral in $K$ in
(\ref{IV-b}) or (\ref{V-b}) has to be approximated by a Riemann sum (using the scheme described in the beginning of \autoref{SecNumericalStudies}) for each
$z_{i}$ for a total of $m$ times. This greatly increases the computational
complexity of applying $K$ to $\left\{  z_{i}\right\}  _{i=1}^{m} $ when $m$
is very large. So, we only consider $4$ values for $m$, i.e., $m=10^{3}$,
$5\times10^{3}$, $10^{4}$ or $5 \times10^{4}$, together with $2$ sparsity
levels $\pi_{1,m}=0.2$ \textcolor{black}{(indicating the ``dense regime")} or $\left(  \ln\ln m\right)  ^{-1}$ \textcolor{black}{(indicating the ``moderately sparse regime")}. We set $\sigma=4$
for the simulated Gamma random variables. The speed of the proposed estimators
$t_{m}=\sqrt{0.25\sigma^{-1}u_{3,m}\ln m}$ (i.e., $\gamma=1$ is set for
$t_{m}$) for a bounded null and $t_{m}=2^{-5/4}\sigma^{-1/2}\sqrt{u_{3,m}\ln
m}$ (i.e., $\gamma=1$ is set for $t_{m}$) for a one-side null, both with
$u_{3,m}=0.2/\ln\ln m$, so that the consistency conditions in
\autoref{ConcentrationIII} and \autoref{ThmVNEFConsistency} are satisfied. The
simulated data are generated as follows:

\begin{itemize}
\item Scenario {I} ``estimating $\pi_{1,m}$ for a bounded
null'': set $\theta_{a}=0$, $\theta_{b}=0.35$, $\theta_{\ast
}=-0.2$ and $\theta^{\ast}=0.55$; generate $m_{0}$ $\theta_{i}$'s
independently from $\mathsf{U}\left(  \theta_{a}+u_{3,m},\theta_{b}%
-u_{3,m}\right)  $, $m_{11}$ $\theta_{i}$'s independently from $\mathsf{U}%
\left(  \theta_{b}+u_{\textcolor{black}{3,m}},\theta^{\ast}\right)  $, and $m_{11}$ $\theta_{i}$'s
independently from $\mathsf{U}\left(  \theta_{\ast},\theta_{a}-u_{\textcolor{black}{3,m}}\right)
$, where $m_{11}=\max \{  1,\textcolor{black}{\lfloor}  0.5m_{1} \textcolor{black}{\rfloor}  - \textcolor{black}{\lfloor}  m/\ln\ln
m \textcolor{black}{\rfloor}   \}  $ \textcolor{black}{and $\lfloor x \rfloor$ is the integer part of $x \in \mathbb{R}$}; set half of the remaining $m-m_{0}-2m_{11}$ $\theta
_{i}$'s to be $\theta_{a}$, and the rest to be $\theta_{b}$.

\item Scenario {II} ``estimating $\pi_{1,m}$ for a one-sided
null'': generate $m_{0}$ $\mu_{i}$'s independently from
$\mathsf{U} (  \theta_{\ast},\theta_{b}-u_{\textcolor{black}{3,m}} )  $, and $\textcolor{black}{\lfloor}
0.9m_{1} \textcolor{black}{\rfloor}  $ $\mu_{i}$'s independently from $\mathsf{U}\left(
\theta_{b}+u_{\textcolor{black}{3,m}},\theta^{\ast}\right)  $; set the rest $\theta_{i}$'s to be
$\theta_{b}$.
\end{itemize}

Each triple of $\left(  m,\pi_{1,m},\Theta_{0}\right)  $ determines an
experiment, and there are $16$ experiments in total. Each experiment is
repeated independently $100$ times. The assessment method for an estimator
$\hat{\pi}_{1,m}$ of $\pi_{1,m}$ is again based on the mean and standard
deviation of the excess $\tilde{\delta}_{m}=\hat{\pi}_{1,m}\pi_{1,m}^{-1}-1$.
\textcolor{black}{
As mentioned earlier in this section, to numerically approximate $K$ and hence the new estimators (since they are defined by integrals and power series), we computed a $26$-term partial sum of each of those power series and computed Riemann sums based on an equally spaced one-dimensional domain partition with norm $0.01$ for those integrals.
Namely, the new estimators are implemented by this scheme of numerical approximation, which we call ``numerical versions''. Even with this simple approximation scheme of relatively low computational complexity, the sequential nature of computing an approximation to $K$ and evaluating this approximation at each of the $m$ observations via ``for'' loops and the sequential nature of repeating an experiment via a ``for'' loop
took much time, and the simulations took around $45$ days to complete on a computer with $8$-core CPU and $64$GB of RAM.
Note that this numerical implementation/approximation causes numerical error and that the simulations are for the ``numerical versions'' of the new estimators rather than the new estimators themselves.}

\begin{figure}[th]
\centering
\includegraphics[height=0.65\textheight,width=0.95\textwidth]{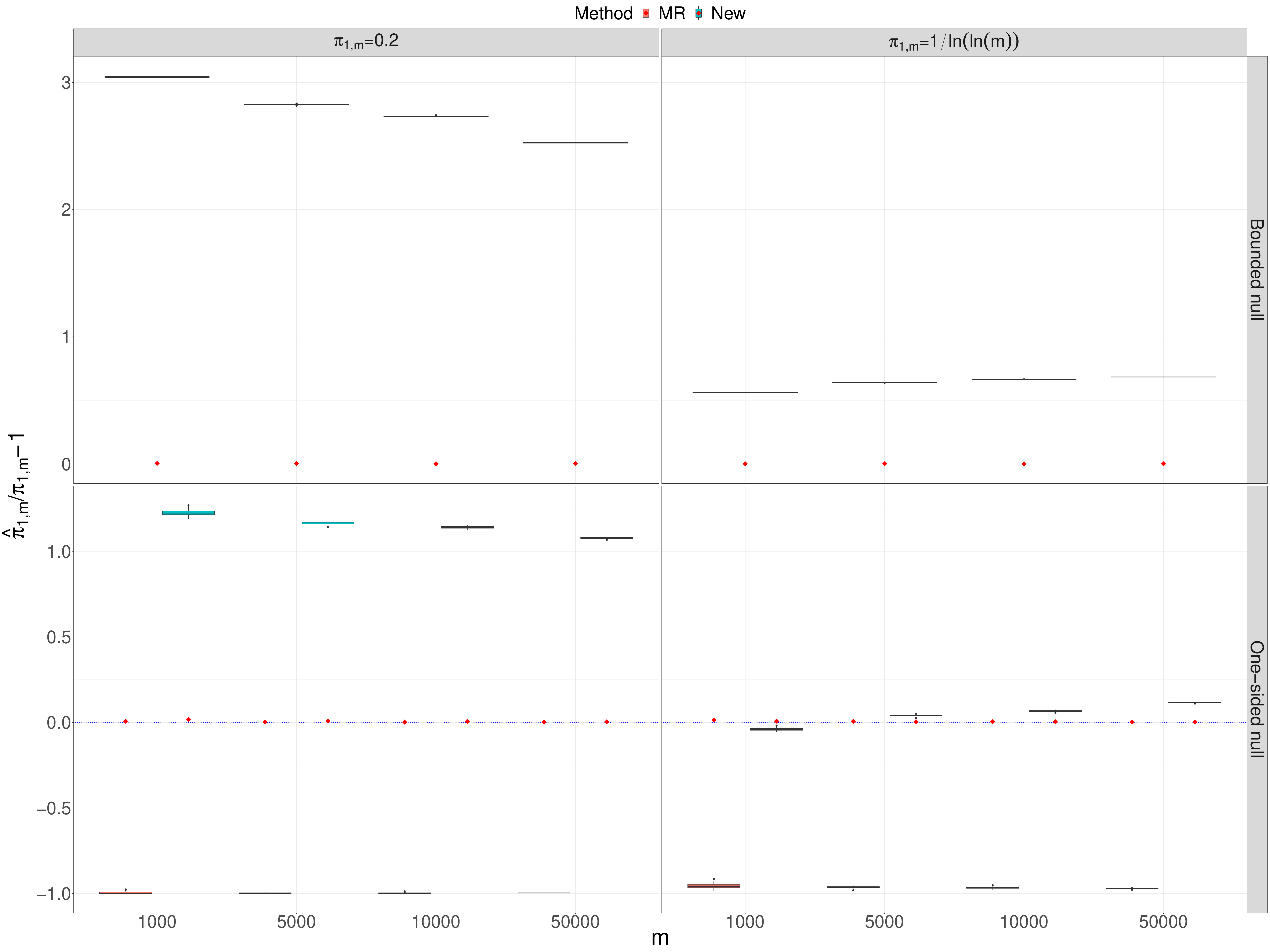}
\vspace{-0.4cm}\caption[Gaussian]{Boxplot of the excess $\tilde{\delta}_{m}$
(on the vertical axis) of an estimator $\hat{\pi}_{1,m}$ of ${\pi}_{1,m}$ {as $\tilde{\delta}_{m}=\hat{\pi}_{1,m}\pi_{1,m}^{-1}-1$}. The
thick horizontal line and the diamond in each boxplot are respectively the
mean and standard deviation of $\tilde{\delta}_{m}$, and the dotted horizontal
line is the reference for $\tilde{\delta}_{m}=0$. {An estimator with a narrower boxplot that is closer to the dotted horizontal line is better.}
All estimators have been
applied to Gamma family. For the case of a one-side null, the right one for
each pair of boxplots for each $m$ is for the proposed estimator ``New'' and
the left one is for the ``MR'' estimator. No simulation was done for the
``MR'' estimator for a bounded null. \textcolor{black}{Note that the ``Method" legend for boxplots is basically invisible
in the subplots since each boxplot contains observations that vary so little and are hence very narrow vertically.}}%
\label{fig2}%
\end{figure}

\autoref{fig2} visualizes the simulation results, for which Storey's estimator
is not shown since it is always $0$ for all experiments in Scenario {II}, \textcolor{black}{and Table 1 provides numerical summaries that complement the visualizations in \autoref{fig2}. Please note again that these results are for the numerical implementation, i.e., numerical approximation, of the new estimators, rather than the new estimators themselves, even though the interpretations of the results will be for the new estimators.}
The following three observations can be made: (i) for estimating the
alternative proportion for a one-sided null, the proposed estimator is
more accurate than the MR estimator, is very stable, and in the dense regime shows a clear trend
of convergence towards consistency. In contrast, the MR estimator is always
very close to $0$ \textcolor{black}{regardless of the sparsity regime for $\pi_{1,m}$}, either failing to detect the existence of alternative
hypotheses or very inaccurately estimating the alternative proportion. \textcolor{black}{This largely explains why for a one-side null in the dense regime, our
new estimators have slightly larger standard deviation than the MR estimator (because the latter almost always gives an estimate that is very close to $0$)}. (ii)
for estimating the alternative proportion for a bounded null, the proposed
estimator is stable and reasonably accurate, and in the dense regime shows a clear trend of
convergence towards consistency. (iii) the proposed estimator seems to be much
more accurate in the moderately sparse regime than in the dense regime.
{We remark that} the accuracy and speed of
convergence of the proposed estimators can be improved by employing more
accurate Riemann sums for the integrals and more accurate partial sums of the power series
in the computation of the matching function than currently used.
\textcolor{black}{(iv) non-asymptotically the new estimator $\hat{\varphi}_m\left(t,\boldsymbol{\mu}\right)$ of the proportion of false nulls $\pi_{1,m}$ often over-estimates
$\pi_{1,m}$, meaning that its dual $\hat{\psi}_m\left(t,\boldsymbol{\mu}\right)$, which estimates the proportion of true nulls $\pi_{0,m}=1-\pi_{1,m}$, usually
under-estimates $\pi_{0,m}$. In terms of false discovery rate (FDR) control in nonasymptotic settings, an adaptive FDR procedure that uses the new estimators $\hat{\psi}_m\left(t,\boldsymbol{\mu}\right)$ may fail to maintain a prespecified nominal FDR, even though such a procedure may have larger power compared to its non-adaptive counterparts.
}

\textcolor{black}{
Now let us explain why in the moderately sparse regime, i.e., $\pi_{1,m}=1/\ln{(\ln{m})}$, \autoref{fig2} does not provide numerical evidence that our ``New'' estimators
are consistent but does not undermine our rigorous theory, and why in the dense regime, i.e., $\pi_{1,m}=0.2$, \autoref{fig2} provides numerical evidence that our ``New'' estimators are consistent but not with a fast enough speed of convergence.
Recall that we have truncated the power series (that define $K$ and our ``New" estimators) to a $26$-term finite sum and used Riemannian sums of one-dimensional partition norm $0.01$ to approximate integrals (that define $K$ and our ``New" estimators) when implementing these estimators, which gives the actual ``New" estimators as, e.g., $\hat{\pi}^{\dagger,\textrm{New}}_{1,m}$.
Let $\hat{\pi}_{1,m}^{\textrm{New}}$ specifically denote our ``New" estimators. Then,
the numerical error of the ``New" estimators $\hat{\pi}_{1,m}^{\textrm{new}}$ is $\tilde{e}_m^{\textrm{New}} = \hat{\pi}_{1,m}^{\textrm{New}} - \hat{\pi}^{\dagger,\textrm{New}}_{1,m}$.
\\
Recall $\tilde{\delta}_m = \hat{\pi}_{1,m}/\pi_{1,m} -1$, where $\hat{\pi}_{1,m}$ is an estimate of $\pi_{1,m}$, and that $\tilde{\delta}_m$ converges to $0$ as $m \to \infty$ is equivalent to the consistency of the estimator $\hat{\pi}_{1,m}$. Due to our numerical approximation, for our ``New" estimator, $\tilde{\delta}_m$ is actually computed as
\begin{equation*}
  \tilde{\delta}_m = \frac{\hat{\pi}^{\dagger,\textrm{New}}_{1,m}}{\pi_{1,m}} -1
  = \frac{\hat{\pi}^{\textrm{New}}_{1,m}}{\pi_{1,m}} - \frac{\tilde{e}_m^{\textrm{New}}}{\pi_{1,m}}  -1.
\end{equation*}
Since our theory has rigorously proved that $\hat{\pi}^{\textrm{New}}_{1,m}/\pi_{1,m} -1$ converges to $0$ in probability as $m \to \infty$, we see the actual $\tilde{\delta}_m$ computed for our ``New" estimator $\hat{\pi}^{\textrm{New}}_{1,m}$, as given above, satisfies
\begin{equation*}
  \tilde{\delta}_m  \approx \frac{\tilde{e}_m^{\textrm{New}}}{\pi_{1,m}} \quad \text{with high probability for large } m.
\end{equation*}
However, the numerical error $\tilde{e}_m^{\textrm{New}}$ may not converge to $0$ as $m \to \infty$.
So, in the dense regime when $\pi_{1,m} = 0.2$, we will see a trend of convergence for $\tilde{\delta}_m $ as $m$ increases. But such a convergence may stall if $\tilde{e}_m^{\textrm{New}}$ does not decrease with $m$. Non-monotone decreasing or not small enough $\tilde{e}_m^{\textrm{New}}$ also creates a feeling that $\vert\tilde{\delta}_m\vert$ is larger for the ``New'' estimator than the ``MR'' estimator, which is not true for all $m$ but may be true for small $m$. This is exactly what happened
for the dense regime in \autoref{fig2}. In contrast, in the moderately sparse regime when $\pi_{1,m}=1/\ln{(\ln{m})}$, $\pi_{1,m}$ monotonically decreases as $m$ increases and $\pi_{1,m}$ converges to $0$ as $m \to \infty$. So, when $\tilde{e}_m^{\textrm{New}}$ is not of smaller order than $\pi_{1,m}=1/\ln{(\ln{m})}$ as $m$ increases, we may see the actual $\tilde{\delta}_m$ on average increases with $m$. This is exactly what happened
for the moderately sparse regime in \autoref{fig2}.\\
Unless we increase the numerical precision or equivalently reduce the numerical error $\tilde{e}_m^{\textrm{New}}$ (dynamically also with respect to $m$), increasing $m$ but keeping the current numerical approximation scheme as described earlier will not allow us to see a clear trend of convergence of our
``New" estimators in the moderately sparse regime where $\lim_{m \to \infty}\pi_{1,m}=0$. In fact, how to rigorously and precisely control the numerical error when truncating a power series and using Riemann sums to approximate integrals when implementing our proposed estimator requires very delicate analysis, may well form another manuscript, and unfortunately cannot be fully numerically explored in this work. Nevertheless, for practical applications where we do not need to repeat an experiment many times as is done in the simulations here, we recommend keeping as many terms and using as fine partitions as one's computational recourses allow when respectively truncating the power series and forming Riemann sums that numerically approximate the definitions of the new estimators.
}

\section{Estimators for closed or half-closed nulls and their consistency}
\label{SecClosedNulls}

Let us discuss how to adapt the constructions, the estimators, their
concentration inequalities, and their consistency results to estimating the
proportion $\pi_{1,m}$ when the null hypotheses are closed or half-closed
sets. For the Gamma family, the mean parameter $\mu$ is a function of the
natural parameter $\theta$, such that $\mu=\mu\left(  \theta\right)
=\sigma\left(  1-\theta\right)  ^{-1}$ for $\theta<1$, and this function is
differentiable and strictly monotone with inverse $\theta=\theta\left(
\mu\right)  =1-\sigma\mu^{-1}$. Here $\sigma>0$ is the scale parameter.
Further, $\xi\left(  \theta\right)  =\left(  1-\theta\right)  ^{-1}$,
$\xi\left(  \theta\right)  =\sigma^{-1}\mu\left(  \theta\right)  $, $\mu
_{i}=\mu\left(  \theta_{i}\right)  $ and $\theta_{\tau}$ is such that
$\mu\left(  \theta_{\tau}\right)  =\tau$ for $\tau\in\left\{  a,b\right\}  $.

Recall%
\[
\left\{
\begin{array}
[c]{l}%
\psi_{3,0}\left(  t,\theta;\theta^{\prime}\right)  ={\int_{\left[
-1,1\right]  }}\cos\left[  ts\left\{  \xi\left(  \theta^{\prime}\right)
-\xi\left(  \theta\right)  \right\}  \right]  \omega\left(  s\right)  ds\\
\tilde{\psi}_{1,0}\left(  t,\mu;\mu^{\prime}\right)  ={\int_{\left[
-1,1\right]  }}\omega\left(  s\right)  \cos\left\{  ts\sigma^{-1}\left(
\mu-\mu^{\prime}\right)  \right\}  ds
\end{array}
\right.  ,
\]
where $\tilde{\psi}_{1,0}\left(  t,\mu;\mu^{\prime}\right)$ is defined in \autoref{OracleSpeeds}.
We see that $\mu=\mu\left(  \theta\right)  =\sigma\left(  1-\theta\right)
^{-1}$ implies $\psi_{3,0}\left(  t,\theta;\theta^{\prime}\right)
=\tilde{\psi}_{1,0}\left(  t,\mu;\mu^{\prime}\right)  $. This fact will be
used in our discussion on the quantity $\varphi_{m}\left(
t,\boldsymbol{\theta}\right)  $ or its equivalent $\varphi_{m}\left(
t,\boldsymbol{\mu}\right)  $, where $\boldsymbol{\theta}=\left(  \theta
_{1},\ldots,\theta_{m}\right)  $, $\theta_{i}=\theta\left(  \mu_{i}\right)  $,
$\theta\left(  \boldsymbol{\mu}\right)  =\left(  \theta\left(  \mu_{1}\right)
,\ldots,\theta\left(  \mu_{m}\right)  \right)  $ and $\boldsymbol{\mu}=\left(
\mu_{1},\ldots,\mu_{m}\right)  $.

\subsection{The case of a bounded null}

\label{GammaClosedBoundNull}

When $\Theta_{0}=\left[  a,b\right]  $, we can just set%
\begin{equation}
\left\{
\begin{array}
[c]{l}%
K\left(  t,x\right)  =K_{1}\left(  t,x\right)  +2^{-1}\left\{  K_{3,0}\left(
t,x;\theta_{a}\right)  +K_{3,0}\left(  t,x;\theta_{b}\right)  \right\}  \\
\psi\left(  t,\theta\right)  =\psi_{1}\left(  t,\theta\right)  +2^{-1}\left\{
\psi_{3,0}\left(  t,\theta;\theta_{a}\right)  +\psi_{3,0}\left(
t,\theta;\theta_{b}\right)  \right\}
\end{array}
\right.  \label{GammaCloseBndNull}%
\end{equation}
in comparison to the construction when $\Theta_{0}=\left(  a,b\right)  $ as%
\begin{equation}
\left\{
\begin{array}
[c]{l}%
K\left(  t,x\right)  =K_{1}\left(  t,x\right)  -2^{-1}\left\{  K_{3,0}\left(
t,x;\theta_{a}\right)  +K_{3,0}\left(  t,x;\theta_{b}\right)  \right\}  \\
\psi\left(  t,\theta\right)  =\psi_{1}\left(  t,\theta\right)  -2^{-1}\left\{
\psi_{3,0}\left(  t,\theta;\theta_{a}\right)  +\psi_{3,0}\left(
t,\theta;\theta_{b}\right)  \right\}
\end{array}
\right.  .\label{GammaBndOpenNull}%
\end{equation}
The definitions of the estimator and its expectation for either $\Theta
_{0}=\left(  a,b\right)  $ or $\Theta_{0}=\left[  a,b\right]  $ remain
identical as%
\[
\hat{\varphi}_{m}\left(  t,\mathbf{z}\right)  =m^{-1}\sum_{i=1}^{m}\left\{
1-K\left(  t,z_{i}\right)  \right\}  \text{ \ and\ }\varphi_{m}\left(
t,\boldsymbol{\theta}\right)  =m^{-1}\sum_{i=1}^{m}\left\{  1-\psi\left(
t,\mu_{i}\right)  \right\}  .
\]

When $\Theta_{0}=\left(  a,b\right)  $, in the proofs for the estimator
$\hat{\varphi}_{m}\left(  t,\mathbf{z}\right)  $, we have used $e_{1,m}\left(
t\right)  :=\hat{\varphi}_{1,m}\left(  t,\mathbf{z}\right)  -\varphi
_{1,m}\left(  t,\boldsymbol{\theta}\right)  $, where%
\[
\hat{\varphi}_{1,m}\left(  t,\mathbf{z}\right)  =m^{-1}\sum_{i=1}^{m}%
K_{1}\left(  t,z_{i}\right)  \ \text{and }\varphi_{1,m}\left(
t,\boldsymbol{\theta}\right)  =\mathbb{E}\left\{  \hat{\varphi}_{1,m}\left(
t,\mathbf{z}\right)  \right\}  ,
\]
$e_{3,0,m}\left(  t,\tau\right)  :=\hat{\varphi}_{3,0,m}\left(  t,\mathbf{z}%
;\tau\right)  -$ $\varphi_{3,0,m}\left(  t,\boldsymbol{\theta};\tau\right)  $,
$\tau\in\left\{  a,b\right\}  $, where%
\[
\hat{\varphi}_{3,0,m}\left(  t,\mathbf{z};\tau\right)  =m^{-1}\sum_{i=1}%
^{m}K_{3,0}\left(  t,z_{i};\theta_{\tau}\right)  \text{ and }\varphi
_{3,0,m}\left(  t,\boldsymbol{\theta};\tau\right)  =\mathbb{E}\left\{
\hat{\varphi}_{3,0,m}\left(  t,\mathbf{z};\tau\right)  \right\}  ,
\]
$e_{m}\left(  t\right)  :=\hat{\varphi}_{m}\left(  t,\mathbf{z}\right)
-\varphi_{m}\left(  t,\boldsymbol{\theta}\right)  $ and%
\begin{equation}
e_{m}\left(  t\right)    =-e_{1,m}\left(  t\right)
+2^{-1}e_{3,0,m}\left(  t,a\right)  +2^{-1}e_{3,0,m}\left(  t,b\right)
.\label{ErrorGammaOpenBndNull}%
\end{equation}
Further, to upper bound the variance of $-e_{m}\left(  t\right)  $, which is
also the variance of $e_{m}\left(  t\right)  $, we have upper bounded the
variances of $e_{1,m}\left(  t\right)  $, $e_{3,0,m}\left(  t,a\right)  $ and
$e_{3,0,m}\left(  t,b\right)  $ individually, and then directly replaced each
variance in each summand on the right-hand side of the inequality%
\begin{equation}
\mathbb{V}\left[  e_{m}\left(  t\right)  \right]  \leq2\mathbb{V}\left\{
e_{1,m}\left(  t\right)  \right\}  +\mathbb{V}\left[  e_{3,0,m}\left(
t,a\right)  \right]  +\mathbb{V}\left[  e_{3,0,m}\left(  t,b\right)  \right]
\label{VarBndGammaBndNull}%
\end{equation}
with these individual variance upper bounds. In addition, concentration of
$\left\vert e_{m}\left(  t\right)  \right\vert $ is derived by Chebyshev's
inequality based on the upper bound for the variance of $e_{m}\left(
t\right)  $.

In the setting of the closed null $\Theta_{0}=\left[  a,b\right]  $,
(\ref{GammaCloseBndNull}) implies that (\ref{ErrorGammaOpenBndNull}) becomes%
\begin{equation}
e_{m}\left(  t\right)  =\hat{\varphi}_{m}\left(  t,\mathbf{z}\right)
-\varphi_{m}\left(  t,\boldsymbol{\theta}\right)  =-e_{1,m}\left(  t\right)
-2^{-1}e_{3,0,m}\left(  t,a\right)  -2^{-1}e_{3,0,m}\left(  t,b\right)
.\label{ErrorGammaClosedBndNull}%
\end{equation}
However, (\ref{VarBndGammaBndNull}) remains valid for $e_{m}\left(  t\right)
$ in (\ref{ErrorGammaClosedBndNull}). So, the upper bound on the variance of
$e_{m}\left(  t\right)  $ and the concentration of $\left\vert e_{m}\left(
t\right)  \right\vert $ we have derived for the setting $\Theta_{0}=\left(
a,b\right)  $ and construction (\ref{GammaBndOpenNull}) remain valid for
$e_{m}\left(  t\right)  $ and $\left\vert e_{m}\left(  t\right)  \right\vert $
for the setting $\Theta_{0}=\left[  a,b\right]  $ and construction
(\ref{GammaCloseBndNull}).

Now let us discuss $\varphi_{m}\left(  t,\boldsymbol{\mu}\right)  $, which is
equivalent to $\varphi_{m}\left(  t,\boldsymbol{\theta}\right)  $ due to the
mapping $\mu=\mu\left(  \theta\right)  =\sigma\left(  1-\theta\right)  ^{-1}$
for $\theta<1$. When $\Theta_{0}=\left(  a,b\right)  $,%
\begin{equation}
\varphi_{m}\left(  t,\boldsymbol{\mu}\right)  =1-\varphi_{1,m}\left(
t,\boldsymbol{\mu}\right)  +2^{-1}\varphi_{1,0,m}\left(  t,\boldsymbol{\mu
};a\right)  +2^{-1}\varphi_{1,0,m}\left(  t,\boldsymbol{\mu};b\right)
=\sum_{i=1}^{5}\widetilde{d}_{1,m}\label{OracleBndG}%
\end{equation}
and%
\begin{equation*}
\pi_{1,m}^{-1}{\varphi_{m}\left(  t,\boldsymbol{\mu}\right)  -1}=\pi
_{1,m}^{-1}\widetilde{d}_{1,m}-1+\pi_{1,m}^{-1}\widetilde{d}_{2,m}+\pi
_{1,m}^{-1}\widetilde{d}_{3,m}+\pi_{1,m}^{-1}\widetilde{d}_{4,m}+\pi
_{1,m}^{-1}\widetilde{d}_{5,m},
\end{equation*}
where%
\[
\left\{
\begin{array}
[c]{l}%
\widetilde{d}_{1,m}=1-m^{-1}\sum\nolimits_{\left\{  j:\mu_{j}\in\left(
a,b\right)  \right\}  }\psi_{1}\left(  t,\mu_{j}\right)  \\
\widetilde{d}_{2,m}=-m^{-1}\sum\nolimits_{\left\{  j:\mu_{j}=a\right\}  }%
\psi_{1}\left(  t,\mu_{j}\right)  +2^{-1}m^{-1}\sum\nolimits_{\left\{
j:\mu_{j}=a\right\}  }{\tilde{\psi}}_{1,0}\left(  t,\mu_{j};a\right)  \\
\widetilde{d}_{3,m}=-m^{-1}\sum\nolimits_{\left\{  j:\mu_{j}=b\right\}  }%
\psi_{1}\left(  t,\mu_{j}\right)  +2^{-1}m^{-1}\sum\nolimits_{\left\{
j:\mu_{j}=b\right\}  }{\tilde{\psi}}_{1,0}\left(  t,\mu_{j};b\right)  \\
\widetilde{d}_{4,m}=2^{-1}m^{-1}\sum\nolimits_{\left\{  j:\mu_{j}\neq
a\right\}  }{\tilde{\psi}}_{1,0}\left(  t,\mu_{j};a\right)  +2^{-1}m^{-1}%
\sum\nolimits_{\left\{  j:\mu_{j}\neq b\right\}  }{\tilde{\psi}}_{1,0}\left(
t,\mu_{j};b\right)  \\
\widetilde{d}_{5,m}=-m^{-1}\sum\nolimits_{\left\{  j:\mu_{j}<a\right\}  }%
\psi_{1}\left(  t,\mu_{j}\right)  -m^{-1}\sum\nolimits_{\left\{  j:\mu
_{j}>b\right\}  }\psi_{1}\left(  t,\mu_{j}\right)
\end{array}
\right.  .
\]
Further, when $\Theta_{0}=\left(  a,b\right)  $, to upper bound $\left\vert
\pi_{1,m}^{-1}{\varphi_{m}\left(  t,\boldsymbol{\mu}\right)  -1}\right\vert $,
we have replaced each $\left\vert \widetilde{d}_{j,m}\right\vert ,2\leq
j\leq5$ by its upper bound $\hat{d}_{j,m},2\leq j\leq5$ and replaced
$\left\vert \pi_{1,m}^{-1}\widetilde{d}_{1,m}-1\right\vert $ by its upper
bound $\hat{d}_{0,m}$ in the inequality%
\begin{equation}
\left\vert \pi_{1,m}^{-1}{\varphi_{m}\left(  t,\boldsymbol{\mu}\right)
-1}\right\vert \leq\left\vert \pi_{1,m}^{-1}\widetilde{d}_{1,m}-1\right\vert
+\pi_{1,m}^{-1}\left\vert \widetilde{d}_{2,m}\right\vert +\pi_{1,m}%
^{-1}\left\vert \widetilde{d}_{3,m}\right\vert +\pi_{1,m}^{-1}\left\vert
\widetilde{d}_{4,m}\right\vert +\pi_{1,m}^{-1}\left\vert \widetilde{d}%
_{5,m}\right\vert .\label{OracleBndUpperBoundG}%
\end{equation}

In case $\Theta_{0}=\left[  a,b\right]  $, (\ref{OracleBndG}) becomes%
\begin{align*}
\varphi_{m}\left(  t,\boldsymbol{\mu}\right)   &  =1-\varphi_{1,m}\left(
t,\boldsymbol{\mu}\right)  -2^{-1}\varphi_{1,0,m}\left(  t,\boldsymbol{\mu
};a\right)  -2^{-1}\varphi_{1,0,m}\left(  t,\boldsymbol{\mu};b\right) \\
&  =\widetilde{d}_{1,m}+\widetilde{d}_{2,m}^{\ast}+\widetilde{d}_{3,m}^{\ast
}-\widetilde{d}_{4,m}+\widetilde{d}_{5,m}%
\end{align*}
and%
\begin{equation}
\left\vert \pi_{1,m}^{-1}{\varphi_{m}\left(  t,\boldsymbol{\mu}\right)
-1}\right\vert \leq\left\vert \pi_{1,m}^{-1}\left(  \widetilde{d}%
_{1,m}+\widetilde{d}_{2,m}^{\ast}+\widetilde{d}_{3,m}^{\ast}\right)
-1\right\vert +\pi_{1,m}^{-1}\left\vert \widetilde{d}_{4,m}\right\vert
+\pi_{1,m}^{-1}\left\vert \widetilde{d}_{5,m}\right\vert
\label{OracleBndClosedG}%
\end{equation}
where%
\[
\left\{
\begin{array}
[c]{l}%
\widetilde{d}_{2,m}^{\ast}=-m^{-1}\sum\nolimits_{\left\{  j:\mu_{j}=a\right\}
}\psi_{1}\left(  t,\mu_{j}\right)  -2^{-1}m^{-1}\sum\nolimits_{\left\{
j:\mu_{j}=a\right\}  }{\tilde{\psi}}_{1,0}\left(  t,\mu_{j};a\right) \\
\widetilde{d}_{3,m}^{\ast}=-m^{-1}\sum\nolimits_{\left\{  j:\mu_{j}=b\right\}
}\psi_{1}\left(  t,\mu_{j}\right)  -2^{-1}m^{-1}\sum\nolimits_{\left\{
j:\mu_{j}=b\right\}  }{\tilde{\psi}}_{1,0}\left(  t,\mu_{j};b\right)
\end{array}
\right.  .
\]
However,
\[
\pi_{1,m}^{-1}\left\vert \widetilde{d}_{1,m}+\widetilde{d}_{2,m}^{\ast
}+\widetilde{d}_{3,m}^{\ast}-1\right\vert \leq\hat{d}_{0,m}+\pi_{1,m}^{-1}%
\hat{d}_{2,m}+\pi_{1,m}^{-1}\hat{d}_{3,m}.
\]
So, the upper bound for $\left\vert \pi_{1,m}^{-1}{\varphi_{m}\left(
t,\boldsymbol{\mu}\right)  -1}\right\vert $ in (\ref{OracleBndUpperBoundG})
when $\Theta_{0}=\left(  a,b\right)  $ is also an upper bound for $\left\vert
\pi_{1,m}^{-1}{\varphi_{m}\left(  t,\boldsymbol{\mu}\right)  -1}\right\vert $
in (\ref{OracleBndClosedG}) when $\Theta_{0}=\left[  a,b\right]  $.

Therefore, all results we have derived for the estimator $\hat{\varphi}%
_{m}\left(  t,\mathbf{z}\right)  $ when $\Theta_{0}=\left(  a,b\right)  $ for
the construction (\ref{GammaBndOpenNull}) remain valid for the estimator
$\hat{\varphi}_{m}\left(  t,\mathbf{z}\right)  $ when $\Theta_{0}=\left[
a,b\right]  $ for the construction (\ref{GammaCloseBndNull}).

\subsection{The case of a one-sided null}

When $\Theta_{0}=(-\infty,b]$, we can just set%
\begin{equation}
\left\{
\begin{array}
[c]{l}%
K\left(  t,x\right)  =2^{-1}-K_{1}\left(  t,x\right)  +2^{-1}K_{3,0}\left(
t,x;\theta_{b}\right)  \\
\psi\left(  t,\theta\right)  =2^{-1}-\psi_{1}\left(  t,\theta\right)
+2^{-1}\psi_{3,0}\left(  t,\theta;\theta_{b}\right)
\end{array}
\right.  \label{GammaOneSideClosedNull}%
\end{equation}
in comparison to the construction when $\Theta_{0}=\left(  -\infty,b\right)  $
as%
\begin{equation}
\left\{
\begin{array}
[c]{l}%
K\left(  t,x\right)  =2^{-1}-K_{1}\left(  t,x\right)  -2^{-1}K_{3,0}\left(
t,x;\theta_{b}\right)  \\
\psi\left(  t,\theta\right)  =2^{-1}-\psi_{1}\left(  t,\theta\right)
-2^{-1}\psi_{3,0}\left(  t,\theta;\theta_{b}\right)
\end{array}
\right.  .\label{GammaOneSideOpenNull}%
\end{equation}
The definitions of the estimator and its expectation for either $\Theta
_{0}=(-\infty,b]$ or $\Theta_{0}=\left(  -\infty,b\right)  $ remain identical
as%
\[
\hat{\varphi}_{m}\left(  t,\mathbf{z}\right)  =m^{-1}\sum_{i=1}^{m}\left\{
1-K\left(  t,z_{i}\right)  \right\}  \text{ \ and\ }\varphi_{m}\left(
t,\boldsymbol{\mu}\right)  =m^{-1}\sum_{i=1}^{m}\left\{  1-\psi\left(
t,\mu_{i}\right)  \right\}  .
\]

We will reuse the definitions of $e_{m}\left(  t\right)  $, $e_{1,m}\left(
t\right)  $ and $e_{3,0,m}\left(  t,b\right)  $ introduced previously in
\autoref{GammaClosedBoundNull}. Then $e_{m}\left(  t\right)  =e_{1,m}\left(
t\right)  +2^{-1}e_{3,m,0}\left(  t,b\right)  $ when $\Theta_{0}=\left(
-\infty,b\right)  $ becomes%
\[
e_{m}\left(  t\right)  =e_{1,m}\left(  t\right)  -2^{-1}e_{3,0,m}\left(
t,b\right)  \text{ when }\Theta_{0}=(-\infty,b].
\]
Again, to upper bound the variance of $-e_{m}\left(  t\right)  $, which is
also the variance of $e_{m}\left(  t\right)  $, we have upper bounded the
variances of $e_{1,m}\left(  t\right)  $ and $e_{3,0,m}\left(  t,b\right)  $
individually, and then directly replaced each variance in each summand on the
right-hand side of the inequality%
\begin{equation}
\mathbb{V}\left[  e_{m}\left(  t\right)  \right]  \leq2\mathbb{V}\left\{
e_{1,m}\left(  t\right)  \right\}  +2^{-1}\mathbb{V}\left[  e_{3,0,m}\left(
t,b\right)  \right]  \label{VarianceBndOneSideNullGamma}%
\end{equation}
with these individual variance upper bounds. In addition, concentration of
$\left\vert e_{m}\left(  t\right)  \right\vert $ is derived by Chebyshev's
inequality based on the upper bound for the variance of $e_{m}\left(
t\right)  $. However, (\ref{VarianceBndOneSideNullGamma}) remains valid for
$e_{m}\left(  t\right)  $ when $\Theta_{0}=(-\infty,b]$.
So, the upper bound on the variance of
$e_{m}\left(  t\right)  $ and the concentration of $\left\vert e_{m}\left(
t\right)  \right\vert $ we have derived for the setting $\Theta_{0}=\textcolor{black}{(
-\infty,b)}  $ remain valid for
$e_{m}\left(  t\right)  $ and $\left\vert e_{m}\left(  t\right)  \right\vert $
for the setting $\Theta_{0}=\textcolor{black}{(-\infty  ,b]}  $.

Now let us discuss $\varphi_{m}\left(  t,\boldsymbol{\mu}\right)  $, which is
equivalent to $\varphi_{m}\left(  t,\boldsymbol{\theta}\right)  $ due to the
mapping $\mu=\mu\left(  \theta\right)  =\sigma\left(  1-\theta\right)  ^{-1}$
for $\theta<1$. When $\Theta_{0}=\left(  -\infty,b\right)  $, we have%
\[
\varphi_{m}\left(  t,\boldsymbol{\mu}\right)  =2^{-1}+\varphi_{1,m}\left(
t,\boldsymbol{\mu}\right)  +2^{-1}\varphi_{1,0,m}\left(  t,\boldsymbol{\mu
};{b}\right)  =\bar{d}_{1,m}+\bar{d}_{2,m}+\bar{d}_{3,m}+\bar{d}_{4,m},
\]
where%
\[
\left\{
\begin{array}
[c]{l}%
\bar{d}_{1,m}=m^{-1}\sum\nolimits_{\left\{  i:\mu_{i}>{b}\right\}  }\left(
2^{-1}+\psi_{1}\left(  t,\mu_{i}\right)  \right)  \\
\bar{d}_{2,m}=m^{-1}\sum\nolimits_{\left\{  i:\mu_{i}={b}\right\}  }\left(
2^{-1}+\psi_{1}\left(  t,\mu_{i}\right)  +2^{-1}\tilde{\psi}_{1,0}\left(
t,\mu_{i};{b}\right)  \right)  \\
\bar{d}_{3,m}=m^{-1}\sum\nolimits_{\left\{  i:\mu_{i}<{b}\right\}  }\left(
2^{-1}+\psi_{1}\left(  t,\mu_{i}\right)  \right)  \\
\bar{d}_{4,m}=2^{-1}m^{-1}\sum\nolimits_{\left\{  i:\mu_{i}\neq{b}\right\}
}\tilde{\psi}_{1,0}\left(  t,\mu_{i};{b}\right)
\end{array}
\right.
\]
and specifically $\bar{d}_{2,m}=m^{-1}\sum\nolimits_{\left\{  i:\mu_{i}%
={b}\right\}  }1$. Further, when $\Theta_{0}=\left(  -\infty,b\right)  $, to
upper bound $\left\vert \pi_{1,m}^{-1}{\varphi_{m}\left(  t,\boldsymbol{\mu
}\right)  -1}\right\vert $, we have replaced each $\left\vert \bar{d}%
_{j,m}\right\vert ,3\leq j\leq4$ by its upper bound and replaced $\left\vert
\pi_{1,m}^{-1}\left(  \bar{d}_{1,m}+\bar{d}_{2,m}\right)  -1\right\vert $ by
its upper bound in the inequality%
\[
\left\vert \pi_{1,m}^{-1}{\varphi_{m}\left(  t,\boldsymbol{\mu}\right)
-1}\right\vert \leq\left\vert \pi_{1,m}^{-1}\left(  \bar{d}_{1,m}+\bar
{d}_{2,m}\right)  -1\right\vert +\pi_{1,m}^{-1}\left\vert \bar{d}%
_{3,m}\right\vert +\pi_{1,m}^{-1}\left\vert \bar{d}_{4,m}\right\vert ,
\]
where%
\begin{equation}
\pi_{1,m}^{-1}\left(  \bar{d}_{1,m}+\bar{d}_{2,m}\right)  -1=\pi_{1,m}%
^{-1}m^{-1}\sum\nolimits_{\left\{  i:\mu_{i}>{b}\right\}  }\left(  \psi
_{1}\left(  t,\mu_{i}\right)  -2^{-1}\right)  .\label{OracleTrickOneSideG}%
\end{equation}
Specifically, the upper bound on $\left\vert \pi_{1,m}^{-1}\left(  \bar
{d}_{1,m}+\bar{d}_{2,m}\right)  -1\right\vert $ is directly based on the
inequality $$\left\vert \psi_{1}\left(  t,\mu_{i}\right)  -2^{-1}\right\vert
\leq2\left(  t\tilde{u}_{m}\right)  ^{-1} \text{\ for \ } \mu_{i}>b,$$ where $\tilde
{u}_{m}=\min_{\left\{  j:\mu_{j}\neq b\right\}  }\left\vert \mu_{j}%
-b\right\vert $.

In contrast, when $\Theta_{0}=(-\infty,b]$, we have%
\[
\varphi_{m}\left(  t,\boldsymbol{\mu}\right)  =2^{-1}+\varphi_{1,m}\left(
t,\boldsymbol{\mu}\right)  -2^{-1}\varphi_{1,0,m}\left(  t,\boldsymbol{\mu
};{b}\right)  =\bar{d}_{1,m}+\bar{d}_{2,m}^{\ast}+\bar{d}_{3,m}-\bar{d}_{4,m},
\]
and%
\[
\left\vert \pi_{1,m}^{-1}{\varphi_{m}\left(  t,\boldsymbol{\mu}\right)
-1}\right\vert \leq\left\vert \pi_{1,m}^{-1}\bar{d}_{1,m}-1\right\vert
+\pi_{1,m}^{-1}\left\vert \bar{d}_{3,m}\right\vert +\pi_{1,m}^{-1}\left\vert
\bar{d}_{4,m}\right\vert ,
\]
where%
\[
\bar{d}_{2,m}^{\ast}=m^{-1}\sum\nolimits_{\left\{  i:\mu_{i}={b}\right\}
}\left(  2^{-1}+\psi_{1}\left(  t,\mu_{i}\right)  -2^{-1}\tilde{\psi}%
_{1,0}\left(  t,\mu_{i};{b}\right)  \right)  =0.
\]
However, again%
\[
\pi_{1,m}^{-1}\bar{d}_{1,m}-1=\pi_{1,m}^{-1}m^{-1}\sum\nolimits_{\left\{
i:\mu_{i}>{b}\right\}  }\left(  \psi_{1}\left(  t,\mu_{i}\right)
-2^{-1}\right)  ,
\]
whose right-hand side is identical for that of (\ref{OracleTrickOneSideG}).
Therefore, the upper bound we have derived for $\left\vert \pi_{1,m}%
^{-1}{\varphi_{m}\left(  t,\boldsymbol{\mu}\right)  -1}\right\vert $ when
$\Theta_{0}=\left(  -\infty,b\right)  $ is also an upper bound for $\left\vert
\pi_{1,m}^{-1}{\varphi_{m}\left(  t,\boldsymbol{\mu}\right)  -1}\right\vert $
when $\Theta_{0}=(-\infty,b]$.

Therefore, results we have derived for the\ estimator $\hat{\varphi}%
_{m}\left(  t,\mathbf{z}\right)  $ when $\Theta_{0}=(-\infty,b)$ for the
construction (\ref{GammaOneSideOpenNull}) remain valid for the estimator
$\hat{\varphi}_{m}\left(  t,\mathbf{z}\right)  $ when $\Theta_{0}=(-\infty,b]$
for the construction (\ref{GammaOneSideClosedNull}).

\subsection{The case of the extensions}

When $\Theta_{0}=\left[  a,b\right]  $, we can just set%
\begin{equation}
\left\{
\begin{array}
[c]{l}%
K\left(  t,x\right)  =K_{1}\left(  t,x\right)  +2^{-1}\left\{  \phi\left(
a\right)  K_{3,0}\left(  t,x;a\right)  +\phi\left(  b\right)  K_{3,0}\left(
t,x;b\right)  \right\}  \\
\psi\left(  t,\mu\right)  =\psi_{1}\left(  t,\mu\right)  +2^{-1}\left\{
\phi\left(  a\right)  \psi_{3,0}\left(  t,\mu;a\right)  +\phi\left(  b\right)
\psi_{3,0}\left(  t,\mu;b\right)  \right\}
\end{array}
\right.  \label{GammaExtClosedNull}%
\end{equation}
in comparison to the construction for $\Theta_{0}=\left(  a,b\right)  $ as%
\begin{equation}
\left\{
\begin{array}
[c]{l}%
K\left(  t,x\right)  =K_{1}\left(  t,x\right)  -2^{-1}\left\{  \phi\left(
a\right)  K_{3,0}\left(  t,x;a\right)  +\phi\left(  b\right)  K_{3,0}\left(
t,x;b\right)  \right\}  \\
\psi\left(  t,\mu\right)  =\psi_{1}\left(  t,\mu\right)  -2^{-1}\left\{
\phi\left(  a\right)  \psi_{3,0}\left(  t,\mu;a\right)  +\phi\left(  b\right)
\psi_{3,0}\left(  t,\mu;b\right)  \right\}
\end{array}
\right.  .\label{GammaExtOpenNull}%
\end{equation}
Again the definitions of the estimator and its expectation for either
$\Theta_{0}=\left(  a,b\right)  $ or $\Theta_{0}=\left[  a,b\right]  $ remain
identical as%
\[
\hat{\varphi}_{m}\left(  t,\mathbf{z}\right)  =m^{-1}\sum_{i=1}^{m}
K\left(  t,z_{i}\right)   \text{ \ and\ }\varphi_{m}\left(
t,\boldsymbol{\mu}\right)  =m^{-1}\sum_{i=1}^{m} \psi\left(
t,\mu_{i}\right)   ,
\]
and we will reuse the definitions of $e_{m}\left(  t\right)  $, $e_{1,m}%
\left(  t\right)  $ and $e_{3,0,m}\left(  t,b\right)  $ introduced previously
in \autoref{GammaClosedBoundNull}.

For $\Theta_{0}=\left(  a,b\right)  $ we have%
\[
e_{m}\left(  t\right)  =\hat{\varphi}_{m}\left(  t,\mathbf{z}\right)
-\varphi_{m}\left(  t,\boldsymbol{\mu}\right)  =e_{1,m}\left(  t\right)
-2^{-1}\phi\left(  a\right)  e_{3,0,m}\left(  t,a\right)  -2^{-1}\phi\left(
a\right)  e_{3,0,m}\left(  t,b\right)  ,
\]
whereas for $\Theta_{0}=\left[  a,b\right]  $ we have%
\[
e_{m}\left(  t\right)  =\hat{\varphi}_{m}\left(  t,\mathbf{z}\right)
-\varphi_{m}\left(  t,\boldsymbol{\mu}\right)  =e_{1,m}\left(  t\right)
+2^{-1}\phi\left(  a\right)  e_{3,0,m}\left(  t,a\right)  +2^{-1}\phi\left(
a\right)  e_{3,0,m}\left(  t,b\right)  .
\]
Since we have employed the inequality%
\begin{equation}
\mathbb{V}\left[  e_{m}\left(  t\right)  \right]  \leq2\mathbb{V}\left\{
e_{1,m}\left(  t\right)  \right\}  +\left\Vert \omega\right\Vert _{\infty}%
^{2}\mathbb{V}\left[  e_{1,0,m}\left(  t,a\right)  \right]  +\left\Vert
\omega\right\Vert _{\infty}^{2}\mathbb{V}\left[  e_{1,0,m}\left(  t,b\right)
\right]  \label{GammaExtVarBound}%
\end{equation}
and then directly replaced each variance in each summand on the right-hand
side of (\ref{GammaExtVarBound}) by their individual upper bounds, the upper
bound we have obtained on for the variance of $e_{m}\left(  t\right)  $ when
$\Theta_{0}=\left(  a,b\right)  $ is also an upper bound for the variance of
$e_{m}\left(  t\right)  $ when $\Theta_{0}=\left[  a,b\right]  $. Further,
since concentration inequalities for $e_{m}\left(  t\right)  $ when
$\Theta_{0}=\left(  a,b\right)  $ have been derived by Chebyshev's inequality
based on the upper bound for the variance of $e_{m}\left(  t\right)  $ when
$\Theta_{0}=\left(  a,b\right)  $, these concentration inequalities are also
valid for $e_{m}\left(  t\right)  $ when $\Theta_{0}=\left[  a,b\right]  $.

Now let us discuss $\varphi_{m}\left(  t,\boldsymbol{\mu}\right)  $, which is
equivalent to $\varphi_{m}\left(  t,\boldsymbol{\theta}\right)  $ due to the
mapping $\mu=\mu\left(  \theta\right)  =\sigma\left(  1-\theta\right)  ^{-1}$
for $\theta<1$. When $\Theta_{0}=\left(  a,b\right)  $, we have%
\[
\varphi_{m}\left(  t,\boldsymbol{\mu}\right)  =m^{-1}\sum_{i=1}^{m}\left[
\psi_{1}\left(  t,\mu_{i}\right)  -2^{-1}\left\{  \phi\left(  a\right)
\psi_{1,0}\left(  t,\mu;a\right)  +\phi\left(  b\right)  {\tilde{\psi}}%
_{1,0}\left(  t,\mu;b\right)  \right\}  \right]  ,
\]
$\varphi_{m}\left(  t,\boldsymbol{\mu}\right)  =\sum_{j=1}^{5}d_{\phi
,j}\left(  t,\boldsymbol{\mu}\right)  $ and%
\begin{equation}
\left\vert \check{\pi}_{0,m}^{-1}\varphi_{m}\left(  t,\boldsymbol{\mu}\right)
-1\right\vert \leq\left\vert \check{\pi}_{0,m}^{-1}d_{\phi,1}\left(
t,\boldsymbol{\mu}\right)  -1\right\vert +\sum_{j=2}^{5}\check{\pi}_{0,m}%
^{-1}\left\vert d_{\phi,j}\left(  t,\boldsymbol{\mu}\right)  \right\vert
,\label{OracleExtensionInequG}%
\end{equation}
where%
\[
\left\{
\begin{array}
[c]{l}%
d_{\phi,1}\left(  t,\boldsymbol{\mu}\right)  =m^{-1}\sum\nolimits_{\left\{
i:\mu_{i}\in\left(  a,b\right)  \right\}  }\psi_{1}\left(  t,\mu_{i}\right)
\\
d_{\phi,2}\left(  t,\boldsymbol{\mu}\right)  =m^{-1}\sum\nolimits_{\left\{
i:\mu_{i}=a\right\}  }\left(  \psi_{1}\left(  t,\mu_{i}\right)  -2^{-1}%
\phi\left(  a\right)  {\tilde{\psi}}_{1,0}\left(  t,\mu_{i};a\right)  \right)
\\
d_{\phi,3}\left(  t,\boldsymbol{\mu}\right)  =m^{-1}\sum\nolimits_{\left\{
i:\mu_{i}=b\right\}  }\left(  \psi_{1}\left(  t,\mu_{i}\right)  -2^{-1}%
\phi\left(  b\right)  {\tilde{\psi}}_{1,0}\left(  t,\mu_{i};b\right)  \right)
\\
d_{\phi,4}\left(  t,\boldsymbol{\mu}\right)  =-m^{-1}\left(  \sum_{\left\{
i:\mu_{i}\neq a\right\}  }+\sum_{\left\{  i:\mu_{i}\neq b\right\}  }\right)
\left\{  2^{-1}\left[  \phi\left(  a\right)  {\tilde{\psi}}_{1,0}\left(
t,\mu_{i};a\right)  +\phi\left(  b\right)  {\tilde{\psi}}_{1,0}\left(
t,\mu_{i};b\right)  \right]  \right\}  \\
d_{\phi,5}\left(  t,\boldsymbol{\mu}\right)  =m^{-1}\sum_{\left\{  i:\mu
_{i}<b\right\}  }\psi_{1}\left(  t,\mu_{i}\right)  +m^{-1}\sum_{\left\{
i:\mu_{i}>b\right\}  }\psi_{1}\left(  t,\mu_{i}\right)
\end{array}
\right.  .
\]
Further, when $\Theta_{0}=\left(  a,b\right)  $, to upper bound $\left\vert
\check{\pi}_{0,m}^{-1}\varphi_{m}\left(  t,\boldsymbol{\mu}\right)
-1\right\vert $, we have replaced each $\left\vert d_{\phi,j}\left(
t,\boldsymbol{\mu}\right)  \right\vert ,2\leq j\leq5$ by its upper bound
$\hat{d}_{\phi,j}\left(  t,\boldsymbol{\mu}\right)  ,2\leq j\leq5$ and
replaced $\left\vert \check{\pi}_{0,m}^{-1}d_{\phi,1}\left(
t,\boldsymbol{\mu}\right)  -1\right\vert $ by its upper bound $\hat{d}%
_{\phi,0}\left(  t,\boldsymbol{\mu}\right)  $ directly in
(\ref{OracleExtensionInequG}).

When $\Theta_{0}=\left[  a,b\right]  $, we have%
\[
\varphi_{m}\left(  t,\boldsymbol{\mu}\right)  =m^{-1}\sum_{i=1}^{m}\left[
\psi_{1}\left(  t,\mu_{i}\right)  +2^{-1}\left\{  \phi\left(  a\right)
\psi_{1,0}\left(  t,\mu;a\right)  +\phi\left(  b\right)  {\tilde{\psi}}%
_{1,0}\left(  t,\mu;b\right)  \right\}  \right]
\]
and%
\[
\varphi_{m}\left(  t,\boldsymbol{\mu}\right)  =d_{\phi,1}\left(
t,\boldsymbol{\mu}\right)  +d_{\phi,2}^{\ast}\left(  t,\boldsymbol{\mu
}\right)  +d_{\phi,3}^{\ast}\left(  t,\boldsymbol{\mu}\right)  -d_{\phi
,4}\left(  t,\boldsymbol{\mu}\right)  +d_{\phi,5}\left(  t,\boldsymbol{\mu
}\right)
\]
where%
\[
\left\{
\begin{array}
[c]{l}%
d_{\phi,2}^{\ast}\left(  t,\boldsymbol{\mu}\right)  =m^{-1}\sum
\nolimits_{\left\{  i:\mu_{i}=a\right\}  }\left(  \psi_{1}\left(  t,\mu
_{i}\right)  +2^{-1}\phi\left(  a\right)  {\tilde{\psi}}_{1,0}\left(
t,\mu_{i};a\right)  \right)  \\
d_{\phi,3}^{\ast}\left(  t,\boldsymbol{\mu}\right)  =m^{-1}\sum
\nolimits_{\left\{  i:\mu_{i}=b\right\}  }\left(  \psi_{1}\left(  t,\mu
_{i}\right)  +2^{-1}\phi\left(  b\right)  {\tilde{\psi}}_{1,0}\left(
t,\mu_{i};b\right)  \right)
\end{array}
\right.  .
\]
Then%
\begin{align*}
\left\vert \check{\pi}_{0,m}^{-1}\varphi_{m}\left(  t,\boldsymbol{\mu}\right)
-1\right\vert  &  \leq\left\vert \check{\pi}_{0,m}^{-1}\left[  d_{\phi
,1}\left(  t,\boldsymbol{\mu}\right)  +d_{\phi,2}^{\ast}\left(
t,\boldsymbol{\mu}\right)  +d_{\phi,3}^{\ast}\left(  t,\boldsymbol{\mu
}\right)  \right]  -1\right\vert \\
&  \text{ \ \ }+\check{\pi}_{0,m}^{-1}\left\vert d_{\phi,4}\left(
t,\boldsymbol{\mu}\right)  \right\vert +\check{\pi}_{0,m}^{-1}\left\vert
d_{\phi,5}\left(  t,\boldsymbol{\mu}\right)  \right\vert .
\end{align*}
However,%
\[
\left\vert \check{\pi}_{0,m}^{-1}\left[  d_{\phi,1}\left(  t,\boldsymbol{\mu
}\right)  +d_{\phi,2}^{\ast}\left(  t,\boldsymbol{\mu}\right)  +d_{\phi
,3}^{\ast}\left(  t,\boldsymbol{\mu}\right)  \right]  -1\right\vert \leq
\hat{d}_{\phi,0}\left(  t,\boldsymbol{\mu}\right)  +\check{\pi}_{0,m}^{-1}%
\hat{d}_{\phi,2}\left(  t,\boldsymbol{\mu}\right)  +\check{\pi}_{0,m}^{-1}%
\hat{d}_{\phi,3}\left(  t,\boldsymbol{\mu}\right)  .
\]
Therefore, the upper bound we have derived for $\left\vert \check{\pi}%
_{0,m}^{-1}\varphi_{m}\left(  t,\boldsymbol{\mu}\right)  -1\right\vert $ when
$\Theta_{0}=\left(  a,b\right)  $ is also an upper bound for $\left\vert
\check{\pi}_{0,m}^{-1}\varphi_{m}\left(  t,\boldsymbol{\mu}\right)
-1\right\vert $ when $\Theta_{0}=\left[  a,b\right]  $.

In summary, results we have derived for the estimator $\hat{\varphi}%
_{m}\left(  t,\mathbf{z}\right)  $ when $\Theta_{0}=\left(  a,b\right)  $ for
the construction (\ref{GammaExtOpenNull}) remain valid for the estimator
$\hat{\varphi}_{m}\left(  t,\mathbf{z}\right)  $ when $\Theta_{0}=\left[
a,b\right]  $ for the construction (\ref{GammaExtClosedNull}).

\bibliographystyle{dcu}



\begin{table}[ht]
\centering
\begin{tabular}{llllrr}
  \hline
  \hline
  $m$ & Method & Sparsity & Null Type & $\hat{E}(\tilde{\delta}_m)$  & $\hat{\sigma}(\tilde{\delta}_m)$ \\
  \hline
  \hline
  1000 & MR & $\pi_{1,m}=0.2$ & One-sided null & -0.9954 & 0.0061 \\
   1000 & New & $\pi_{1,m}=0.2$ & One-sided null & 1.2245 & 0.0161 \\
 \hline
  5000 & MR & $\pi_{1,m}=0.2$ & One-sided null & -0.9973 & 0.0024 \\
   5000 & New & $\pi_{1,m}=0.2$ & One-sided null & 1.1659 & 0.0090 \\
  \hline
   10000 & MR & $\pi_{1,m}=0.2$ & One-sided null & -0.9969 & 0.0022 \\
   10000 & New & $\pi_{1,m}=0.2$ & One-sided null & 1.1397 & 0.0067 \\
\hline
    50000 & MR & $\pi_{1,m}=0.2$ & One-sided null & -0.9964 & 0.0012 \\
    50000 & New & $\pi_{1,m}=0.2$ & One-sided null & 1.0782 & 0.0042 \\

   \hline
   \hline
   1000 & MR & $\pi_{1,m}=1/\ln{(\ln{m})}$ & One-sided null & -0.9558 & 0.0139 \\
   1000 & New & $\pi_{1,m}=1/\ln{(\ln{m})}$ & One-sided null & -0.0399 & 0.0077 \\
\hline
   5000 & MR & $\pi_{1,m}=1/\ln{(\ln{m})}$ & One-sided null & -0.9638 & 0.0063 \\
   5000 & New & $\pi_{1,m}=1/\ln{(\ln{m})}$ & One-sided null & 0.0393 & 0.0046 \\
\hline
   10000 & MR & $\pi_{1,m}=1/\ln{(\ln{m})}$ & One-sided null & -0.9666 & 0.0050 \\
   10000 & New & $\pi_{1,m}=1/\ln{(\ln{m})}$ & One-sided null & 0.0669 & 0.0036 \\
\hline
      50000 & MR & $\pi_{1,m}=1/\ln{(\ln{m})}$ & One-sided null & -0.9721 & 0.0021 \\
   50000 & New & $\pi_{1,m}=1/\ln{(\ln{m})}$ & One-sided null & 0.1163 & 0.0023 \\

   \hline
   \hline

   1000 & New & $\pi_{1,m}=0.2$ & Bounded null & 3.0421 & 0.0045 \\
   5000 & New & $\pi_{1,m}=0.2$ & Bounded null & 2.8260 & 0.0034 \\
   10000 & New & $\pi_{1,m}=0.2$ & Bounded null & 2.7343 & 0.0028 \\
   50000 & New & $\pi_{1,m}=0.2$ & Bounded null & 2.5244 & 0.0018 \\

   \hline
   \hline
   1000 & New & $\pi_{1,m}=1/\ln{(\ln{m})}$ & Bounded null & 0.5625 & 0.0019 \\
   5000 & New & $\pi_{1,m}=1/\ln{(\ln{m})}$ & Bounded null & 0.6416 & 0.0016 \\
   10000 & New & $\pi_{1,m}=1/\ln{(\ln{m})}$ & Bounded null & 0.6620 & 0.0014 \\
   50000 & New & $\pi_{1,m}=1/\ln{(\ln{m})}$ & Bounded null & 0.6842 & 0.0010 \\
   \hline
   \hline
\end{tabular}
\label{TabSim}
\caption{In the table, $\tilde{\delta}_m = \hat{\pi}_{1,m}/\pi_{1,m} -1$ (where $\hat{\pi}_{1,m}$ is an estimate of $\pi_{1,m}$),
 $\hat{E} (\tilde{\delta}_m )$ is the sample mean of $\tilde{\delta}_m$, and
$\hat{\sigma} (\tilde{\delta}_m )$ the sample standard deviation of $\tilde{\delta}_m$. When $\pi_{1,m}=0.2$, our proposed estimators
``New" show a clear trend of convergence to $0$ as $m$ increases. For $\pi_{1,m}=1/\ln{(\ln{m})}$ though, $\tilde{\delta}_m$ for our ``New'' estimators does not
show a clear trend of convergence to $0$ as $m$ increases. However, this is an artifact of the numerical error when implementing our ``New'' estimators,
as explained in \autoref{simDesignG}.}
\end{table}

\end{document}